\documentclass{amsart}
\usepackage{bm}
\usepackage{mathrsfs}
\usepackage{amsmath}
\usepackage{amsfonts}
\usepackage{mathtools}
\usepackage{amssymb}
\usepackage{color}
\usepackage{graphicx}
\usepackage{hyperref}
\usepackage{cleveref}
\usepackage{tikz}
\usepackage{float}

\usepackage{enumitem}
\setlist[enumerate,1]{label=\textup{(\arabic*)}}

\hypersetup{
	colorlinks=true,
	linkcolor=blue,
	urlcolor=magenta,
	citecolor=green,
}

 \newtheorem{Theorem}{Theorem}[section]
 \newtheorem{Corollary}[Theorem]{Corollary}
 \newtheorem{Lemma}[Theorem]{Lemma}
 \newtheorem{Proposition}[Theorem]{Proposition}

 \newtheorem{Definition}[Theorem]{Definition}
\newtheorem{Question}[Theorem]{Question}
 
 \newtheorem{Problem}[Theorem]{Problem}

 \newtheorem{Remark}[Theorem]{Remark}

 \newtheorem{Example}[Theorem]{Example}
 \numberwithin{equation}{section}

\DeclareMathOperator{\Tn}{Tn}
\DeclareMathOperator{\lct}{lct}

\DeclareMathOperator{\Val}{Val}
\DeclareMathOperator{\ZVal}{ZVal}
\DeclareMathOperator{\QM}{QM}

\def\i{\mathrm{i}}
\def\V{\mathcal{V}}
\def\qm{\mathrm{qm}}

\begin{document}

\title[The existence of valuative interpolation]
 {The existence of valuative interpolation}

\author{Shijie Bao}
\address{Shijie Bao: Academy of Mathematics
	and Systems Science, Chinese Academy of Sciences, Beijing 100190, China.}
\email{bsjie@amss.ac.cn}

\author{Qi'an Guan}
\address{Qi'an Guan: School of Mathematical Sciences,
	Peking University, Beijing, 100871, China.}
\email{guanqian@math.pku.edu.cn}

\author{Zhitong Mi}
\address{Zhitong Mi: School of Mathematics and Statistics, Beijing Jiaotong University, Beijing,
	100044, China.
}
\email{zhitongmi@amss.ac.cn}

\author{Zheng Yuan}
\address{Zheng Yuan: State Key Laboratory of Mathematical Sciences, Academy of Mathematics
	and Systems Science, Chinese Academy of Sciences, Beijing 100190, China.}
\email{yuanzheng@amss.ac.cn}




\date{\today}



\begin{abstract}
In this article, using key tools including Zhou valuations, Tian functions and a convergence result for relative types, we establish necessary and sufficient conditions for the existence of valuative interpolations on the rings of germs of holomorphic functions and real analytic functions at the origin in $\mathbb{C}^{n}$ and $\mathbb{R}^{n}$, respectively. For the cases of polynomial rings with complex and real coefficients, we establish separate necessary conditions and sufficient conditions, which become both necessary and sufficient when the intersection of the zero sets of the given polynomials is the set of the origin in $\mathbb{C}^{n}$.

Furthermore, we obtain a necessary and sufficient condition for a valuation to be of the form given by a relative type with respect to a tame maximal weight. We demonstrate a result of Boucksom--Favre--Jonsson on quasimonomial valuations also holds for quasimonomial Zhou valuations. Finally, we obtain a relationship between Zhou valuations and the differentiable points of Tian functions.
\end{abstract}

\maketitle

\tableofcontents

\section{Introduction}
Denote the ring of germs of holomorphic functions at the origin $o$ in $\mathbb{C}^n$ by $\mathcal{O}_o$, and $\mathcal{O}_{o}^*\coloneqq \mathcal{O}_{o}\setminus\{0\}$.
A \emph{valuation} on $\mathcal{O}_o$ is a non-constant map $\nu\colon\mathcal{O}_o^*\rightarrow\mathbb{R}_{\ge0}$ satisfying the following:

\begin{enumerate}
    \item $\nu(f g)=\nu(f)+\nu(g)$;
    \item $\nu(f+g)\ge\min\{\nu(f),\nu(g)\}$;
    \item $\nu(c)=0$, where $c\not=0$ is a constant function.
\end{enumerate}
 We set $\nu(0)=+\infty$ for any valuation $\nu$. We call the valuation $\nu$ is \emph{centered} at $o$ if $\mathcal{O}_{\nu}\coloneqq \{(f,o)\in\mathcal{O}_o^*\colon \nu(f)>0\}$ coincides with the maximal ideal $\mathfrak{m}$ of $\mathcal{O}_o$. Similarly, we can define the valuations on the polynomial ring $\mathbb{C}[z_1,\ldots,z_n]$.

A natural problem is how to determine the existence of valuative interpolation:

\begin{Question}
	\label{q:1}
Given any positive integer $m$, any finite elements $\{f_j\}_{1\le j\le m}$ in $\mathcal{O}_o$ (or in $\mathbb{C}[z_1,\ldots,z_n]$), and any finite positive real numbers $\{a_j\}_{1\le j\le m}$, can one find a valuation $\nu$ on $\mathcal{O}_o$ (or on $\mathbb{C}[z_1,\ldots,z_n]$) such that $\nu(f_j)=a_j$ for all $j\in\{1,\ldots, m\}$?
\end{Question}

There are many examples where valuative interpolation fails to exist. For instance, there is no valuation $v$ such that $\nu(z_1)=1$, $\nu(z_2)=2$, and $\nu(z_1+z_2)=3$.

The valuation theory is closely related to the singularity theory in several complex variables and algebraic geometry (see e.g. \cite{FJ04,FJ05,FJ05b,BFJ08,JON-Mus2012,BGMY-valuation}).
The concept of \emph{Lelong number} is an analytic generalization of the algebraic notion of multiplicity (see \cite{lelong57,demailly2010}), which is  defined to be $$\nu(u,o):=\sup\big\{c\ge0:u(z)\le c\log|z|+O(1)\text{ near }o\big\},$$
where $u$ is a plurisubharmonic function near $o$.
Replacing $\log|z|$ by $\log\max_{1\le j\le n}|z_j|^{\frac{1}{a_j}}$ ($a=(a_1,\ldots,a_n)\in\mathbb{R}_{>0}^n$) in the definition of Lelong numbers, the directional Lelong numbers $\nu_{a}(u,o)$ (also known as \emph{Kiselman numbers} \cite{kiselman}) were introduced by Kiselman. Both the Lelong number and Kiselman number are used to measure the singularity of $u$ at $o$ (see \cite{demailly-book}).
It is clear that Kiselman numbers induce a class of important valuations $\nu_{a}(\log|\cdot|,o)$ on $\mathcal{O}_o$ (see \cite{FJ05b, Rash06}).

A plurisubharmonic function $\varphi$ on a neighborhood of $o$ is called a \emph{maximal weight}, if $u(o)=-\infty$, $u$ is locally bounded on  $U\setminus\{o\}$ and $(dd^c \varphi)^n=0$ on $U\setminus\{o\}$ (see \cite{Rash06,BFJ08}), where $U$ is a neighborhood of $o$.

Motivated by the definition of Lelong numbers, Rashkovskii \cite{Rash06} introduced the \emph{relative type}
\begin{equation}
    \label{eq:0917e}
    \sigma(u,\varphi) \coloneqq \sup\big\{c \ge 0 \colon u \le c\varphi + O(1) \text{ near } o\big\}
\end{equation}
for a plurisubharmonic function $u$ with respect to a maximal weight $\varphi$, which also generalizes the Lelong numbers. Nevertheless, restricted to $\mathcal{O}_o$, the relative type $\sigma(\log|\cdot|,\varphi)$ is not a valuation in general. In this paper, we simply use the notation of \eqref{eq:0917e} and the name "relative type" for any plurisubharmonic function $\varphi$ near $o$, even if it is not maximal.

In \cite{BFJ08}, Boucksom--Favre--Jonsson established that every quasimonomial valuation on $\mathcal{O}_o$ can be expressed as a relative type with respect to a tame maximal weight. Furthermore, they established a fundamental result connecting Lelong numbers, multiplier ideal sheaves, relative types, and Demailly's generalized Lelong numbers. In \cite{BGMY-valuation}, we introduced a class of tame maximal weights called Zhou weights and proved that the relative types with respect to them yield valuations (called Zhou valuations) on $\mathcal{O}_o$. We then used these valuations to 
measure the singularities of plurisubharmonic functions.

In this article, we study the valuative interpolation problem on $\mathcal{O}_o$ (see \Cref{q:1}) by employing Zhou valuations, Tian functions, and several results on relative types. Here, a Tian function is defined as the function of jumping numbers with respect to the exponents of a holomorphic function or the multiples of a plurisubharmonic function. 

Recall that Kiselman numbers (also known as directional Lelong numbers) induce a class of valuations. Both Zhou valuations (see \cite{BGMY-valuation}) and quasimonomial valuations (see \cite{BFJ08}) can be expressed as relative types with respect to certain tame maximal weights (for the definition of "tame", see \Cref{sec:1.2}). We discuss the relationship between valuations and relative types, and obtain a necessary and sufficient condition for a valuation to be represented by a relative type with respect to a tame maximal weight.

Finally, we show that a result of Boucksom--Favre--Jonsson on quasimonomial valuations \cite{BFJ08} also holds for quasimonomial Zhou valuations, and we establish a relationship between Zhou valuations and the differentiable points of Tian functions.

\subsection{Valuative interpolation}

\

The following theorem establishes a criterion for the valuative interpolation problem on $\mathcal{O}_o$.

\begin{Theorem}\label{thm:interpolation}
Given $m+1$ holomorphic functions $\{f_j\}_{0\le j\le m}$ near $o$, $a_0=0$ and $m$  positive numbers $\{a_j\}_{1\le j\le m}$,
	the following two statements are equivalent:
	\begin{enumerate}
	    \item There exists a valuation $\nu$  such that $\nu(f_j)=a_j$ for every $0\le j\le m$;
	    \item $\sigma(\log|F|,\varphi)=\sum_{1\le j\le m}a_j$, where \[F\coloneqq \prod_{0\le j\le m}f_j, \quad \varphi\coloneqq \log\Big(\sum_{1\le j\le m}|f_j|^{\frac{1}{a_j}}\Big).\]
	\end{enumerate}
	
Furthermore, if $o$ is an isolated point in $\cap_{1\le j\le m}\{f_j=0\}$ and $\sigma(\log|F|,\varphi)=\sum_{1\le j\le m}a_j$, then there exists a valuation $\nu$ and a tame maximal weight $\varphi_{\nu}$ such that $\nu(f_j)=a_j$ for every $j$ and $\nu(f)=\sigma(\log|f|,\varphi_{\nu})$ for any $(f,o)\in\mathcal{O}_o$.
\end{Theorem}

Note that $\nu(f_1f_2)=0$ $\Leftrightarrow$ $\nu(f_1)=0$ and $\nu(f_2)=0$. Therefore, in the theorem above, it suffices to consider the case where exactly one function $f_0$ satisfies $\nu(f_0)=0$. Since $\nu(1)=0$ for any valuation $\nu$, in the above theorem, we can actually ignore $f_0$ when $f_0\equiv1$. This convention allows the case where all interpolated values are nonzero to be naturally incorporated into our framework. 

We give some corollaries of Theorem \ref{thm:interpolation}.

\begin{Corollary}
	\label{c:interpolation-n functions}
		Let $\{f_j\}_{1\le j\le n}$ be $n$ holomorphic functions near $o$ with $\cap_{1\le j\le n}\{f_j=0\}=\{o\}$. Then for any positive numbers $\{a_j\}_{1\le j\le n}$, there exists a valuation $\nu$ and a tame weight $\varphi_{\nu}$ such that $\nu(f_j)=a_j$ for every $j$ and $\nu(f)=\sigma(\log|f|,\varphi_{\nu})$ for any $(f,o)\in\mathcal{O}_o$.
\end{Corollary}

\begin{Corollary}
	\label{c:valuation-tame valuation}
Let $\nu$ be a valuation on $\mathcal{O}_o$ centered at $o$. Then for any finite collection of holomorphic functions $\{f_1,\ldots,f_m\}$ near $o$ with $f_j(o)=0$ for all $j$, there exists a valuation $\tilde\nu$ and a tame maximal weight $\varphi_{\tilde\nu}$ such that $\tilde\nu(f)=\sigma(\log|f|,\varphi_{\tilde\nu})$ for any $(f,o)\in\mathcal{O}_o$, and
\[\nu(f_j)=\tilde\nu(f_j), \quad j=1,\ldots,m.\]
\end{Corollary}

Let us consider the valuative interpolation problem on the polynomial ring $\mathbb{C}[z_1,\ldots,z_n]$. Let $\{f_j\}_{0\le j\le m}\subset \mathbb{C}[z_1,\ldots,z_n]$, $a_0=0$ and  $\{a_j\}_{1\le j\le m}$ be $m$ positive numbers. Denote $F\coloneqq \prod_{0\le j\le m}f_j$ and $\varphi\coloneqq \log(\sum_{1\le j\le m}|f_j|^{\frac{1}{a_j}})$.

\begin{Corollary}\label{C:interpolation-comp-poly}
If $\sigma(\log|F|,\varphi)=\sum_{1\le j\le m}a_j$, then there exists a valuation $\nu$ on $\mathbb{C}[z_1,\ldots,z_n]$ such that $\nu(f_j)=a_j$ for every $0\le j\le m$.

Conversely, if there exists a valuation $\nu$ on $\mathbb{C}[z_1,\ldots,z_n]$ satisfying that $\nu(f_j)=a_j$ for every $0\le j\le m$ and $\nu(z_l)>0$ for every $1\le l\le n$, then we have $\sigma(\log|F|,\varphi)=\sum_{1\le j\le m}a_j$.
\end{Corollary}


When $\cap_{1\le j\le m}\{f_j=0\}=\{o\}$, \Cref{C:interpolation-comp-poly} gives a necessary and sufficient condition for the existence of valuative interpolations on $\mathbb{C}[z_1,\ldots,z_n]$.

\begin{Corollary}\label{C:interpolation-comp-poly2}
 If $\cap_{1\le j\le m}\{f_j=0\}=\{o\}$, then there exists a valuation $\nu$ on $\mathbb{C}[z_1,\ldots,z_n]$ such that $\nu(f_j)=a_j$ for every $j$ if and only if $\sigma(\log|F|,\varphi)=\sum_{1\le j\le m}a_j$.
\end{Corollary}

The following example shows that the condition ``$\cap_{1\le j\le m}\{f_j=0\}=\{o\}$" can not be removed.
\begin{Example}
    Let $h_1=z_1,$ $h_2=z_2$, $h_3=z_1z_2$, $g_1=z_1-1$ and $g_2=z_2-1$ on $\mathbb{C}^2$. Take $\{f_1,f_2,f_3,f_4,f_5,f_6\}=\{h_1g_1,h_1g_2,h_2g_1,h_2g_2,h_3g_1,h_3g_2\}$ and $a_1=\ldots=a_6=1$.  It is clear that 
    $$\cap_{1\le l\le 6}\{f_l=0\}=\{o,(1,1)\}.$$ 
    As $F=\prod_{1\le l\le 6}f_l=3\log(|z_1-1||z_2-1|)+O(1)$ near $(1,1)$ and $\varphi=\log(\sum_{1\le l\le 6}|f_l|^{\frac{1}{a_l}})=\log(|z_1-1|+|z_2-1|)+O(1)$ near $(1,1)$, we have
    $$\sup\{c\colon\log|F|\le\varphi+O(1)\text{ near $(1,1)$}\}=6=\sum_{1\le l\le 6}a_l.$$
    By \Cref{C:interpolation-comp-poly}, there exists a valuation $\nu$ on $\mathbb{C}[z_1,z_2]$ such that $\nu(f_l)=a_l$ for every $1\le l\le 6$. But we have 
    $$\sigma(\log|F|,\varphi)>6=\sum_{1\le l\le 6}a_l$$
    since $F=4\log|z_1z_2|+O(1)$ near $o$ and $\varphi=\log(|z_1|+|z_2|+|z_1z_2|)+O(1)$ near $o$. Thus, the condition ``$\cap_{1\le j\le m}\{f_j=0\}=\{o\}$" in \Cref{C:interpolation-comp-poly2} can not be removed.
\end{Example}

\subsubsection{Applications on real cases}
Denote the set of all germs of real analytic functions near the origin $o'\in\mathbb{R}^n$ by $C^{\omega}_{o'}$.

There exists an injective ring homomorphism $P\colon C^{\omega}_{o'}\rightarrow\mathcal{O}_o$ such that
$$P\Big(\sum_{\alpha\in\mathbb{Z}_{\ge0}^n}a_{\alpha}x^{\alpha}\Big)=\sum_{\alpha\in\mathbb{Z}_{\ge0}^n}a_{\alpha}z^{\alpha},$$
where $(x_1,\ldots,x_n)$ and $(z_1,\ldots,z_n)$ are the standard coordinates on $\mathbb{R}^n$ and $\mathbb{C}^n$ respectively, and $\sum_{\alpha\in\mathbb{Z}_{\ge0}^n}a_{\alpha}x^{\alpha}$ is the power series expansion of an arbitrarily given real analytic function near $o'$. It is clear that, for any $(h,o)\in \mathcal{O}_o$, there exists a unique pair of real analytic functions $(h_1, h_2)$ near $o'$ such that
\[h=P(h_1)+\i P(h_2) \quad \ \text{near} \ o.\]

In the following, we discuss the valuative interpolation problem for the real cases.

\begin{Corollary}
	\label{C:interpolation-real}
		Given $m+1$ real analytic functions $\{f_j\}_{0\le j\le m}$ near $o'$, $a_0=0$ and $m$  positive numbers $\{a_j\}_{1\le j\le m}$,
	the following two statements are equivalent:
	\begin{enumerate}
	    \item There exists a valuation $\nu$ on $C^{\omega}_{o'}$ such that $\nu(f_j)=a_j$ for every $0\le j\le m$;
	    \item $\sigma(\log|F|,\varphi)=\sum_{1\le j\le m}a_j$, where
	    \[F\coloneqq \prod_{0\le j\le m}P(f_j), \quad \varphi\coloneqq \log\Big(\sum_{1\le j\le m}|P(f_j)|^{\frac{1}{a_j}}\Big).\]
	\end{enumerate}
\end{Corollary}

\begin{Remark}
Theorem \ref{thm:interpolation} and Corollary \ref{C:interpolation-real} show that for any valuation $\nu$ on $C^{\omega}_{o'}$ and finite collection of real analytic functions $\{f_j\}_{0\le j\le m}$ near $o'$, there exists a valuation $\tilde\nu$ on $\mathcal{O}_o$ such that $\tilde\nu(P(f_j))=\nu(f_j)$ for all $j$.
\end{Remark}

Denote  the restriction of $P$ to $\mathbb{R}[x_1,\ldots,x_n]$ also by $P$.
Let $\{f_j\}_{0\le j\le m}\subset \mathbb{R}[x_1,\ldots,x_n]$, $a_0=0$ and  $\{a_j\}_{1\le j\le m}$ be $m$ positive numbers. Set $F\coloneqq \prod_{0\le j\le m}P(f_j)$ and $\varphi\coloneqq \log(\sum_{1\le j\le m}|P(f_j)|^{\frac{1}{a_j}})$. 

\begin{Corollary}\label{C:interpolation-real-poly}
	 If $\sigma(\log|F|,\varphi)=\sum_{1\le j\le m}a_j$, then there exists a valuation $\nu$ on $\mathbb{R}[x_1,\ldots,x_n]$ such that $\nu(f_j)=a_j$ for every $0\le j\le m$.
	
	Conversely, if there exists a valuation $\nu$ on $\mathbb{R}[x_1,\ldots,x_n]$ satisfying that $\nu(f_j)=a_j$ for every $0\le j\le m$ and $\nu(x_l)>0$ for every $1\le l\le n$, then we have $\sigma(\log|F|,\varphi)=\sum_{1\le j\le m}a_j$.
\end{Corollary}

\subsection{Valuation and relative type}\label{sec:1.2}

\

We recall that a plurisubharmonic function $\varphi$ near $o$ is called a \emph{tame weight} (see \cite{BFJ08}), if $\varphi$ has an isolated singularity at $o$,  $e^{\varphi}$ is continuous, and there exists a constant $C>0$ such that
\begin{equation}\label{ineq-tame}
c_o^{f}(\varphi) \le \sigma(\log|f|, \varphi) + C \quad \forall\,(f,o) \in \mathcal{O}_{o}^*.
\end{equation}
Equivalently, the inequality in \eqref{ineq-tame} requires that for any $t > 0$ and every $(f,o) \in \mathcal{I}(t\varphi)_o$, we have $\log|f| \le (t - C)\varphi + O(1)$ near $o$. Here
\[c_o^f (\varphi)\coloneqq \sup\big\{c\ge 0 \colon |f|^2e^{-2c\varphi} \text{ is integrable near } o\big\}\]
denotes the \emph{jumping number} (see \cite{JON-Mus2012,JON-Mus2014}), and
\[\mathcal{I}(t\varphi)_o\coloneqq \big\{(g,o)\in\mathcal{O}_o\colon |g|^2e^{-2t\varphi} \text{ is integrable near } o\big\}\]
is the \emph{multiplier ideal}.

As shown in previous works, both Zhou valuations (\cite{BGMY-valuation}) and quasimonomial valuations (\cite{BFJ08}) can be realized as relative types with respect to some tame maximal weight. These facts motivate the following question:

\textit{Can one characterize which valuations on $\mathcal{O}_o$ are induced by a tame maximal weight? More precisely, for a given valuation $\nu$ on $\mathcal{O}_o$ centered at $o$, what is a necessary and sufficient condition for $\nu$ to be represented as the relative type with respect to a tame maximal weight near $o$?}

In this paper, we provide a complete answer to this question.

\begin{Theorem}\label{thm:valuationexistsweight}
A valuation $\nu$ on $\mathcal{O}_o$ centered at $o$ satisfies $\nu(f) = \sigma(\log|f|, \varphi\nu)$ for some tame maximal weight $\varphi_\nu$ if and only if $\sup_{(f,o)}\big(\lct^{(f)}(\mathfrak{a}^{\nu}_{\bullet})-\nu(f)\big)<+\infty$. Moreover, when this holds, the representing weight $\varphi_{\nu}$ satisfies $c_o^{f}(\varphi_{\nu}) = \lct^{(f)}(\mathfrak{a}^{\nu}_{\bullet})$ for all $(f,o) \in \mathcal{O}_o^*$.
\end{Theorem}

It is worthy of noting that $\sup_{(f,o)}\big(\lct^{(f)}(\mathfrak{a}^{\nu}_{\bullet})-\nu(f)\big)$ coincides with the notion  $\mathscr{A}(\nu)$ defined in \cite{BGY-valuation2}, and the \emph{log discrepancy} (see \cite{JON-Mus2012}) $A(\nu)$ of $\nu$ satisfies $A(\nu)\ge\mathscr{A}(\nu)$. Therefore, the log discrepancy $A(\nu)<+\infty$ implies $\sup_{(f,o)}\big(\lct^{(f)}(\mathfrak{a}^{\nu}_{\bullet})-\nu(f)\big)<+\infty$.

\subsection{Zhou valuation and quasimonomial valuation}\label{subsec-Zhou.qm}

\

Let $f_{0}=(f_{0,1},\cdots,f_{0,m})$ be a vector,
where $f_{0,1},\cdots,f_{0,m}$ are holomorphic functions near $o$.
Denote by $|f_{0}|^{2}=|f_{0,1}|^{2}+\cdots+|f_{0,m}|^{2}$.
Let $\varphi_{0}$ be a plurisubharmonic function near $o$,
such that $|f_{0}|^{2}e^{-2\varphi_{0}}$ is integrable near $o$.

\begin{Definition}[\cite{BGMY-valuation}]
	\label{def:max_relat}
	We call a plurisubharmonic function $\Phi^{f_0,\varphi_0}_{o,\max}$ ($\Phi_{o,\max}$ for short)  near $o$  a \emph{local Zhou weight related to $|f_{0}|^{2}e^{-2\varphi_{0}}$ near $o$},
	if the following three statements hold:
	\begin{enumerate}
	    \item $|f_{0}|^{2}e^{-2\varphi_{0}}|z|^{2N_{0}}e^{-2\Phi_{o,\max}}$ is integrable near $o$ for large enough $N_{0}\gg 0$;
	    \item $|f_{0}|^{2}e^{-2\varphi_{0}}e^{-2\Phi_{o,\max}}$ is not integrable near $o$;
	    \item for any plurisubharmonic function $\varphi'\geq\Phi_{o,\max}+O(1)$ near $o$ such that $|f_{0}|^{2}e^{-2\varphi_{0}}e^{-2\varphi'}$ is not integrable near $o$, $\varphi'=\Phi_{o,\max}+O(1)$ holds.
	\end{enumerate}
\end{Definition}

The existence of local Zhou weights follows from the strong openness property of multiplier ideal sheaves (see \cite{BGMY-valuation} or Remark \ref{rem:max_existence}).

Let us recall some properties of Zhou weights (see \cite{BGMY-valuation}). For a local Zhou weight $\Phi_{o,\max}$, we have $\Phi_{o,\max}\ge N\log|z|$ near $o$ for some $N>0$, and its corresponding \emph{Zhou valuation} $\nu(\cdot,\Phi_{o,\max})$ is defined by $\nu(f,\Phi_{o,\max})\coloneqq \sigma(\log|f|,\Phi_{o,\max})$ for any $(f,o)\in\mathcal{O}_o$. If $\Phi_{o,\max}$ is a local Zhou weight related to $|f_0|^2e^{-2\varphi_0}$ near $o$, then $(1+\sigma(\varphi,\Phi_{o,\max}))\Phi_{o,\max}$ is a local Zhou weight related to $|f_0|^2$ near $o$.

We give a sufficient condition for a Zhou valuation to be quasimonomial.

\begin{Proposition}\label{p:quasi}
	Let $\varphi$ be a plurisubharmonic function near $o$ with analytic singularities such that $c_o^{f_0}(\varphi)=1.$
	 Let $\Phi_{o,\max}$ be a local Zhou weight related to $|f_0|^2$ near $o$ such that $\Phi_{o,\max}\ge\varphi+O(1)$ near $o$. Denote by $\nu$ the Zhou valuation corresponding to the Zhou weight $\Phi_{o,\max}$. Then $\nu$ is a quasimonomial valuation.
\end{Proposition}

In \cite{BFJ08}, Boucksom--Favre--Jonsson established a fundamental result about the relations between multiplier ideal sheaves and relative types by using quasimonomial valuations. In \cite{BGMY-valuation}, we proved that the result also holds for Zhou valuations. Furthermore, the following theorem shows that the result  holds for quasimonomial Zhou valuations

\begin{Theorem}
	\label{thm:multi-valua,quasi}
	Let $u,$ $v$ be two plurisubharmonic functions near $o$. Then the following  statements are equivalent:
	\begin{enumerate}
	    \item there exists a plurisubharmonic function $\varphi_{0}$ near $o$ ($\not\equiv-\infty$) and two sequences of numbers  $\{t_{i,j}\}_{j\in\mathbb{Z}_{\ge0}}$ ($t_{i,j}\rightarrow+\infty$ when $j\rightarrow+\infty$, $i=1,2$) such that $\lim_{j\rightarrow+\infty}\frac{t_{1,j}}{t_{2,j}}=1$ and $$\mathcal{I}(\varphi_{0}+t_{1,j}v)_o\subset\mathcal{I}(\varphi_{0}+t_{2,j}u)_o$$ for every $j$;
	    \item for any plurisubharmonic function $\varphi_0$ near $o$ and any $t>0$, we have $$\mathcal{I}(\varphi_0+tv)_o\subset\mathcal{I}(\varphi_0+tu)_o;$$
	    \item for any local Zhou weight $\Phi_{o,\max}$ near $o$ whose corresponding Zhou valuation is quasimonomial, we have $$\sigma(u,\Phi_{o,\max})\le\sigma(v,\Phi_{o,\max}).$$
	\end{enumerate}
\end{Theorem}

\subsection{Valuations related to Tian functions}
\label{sec:1.1}

\

Let $f_0$ and $\varphi_0$ be as those defined in \Cref{subsec-Zhou.qm}. Let $\varphi$, $\psi$ be any two plurisubharmonic functions defined near $o$. Denote
\[c_{o}(\varphi,t\psi;f_0,\varphi_0)\coloneqq\sup\big\{c\colon |f_{0}|^{2}e^{-2\varphi_{0}}e^{2t\psi}e^{-2c\varphi} \ \text{is integrable near} \ o\big\}.\]

We recall the following notation of \emph{Tian function} (see \cite{BGMY-valuation}):
$$\Tn (t;f_0,\varphi_0,\varphi,\psi)\coloneqq c_{o}(\varphi,t\psi;f_0,\varphi_0), \quad \text{for} \ t\in\mathbb{R}.$$
Note that we allow $t$ to take negative value. We also denote $\Tn (t;f_0,\varphi_0,\varphi,\psi)$ by $\Tn (t;\psi)$ or $\Tn (t)$ for short. It is easy to see that $\Tn (t)$ is non-decreasing. It follows from the H\"{o}lder inequality that $\Tn (t)$ is concave in $t$.

If $\varphi$ is a local Zhou weight related to $|f_0|^2e^{-\varphi_0}$ near $o$, then $\Tn (t)$ is differentiable at $t=0$ and $\Tn'(0)=\sigma(\psi,\varphi)$ for any plurisubharmonic function $\psi$ near $o$ (see \Cref{prop:concave_B}). 

Now we assume $\varphi$ is a general plurisubharmonic function near $o$. If $\Tn (t)$ is differentiable at $t=0$ for any plurisubharmonic function $\psi$, we set
$$\tilde{\sigma}(\psi,\varphi)\coloneqq \Tn'(0)=\lim_{t\to 0}\frac{\Tn (t)-\Tn (0)}{t}.$$
We show that $\tilde{\sigma}(\cdot,\varphi)$ satisfies the tropical multiplicativity and tropical additivity properties in the following theorem, which has been proved in \cite{BGMY-valuation} when $\varphi$ is a local Zhou weight.
\begin{Theorem}
	\label{thm:derivative} Assume that  $\Tn (t)$ is differentiable at  $t= 0$ for any plurisubharmonic function $\psi$ near $o$. Then
	 $\tilde{\sigma}(\psi,\varphi)$ defined above satisfies
	 \begin{enumerate}
	     \item for any $c_{1}\geq 0$ and $c_{2}\geq 0$, $\tilde{\sigma}(c_{1}\psi_{1}+c_{2}\psi_{2},\varphi)=c_{1}\tilde{\sigma}(\psi_{1},\varphi)+c_{2}\tilde{\sigma}(\psi_{2},\varphi)$;
	     \item $\tilde{\sigma}(\max\{\psi_1,\psi_2\},\varphi)=
	\min\{\tilde{\sigma}(\psi_1,\varphi),\tilde{\sigma}(\psi_2,\varphi)\}$, where $\psi_1,\psi_2$ are germs of plurisubharmonic functions near $o$;
	    \item $\tilde{\sigma}(\log|f_{1}+f_{2}|,\varphi)\geq\min\big\{\tilde{\sigma}(\log|f_{1}|,\varphi),\tilde{\sigma}(\log|f_{2}|,\varphi)\big\}$, where $f_{1},f_{2}\in \mathcal{O}_o$.
	 \end{enumerate}
\end{Theorem}

Statements $(1)$ and $(3)$ in Theorem \ref{thm:derivative} show that $\tilde\sigma(\log|\cdot|,\varphi)$ is a valuation on $\mathcal{O}_o$.
For the isolated singularity case, the following proposition shows that the valuation is a Zhou valuation. Thus, $\tilde\sigma(\log|\cdot|,\varphi)$ is a generalization of Zhou valuation.

\begin{Proposition}\label{p:differentiable-zhou valuation}Assume that  there exists $N\gg 0$ such that $|f_0|^2|z|^{2N}e^{-2\varphi_0-2\varphi}$ is integrable near $o$.
	If $\Tn (t;f_0,\varphi_0,\varphi,\log|g|)$ is differentiable at $t=0$ for any holomorphic function $g$ near $o$, then $\nu(g)\coloneqq \Tn '(0;f_0,\varphi_0,\varphi,\log|g|)$ is a Zhou valuation on $\mathcal{O}_o$.
\end{Proposition}


\section{Preliminaries}

\subsection{Strong openness property of multiplier ideal sheaves}

Let $f$ be a holomorphic function near $o$, and let $\varphi$  be a plurisubharmonic function near $o$.
Let us define the minimal $L^2$ integral related to multiplier ideal sheaf (see \cite{GZopen-effect,guan-effect})
\[C_{f,\varphi}(U)\coloneqq \inf\left\{\int_{U} |\tilde{f}|^{2} \colon (\tilde{f}-f,o)\in \mathcal{I}(\varphi)_{o} \ \& \ \tilde{f}\in\mathcal{O}(U)\right\},\]
where $U\subseteq \Delta^{n}$ is a domain with $o\in U$.

Let $\{\phi_{m}\}_{m\in\mathbb{N}^{+}}$ be a sequence of negative plurisubharmonic functions on $\Delta^{n}$,
which is convergent to a negative Lebesgue measurable function $\phi$ on $\Delta^{n}$ in Lebesgue measure.

In \cite{GZopen-effect}, Guan--Zhou presented the following lower semicontinuity property of plurisubharmonic functions with a multiplier.

\begin{Proposition}[{see \cite[Proposition 1.8]{GZopen-effect}}]
	\label{p:effect_GZ}
	Let $f$ be a holomorphic function near $o$, and let $\varphi_0$ be a plurisubharmonic function near $o$.
	Assume that for any small enough neighborhood $U$ of $o$,
	the pairs $(f,\phi_{m})$ $(m\in\mathbb{Z}_{>0})$ satisfies
	\begin{equation}
		\label{equ:A}
		\inf_{m}C_{f,\varphi_0+\phi_{m}}(U)>0.
	\end{equation}
	Then $|f|^{2}e^{-2\varphi_{0}}e^{-2\phi}$ is not integrable near $o$.
\end{Proposition}

The Noetherian property of multiplier ideal sheaves (see \cite{demailly-book}) shows that
\begin{Remark}
	\label{rem:effect_GZ}
	Assume that
	\begin{enumerate}
	    \item $\phi_{m+1}\geq\phi_{m}$ holds for any $m$;
	    \item $|f|^{2}e^{-2\varphi_{0}}e^{-2\phi_{m}}$ is not integrable near $o$ for any $m$.
	\end{enumerate}
	Then the inequality \eqref{equ:A} holds.
\end{Remark}

In particular, Proposition \ref{p:effect_GZ} and Remark \ref{rem:effect_GZ} imply the \emph{strong openness property} (cf. \cite{GZopen-c})
$$\mathcal{I}(\varphi)=\mathcal{I}_+(\varphi)$$ 
of plurisubharmonic functions $\varphi$.

\subsection{Some results on plurisubharmonic functions, multiplier ideal sheaves and valuations}

Firstly, let us recall a basic result on plurisubharmonic functions.

\begin{Lemma}[\cite{Skda72}; see also \cite{demailly2010}]\label{l:Lelong}
	For any plurisubharmonic function $\varphi$ near $o$, if the Lelong number $\nu_o(\varphi)<1$, then $e^{-2\varphi}$ is locally integrable near $o$.
\end{Lemma}

The following is the so-called \emph{Demailly's approximation theorem}.

\begin{Lemma}[see \cite{demailly2010}]\label{l:appro-Berg}
	Let $D\subset\mathbb{C}^n$ be a bounded pseudoconvex domain, and let $\varphi\in \mathrm{PSH}(D)$. For any positive integer $m$, let $\{\sigma_{m,k}\}_{k=1}^{\infty}$ be an orthonormal basis of $A^2(D,2m\varphi)\coloneqq \big\{f\in\mathcal{O}(D)\colon \int_{D}|f|^2e^{-2m\varphi}d\lambda<+\infty\big\}$. Set
	$$\varphi_m\coloneqq \frac{1}{2m}\log\sum_{k=1}^{\infty}|\sigma_{m,k}|^2$$
	on $D$. Then there exist two positive constants $c_1$ (depending only on $n$ and the diameter of $D$) and $c_2$ such that
	\begin{equation}
		\label{eq:0911b} \varphi(z)-\frac{c_1}{m}\le\varphi_m(z)\le\sup_{|\tilde z-z|<r}\varphi(\tilde z)+\frac{1}{m}\log\frac{c_2}{r^n}
	\end{equation}
	for any $z\in D$ satisfying $\{\tilde z\in\mathbb{C}^n\colon|\tilde z-z|<r\}\subset\subset D$. Especially, $\varphi_m$ converges to $\varphi$ pointwise and in $L^1_{\mathrm{loc}}$ on $D$.
\end{Lemma}

The strong openness property implies that the multiplier ideal sheaf is essentially with analytic singularities.

\begin{Theorem}[\cite{GZopen-c}]\label{thm:analytic weight}
	Let $\varphi$ be a plurisubharmonic function near $o$. Then $\mathcal{I}(\varphi)_o=\mathcal{I}(\frac{1}{N}\log|J_{N+1}|)_o$ for large enough $N\gg 0$, where $J_{N+1}\coloneqq \mathcal{I}((N+1)\varphi)_o=(f_1,\ldots,f_M)_{o}$ and $|J_{N+1}|^2=\sum_{j=1}^M|f_j|^2$.
\end{Theorem}
\begin{proof} For the convenience of readers, we give a proof.

There exists a neighborhood $V$ of $o$ such that 
\begin{equation}\nonumber
    \begin{split}
        &\int_{V}e^{-2\varphi}-e^{-2\max\{\varphi,\frac{1}{N}\log|J_{N+1}|\}}\\
        =& \ \int_{\{\log|J_{N+1}|\ge N\varphi\}\cap V}e^{-2\varphi}-e^{-2\max\{\varphi,\frac{1}{N}\log|J_{N+1}|\}}\\
        \le& \ \int_{\{\log|J_{N+1}|\ge N\varphi\}\cap V}e^{-2\varphi}|J_{N+1}|^2e^{-2N\varphi}\\
        \le& \ \int_{V}|J_{N+1}|^2e^{-2(N+1)\varphi}\\
        <& \ +\infty,
    \end{split}
\end{equation}
which implies that
\begin{equation}
    \label{eq:0917c}
    \mathcal{I}(\varphi)_o=\mathcal{I}\Big(\max\big\{\varphi,\frac{1}{N}\log|J_{N+1}|\big\}\Big)_o.
\end{equation}
 Lemma \ref{l:appro-Berg} shows that 
\begin{equation}
    \label{eq:0917a}
    \varphi\le\frac{1}{N}\log|J_{N+1}|+O(1)
\end{equation}
near $o$. By the strong openness property, there exists $\epsilon_0>0$ such that
\begin{equation}
    \label{eq:0917b}
    \mathcal{I}(\varphi)_{o}=\mathcal{I}\big((1+\epsilon_0)\varphi\big)_{o}.
\end{equation}
Taking $N\ge\frac{1}{\epsilon_0}$, we have $\frac{1}{N}\log|J_{N+1}|\ge(1+\epsilon_0)\frac{1}{N+1}\log|J_{N+1}|$.
Combining with the equalities \eqref{eq:0917c}, \eqref{eq:0917a} and \eqref{eq:0917b}, we have
\begin{equation}\nonumber
    \begin{split}
        \mathcal{I}\Big(\frac{1}{N}\log|J_{N+1}|\Big)_o&\subset\mathcal{I}\Big(\max\big\{\frac{1}{N}\log|J_{N+1}|,\varphi\big\}\Big)_o\\
        &=\mathcal{I}(\varphi)_o\\
        &=\mathcal{I}\big((1+\epsilon_0)\varphi\big)_{o}\\
        &\subset\mathcal{I}\Big((1+\epsilon_0)\frac{1}{N}\log|J_{N+1}|\Big)_{o}\\
        &\subset\mathcal{I}\Big(\frac{1}{N}\log|J_{N+1}|\Big)_o,
    \end{split}
\end{equation}
which implies $\mathcal{I}(\varphi)_o=\mathcal{I}(\frac{1}{N}\log|J_{N+1}|)_o$.
Theorem \ref{thm:analytic weight} holds.
\end{proof}

The following lemma gives a convergence property of jumping numbers.

\begin{Lemma}
	[{see \cite[Lemma 6.4]{BGMY-valuation}}]
	\label{l:lct-N}
	Let $\varphi$ be any plurisubharmonic function near $o$, and let $f\not\equiv0$ be any holomorphic function near $o$. Then
	\[\lim_{N\rightarrow+\infty}c_o^{f}\big(\max\{\varphi,N\log|z|\}\big)=c_o^{f}(\varphi).\]
\end{Lemma}

Let us recall two results about the upper envelope of family of plurisubharmonic functions.

\begin{Lemma}[\emph{Choquet's lemma}, see \cite{demailly-book}]
	\label{lem:Choquet} Every family $(u_{\alpha})$ of upper-semicontinuous functions has a countable subfamily $(v_{j})=(u_{\alpha(j)})$,
	such that its upper envelope $v=\sup_{j}v_{j}$ satisfies $v\leq u\leq u^{*} = v^{*}$,
	where $u=\sup_{\alpha}u_{\alpha}$, $u^{*}(z)\coloneqq \lim_{\varepsilon\to0}\sup_{\mathbb{B}^{n}(z,\varepsilon)}u$
	and $v^{*}(z)\coloneqq \lim_{\varepsilon\to0}\sup_{\mathbb{B}^{n}(z,\varepsilon)}v$ are the regularizations of $u$ and $v$.
\end{Lemma}

\begin{Proposition}[{see \cite[Proposition 4.24]{demailly-book}}]
	\label{pro:Demailly}
	If all $(u_{\alpha})$ are subharmonic, the upper regularization $u^{*}$ is subharmonic
	and equals almost everywhere to $u$.
\end{Proposition}

We recall Skoda's division theorem.

\begin{Theorem}[{see \cite[Theorem 12.13]{demailly2010}}]
	\label{thm:Skoda}
	Let $\Omega$ be a complete K\"ahler open subset of $\mathbb{C}^n$ and $\varphi$ be a plurisubharmonic function on $\Omega$. Set $m=\min\{n,r-1\}$. Then for every holomorphic function $f$ on $\Omega$ such that
	$$I\coloneqq \int_{\Omega}|f|^2|g|^{-2(m+1+\epsilon)}e^{-\varphi}<+\infty,$$
	there exist holomorphic functions $(h_1,\ldots,h_r)$ on $\Omega$ such that $f=\sum_{1\le j\le r}g_j h_j$ and
	$$\int_{\Omega}|h|^2|g|^{-2(m+\epsilon)}e^{-\varphi}\le\Big(1+\frac{m}{\epsilon}\Big)I.$$
\end{Theorem}

Using the above division theorem, we give some lemmas on valuations.

\begin{Lemma}\label{l:relativetype}
	For any valuation $\nu$ on $\mathcal{O}_o$ and holomorphic functions $g,f_1,\ldots,f_r$ near $o$ satisfying $\nu(f_j)>0$ for $1\le j\le r$, we have
	$\sigma(\log|g|,\varphi)\le \nu(g)$, where $\varphi\coloneqq \log\big(\sum_{1\le j\le r}|f_j|^{\frac{1}{\nu(f_j)}}\big)$.
\end{Lemma}

\begin{proof}
	We prove this lemma by contradiction: if not, there exists $\delta_1>0$ such that
	$\log|g|\le (\nu(g)+\delta_1)\varphi+O(1)$ near $o$. Then we can find some $\delta_2>0$ and a set of positive integers $\{m_0,\ldots,m_r\}$ such that
	\[\log|g|\le \frac{\nu(g)+\delta_2}{m_0}\tilde\varphi+O(1) \quad  \ \text{near} \ o,\]
	where
	$\tilde\varphi\coloneqq\log\big(\sum_{1\le j\le r}|f_j^{m_j}|\big)$ and $\nu(f_j^{m_j})\ge m_0$ for $1\le j\le r$. For any positive integer $l$, we have
	\[\log|g^l|\le \frac{\nu(g^l)+l\delta_2}{m_0}\tilde\varphi+O(1)\quad \ \text{near} \ o,\]
	which implies that $|g^l|^2e^{-2\frac{\nu(g^l)+l\delta_2+\epsilon}{m_0}\tilde\varphi}$ is integrable near $o$ for small enough $\epsilon>0$ (independent of $l$).
	Denote  by $I$ the ideal generated by $\{(f_j^{m_j},o)\}_{1\le j\le r}$ in $\mathcal{O}_o$. Let $|F|\coloneqq \sum_{1\le j\le r}|f_j^{m_j}|$, $\epsilon_0\coloneqq \frac{\epsilon}{m_0}$ and $m\coloneqq \min\{n,r-1\}$. Since $|g^l|^2e^{-2\frac{\nu(g^l)+l\delta_2+\epsilon}{m_0}\tilde\varphi}$ is integrable on a small open set $U$  containing $o$, we have
	\begin{equation*}
		\begin{split}
			& \int_{U}|g^l|^2|F|^{-2(m+1)-2\epsilon_0}|F|^{-2\big(\frac{\nu(g^l)+l\delta_2}{m_0}-(n+1)\big)} \\
			\le&C \int_{U}|g^l|^2|F|^{-2(n+1)-2\epsilon_0}|F|^{-2\big(\frac{\nu(g^l)+l\delta_2}{m_0}-(n+1)\big)} \\
			= &C \int_{U}|g^l|^2e^{-2\frac{\nu(g^l)+l\delta_2+\epsilon}{m_0}\tilde\varphi}<+\infty.
		\end{split}
	\end{equation*}
	By Theorem \ref{thm:Skoda}, there exist holomorphic functions $(h_1,\ldots,h_r)$ on $U$ such that $g^l=\sum_{1\le j\le r}f_j^{m_j}h_j$ and
	\begin{equation}
		\label{eq:est by skoda}
		\begin{split}
			&\int_{U}|h|^2|F|^{-2(m+1)-2\epsilon_0}|F|^{-2\big(\frac{\nu(g^l)+l\delta_2}{m_0}-1-(n+1)\big)}\\
			=& \int_{U}|h|^2|F|^{-2m-2\epsilon_0}|F|^{-2\big(\frac{\nu(g^l)+l\delta_2}{m_0}-(n+1)\big)}< +\infty.
		\end{split}
	\end{equation}
	Note that $g^l=\sum_{1\le j\le r}f_j^{m_j}h_j$ implies $(g^l,o)\in I$. It follows from the estimate \eqref{eq:est by skoda} that, using Theorem \ref{thm:Skoda} for each $h_j$ $(1\le j\le r)$, we can get $h_j=\sum_{1\le k\le r}h_{j,k}f_j^{m_j}\in I$ (which implies $(g^l,o)\in I^2$) and
	\begin{equation*}
		\int_{U}\sum_{1\le k\le r}|h_{j,k}|^2|F|^{-2m-2\epsilon_0}|F|^{-2\big(\frac{\nu(g^l)+l\delta_2}{m_0}-1-(n+1)\big)}< +\infty.
	\end{equation*}
	Hence, using Theorem \ref{thm:Skoda} repeatedly, we finally have $(g^l,o)\in I^{\big\lfloor\frac{\nu(g^l)+l\delta_2}{m_0}\big\rfloor-n-1}$ for $l\gg 1$. As $\nu(f_j^{m_j})\ge m_0$ for $1\le j\le r$, we get
	\begin{equation}\nonumber
		\begin{split}
			\nu(g^l)&\ge m_0\Big(\Big\lfloor\frac{\nu(g^l)+l\delta_2}{m_0}\Big\rfloor-n-1\Big)
			\\&\ge m_0\Big(\frac{\nu(g^l)+l\delta_2}{m_0}-n-2\Big)\\
			&=\nu(g^l)+l\delta_2-m_0(n+2),
		\end{split}
	\end{equation}
	which implies $l\le\frac{m_0(n+2)}{\delta_2}$. This is a contradiction since we can choose $l$ arbitrarily large. Then we  have $\sigma(\log|g|,\varphi)\le \nu(g)$.
\end{proof}

For valuations on the polynomial ring $\mathbb{C}[z_1,\ldots,z_n]$, we have the following similar result:

\begin{Lemma}\label{l:relativetype2}
	Let $\nu$ be a valuation on $\mathbb{C}[z_1,\ldots,z_n]$ satisfying $\nu(z_j)>0$ for every $1\le j\le n$.  For any polynomials $g,f_1,\ldots,f_r\in\mathbb{C}[z_1,\ldots,z_n]$  satisfying $\nu(f_j)>0$ for $1\le j\le r$, we have
	$\sigma(\log|g|,\varphi)\le \nu(g)$, where $\varphi\coloneqq \log\big(\sum_{1\le j\le r}|f_j|^{\frac{1}{\nu(f_j)}}\big)$.
\end{Lemma}

\begin{proof}
	We prove this lemma by contradiction: if not, there exists $\delta_1>0$ such that
	$\log|g|\le (\nu(g)+\delta_1)\varphi+O(1)$ near $o$. There exist $\delta_2>0$ and positive integers $\{m_0,\ldots,m_r\}$ such that
	$\log|g|\le \frac{\nu(g)+\delta_2}{m_0}\tilde\varphi+O(1)$ near $o$, where
	$\tilde\varphi=\log(\sum_{1\le j\le r}|f_j^{m_j}|)$ and $\nu(f_j^{m_j})\ge m_0$ for $1\le j\le r$.
	For any large positive integer $l$, the proof of Lemma \ref{l:relativetype} shows that $$g^l=\sum_{|\alpha|=\big\lfloor\frac{\nu(g^l)+l\delta_2}{m_0}\big\rfloor-n-1}g_{\alpha}f_{1}^{m_1\alpha_1}\cdots f_{r}^{m_r\alpha_r}$$
	near $o$, where $g_{\alpha}$ is holomorphic function near $o$ and $\alpha=(\alpha_1,\ldots,\alpha_r)\in\mathbb{Z}^r_{\ge0}$.
	
	As $g,f_1,\ldots,f_r$ are polynomials, for any integer $N>0$, there exist polynomials $\{g_{N,\alpha}\}_{\alpha}$ such that
	$$h_N\coloneqq g^l-\sum_{|\alpha|=\big\lfloor\frac{\nu(g^l)+l\delta_2}{m_0}\big\rfloor-n-1}g_{N,\alpha}f_{1}^{m_1\alpha_1}\cdots f_{r}^{m_r\alpha_r}\in \mathbb{C}^{(N)}[z_1,\ldots,z_n],$$
	where $\mathbb{C}^{(N)}[z_1,\ldots,z_n]$ denotes the set of polynomials with degree larger than $N$. Since $\nu(z_j)>0$ for every $1\le j\le n$, we have
	\[\nu(g^l)\ge\min\Big\{\nu(h_N),\sum_{1\le j\le r}m_0\alpha_j\Big\}\ge\nu(g^l)+l\delta_2-m_0(n+2)\]
	by taking large enough $N$, which implies $l\le\frac{m_0(n+2)}{\delta_2}$. It is a contradiction. Then we have $\sigma(\log|g|,\varphi)\le \nu(g)$.
\end{proof}

Denote the set of all germs of real analytic functions near the origin $o'\in\mathbb{R}^n$ by $C^{\omega}_{o'}$. There exists an injective ring homomorphism $P\colon C^{\omega}_{o'}\rightarrow\mathcal{O}_o$, which satisfies
$$P\Big(\sum_{\alpha\in\mathbb{Z}_{\ge0}^n} a_{\alpha}x^{\alpha}\Big)=\sum_{\alpha\in\mathbb{Z}_{\ge0}^n}a_{\alpha}z^{\alpha},$$
where $(x_1,\ldots,x_n)$ and $(z_1,\ldots,z_n)$ are the standard coordinates in $\mathbb{R}^n$ and $\mathbb{C}^n$ respectively, and $\sum_{\alpha\in\mathbb{Z}_{\ge0}^n}a_{\alpha}x^{\alpha}$ is the power series expansion of arbitrary real analytic function near $o'$. It is clear that
for any $(g,o)\in \mathcal{O}_o$, there exists a unique pair of real analytic functions $(g_1, g_2)$ near $o'$ such that
\[g=P(g_1)+\i P(g_2) \quad \ \text{near} \ o.\]

Then we have
\begin{Lemma}
	\label{c:relativetype-real}
	Let $\nu$ be a valuation on $C^{\omega}_{o'}$.  For any real analytic functions $g, f_1, \ldots, f_r$ near $o'$ satisfying $\nu(f_j)>0$ for $1\le j\le r$, we have
	$\sigma(\log|\tilde g|,\varphi)\le \nu(g)$, where $\varphi\coloneqq \log\big(\sum_{1\le j\le r}|\tilde f_j|^{\frac{1}{\nu(f_j)}}\big)$,  $\tilde f_j\coloneqq P(f_j)$ for $1\le j\le r$ and $\tilde g\coloneqq P(g)$.
\end{Lemma}

\begin{proof}
	We prove this lemma by contradiction: if not, there exists $\delta_1>0$ such that
	$\log|\tilde g|\le (\nu(g)+\delta_1)\varphi+O(1)$ near $o$. There exist $\delta_2>0$ and positive integers $\{m_0,\ldots,m_r\}$ such that
	$\log|\tilde g|\le \frac{\nu(g)+\delta_2}{m_0}\tilde\varphi+O(1)$ near $o$, where
	$\tilde\varphi=\log(\sum_{1\le j\le r}|\tilde f_j^{m_j}|)$ and $\nu(f_j^{m_j})\ge m_0$ for $1\le j\le r$.
	For any large positive integer $l$, the proof of Lemma \ref{l:relativetype} shows that
	$$\tilde g^l=\sum_{|\alpha|=\left\lfloor\frac{\nu(g^l)+l\delta_2}{m_0}\right\rfloor-n-1}g_{\alpha}\tilde f_{1}^{m_1\alpha_1}\cdots \tilde f_{r}^{m_r\alpha_r}$$
	near $o$, where $g_{\alpha}$ is holomorphic function near $o$ and $\alpha=(\alpha_1,\ldots,\alpha_r)\in\mathbb{Z}^r_{\ge0}$.
	
	For any $\alpha$, there exists a unique pair of real analytic functions $(g_{1,\alpha}, g_{2,\alpha})$ near $o'$ such that
	$$g_{\alpha}=P(g_{1,\alpha})+\i P(g_{2,\alpha})$$
	near $o$. Note that $\tilde f_j\coloneqq P(f_j)$ for $1\le j\le r$   and $\tilde g\coloneqq P(g)$. It follows from $$P(g^l)=\sum_{|\alpha|=\left\lfloor\frac{\nu(g^l)+l\delta_2}{m_0}\right\rfloor-n-1}(P(g_{1,\alpha})+\i P(g_{2,\alpha}))P(f_{1}^{m_1\alpha_1})\cdots P(f_{r}^{m_r\alpha_r})$$  that
	$$g^l=\sum_{|\alpha|=\big\lfloor\frac{\nu(g^l)+l\delta_2}{m_0}\big\rfloor-n-1}g_{1,\alpha} f_{1}^{m_1\alpha_1}\cdots  f_{r}^{m_r\alpha_r}.$$
	As $\nu(f_j^{m_j})\ge m_0$ for $1\le j\le r$, we get
	\begin{equation}\nonumber
		\begin{split}
			\nu(g^l)&\ge m_0\Big(\Big\lfloor\frac{\nu(g^l)+l\delta_2}{m_0}\Big\rfloor-n-1\Big)
			\\&\ge m_0\Big(\frac{\nu(g^l)+l\delta_2}{m_0}-n-2\Big)\\
			&=\nu(g^l)+l\delta_2-m_0(n+2),
		\end{split}
	\end{equation}
	which implies $l\le\frac{m_0(n+2)}{\delta_2}$. It is a contradiction. Then we have $\sigma(\log|\tilde g|,\varphi)\le \nu(g)$.
\end{proof}



Denote  the restriction of $P$ on $\mathbb{R}[x_1,\ldots,x_n]$ also by $P$.
Combining the proofs of Lemma \ref{l:relativetype2} and Lemma \ref{c:relativetype-real}, we can immediately obtain the following

\begin{Lemma}
	\label{c:relativetype2-real}
	Let $\nu$ be a valuation on $\mathbb{R}[x_1,\ldots,x_n]$ satisfying $\nu(x_j)>0$ for any $1\le j\le n$.  For any real polynomials $g,f_1,\ldots,f_r$ satisfying $\nu(f_j)>0$ for $1\le j\le r$, we have
	$\sigma(\log|\tilde g|,\varphi)\le \nu(g)$, where $\varphi\coloneqq \log\big(\sum_{1\le j\le r}|\tilde f_j|^{\frac{1}{\nu(f_j)}}\big)$, $\tilde g\coloneqq P(g)$ and $\tilde f_j\coloneqq P(f_j)$ for $1\le j\le r$.
\end{Lemma}

\subsection{Tian function}
In this section, we discuss a special condition for Tian functions. 

Firstly, let us recall the definition of Tian functions and a useful lemma.
Let $f_0$ be a vector of holomorphic functions near $o$, and $\varphi$ be a plurisubharmonic function near $o$. Let $f$ be a holomorphic function near $o$.
Recall the Tian function
$$\Tn(t;f_0,f,\varphi)\coloneqq \sup\big\{c\colon |f_0|^2|f|^{2t}e^{-2c\varphi}\text{ is integrable near }o\big\}$$
for any $t\in \mathbb{R}$. The H\"older inequality implies the concavity of $\Tn(t;f_0,f,\varphi)$.
\begin{Lemma}[{see \cite[Lemma 2.14]{BGMY-valuation}}]\label{l:1}
	Assume that $\Tn(t;f_0,f,\varphi)$ is strictly increasing on $(t_0-\delta,t_0+\delta)$ for some fixed $t_0\in\mathbb{R}$ and small $\delta>0$. Then
	$$\Tn\Big(t_0;f_0,f,\max\big\{\varphi,\frac{\log|f|}{a}\big\}\Big)=\Tn(t_0;f_0,f,\varphi)$$
	holds for any $a\in(0,\Tn'_-(t_0;f_0,f,\varphi)]$.
\end{Lemma}

\begin{Example}\label{example}
	Let $f=\sum_{\alpha\in I}b_{\alpha}z^{\alpha}$ and $g=\sum_{\alpha\in \tilde I}\tilde b_{\alpha}z^{\alpha}$ be  holomorphic functions near $o$, and let   
	\[\varphi=\max\{a_1\log|z_1|,\ldots,a_r\log|z_r|\},\]          
	where $a_j>0$, and $I,\tilde I\subset\mathbb{Z}_{\ge0}$ such that $b_{\alpha}\not=0$ for any $\alpha\in I$ and $\tilde b_{\alpha}\not=0$ for any $\alpha\in\tilde I$.    
The following Lemma \ref{l:I(varphi)} shows that 
	\[\Tn(m;g,f,\varphi)=\inf_{\alpha\in\tilde I}\sum_{1\le j\le r}\frac{\alpha_j}{a_j}+m\inf_{\alpha\in I}\sum_{1\le j\le r}\frac{\alpha_j}{a_j}\]
	for any integer $m\ge0$. Due to the concavity of $\Tn(t;g,f,\varphi)$, we have $$\Tn(t;1,f,\varphi)=\inf_{\alpha\in\tilde I}\sum_{1\le j\le r}\frac{\alpha_j}{a_j}+t\inf_{\alpha\in I}\sum_{1\le j\le r}\frac{\alpha_j}{a_j}$$ for any $t\ge0$. By a direct calculation, we have 
	\[\log|f|\le\left(\inf_{\alpha\in I}\sum_{1\le j\le r}\frac{\alpha_j}{a_j}\right)\varphi+O(1)\]
	near $o$, which shows
	\[\sigma(\log|f|,\varphi)\ge \inf_{\alpha\in I}\sum_{1\le j\le r}\frac{\alpha_j}{a_j}.\]
	As $\Tn(t;g,f,\varphi)\ge t\sigma(\log|f|,\varphi)$ for any $t\ge0$, we have 
	\[\sigma(\log|f|,\varphi)= \inf_{\alpha\in I}\sum_{1\le j\le r}\frac{\alpha_j}{a_j}.\]
\end{Example}

\begin{Lemma}[see \cite{guan-20}]\label{l:I(varphi)}
    Let $\varphi=\max\{a_1\log|z_1|,\ldots,a_r\log|z_r|\}$ on $\mathbb{C}^n$. Then \[\mathcal{I}(\varphi)_o=\Big\{(f,o)\in\mathcal{O}_o \colon b_{\alpha}=0 \ \text{if} \ \sum_{1\le j\le r}\frac{\alpha_j+1}{a_j}\le 1, \ \text{where} \ f=\sum_{\alpha\in\mathbb{Z}_{\ge0}} b_{\alpha}z^{\alpha}\Big\}.\]
\end{Lemma}

We say $(f,\varphi;|f_0|^2)$ satisfies \textbf{condition (T)} if the following statements hold:
\begin{enumerate}
    \item $\log|f|\le a\varphi+O(1)$ near $o$ and $\Tn(t;f_0,f,\varphi)=at+b$ for any  $t\ge0$;
    \item $\Tn(t;f_0,f,\varphi)$ is differentiable at $t=0$.
\end{enumerate}

\Cref{example} shows that  $(f,\varphi;f)$ satisfies condition (T), where $f$ is any holomorphic function near $o$ and $\varphi$ is defined as in \Cref{example}.

	If $(f,\varphi;|f_0|^2)$ satisfies condition (T), by definition, we have $\sigma(\log|f|,\varphi)\ge a$. It follows from $\Tn(t;f_0,f,\varphi)=at+b$ that $\sigma(\log|f|,\varphi)\le a$. Then we get $\sigma(\log|f|,\varphi)= a$.

The following lemma shows that condition (T) satisfies a multiplicative property. 

\begin{Lemma}\label{l:2}
	If $(f_i,\varphi;|f_0|^2)$ satisfies condition (T) for $i=1,2$, then we have that $(f_1f_2,\varphi;|f_0|^2)$ also satisfies condition (T) and $$\sigma(\log|f_1f_2|,\varphi)=\sigma(\log|f_1|,\varphi)+\sigma(\log|f_2|,\varphi).$$
\end{Lemma}

\begin{proof}
Set
$$A(s,t)\coloneqq \sup\big\{c\colon |f_0|^2|f_1|^{2s}|f_2|^{2t}e^{-2c\varphi}\text{ is integrable near }o\big\}.$$
It is clear that $A(s,0)=\Tn(s;f_0,f_1,\varphi)$ and $A(0,t)=\Tn(t;f_0,f_2,\varphi)$. Denote $a_i=\sigma(\log|f_i|,\varphi)$ for $i=1,2$. As $\log|f_i|\le a_i\varphi+O(1)$ near $o$, we have
$$\log|f_1f_2|\le (a_1+a_2)\varphi+O(1)$$ near $o$.
Note that $\Tn(t;f_0,f_1f_2,\varphi)$ is concave. Then it suffices to prove 
\[\Tn_-'(0;f_0,f_1f_2,\varphi)\le a_1+a_2.\]

As $A(s,t)$ is concave in $s$ for fixed $t$, we have
\begin{equation}
	\label{eq:0907a}
	\begin{split}
	\limsup_{x\rightarrow0+0}\frac{A(-x,-x)-A(0,-x)}{-x}&\le\limsup_{x\rightarrow0+0}\frac{A(-y,-x)-A(0,-x)}{-y}\\
    &=\frac{A(-y,0)-A(0,0)}{-y}
	\end{split}
\end{equation}
for any $y>0$.
As $(f_i,\varphi;|f_0|^2)$ satisfies condition (T) for $i=1,2$, it follows from \eqref{eq:0907a} that
\begin{equation*}
	\begin{split}
		\Tn_-'(0;f_0,f_1f_2,\varphi)&=\limsup_{x\rightarrow0+0}\frac{A(-x,-x)-A(0,0)}{-x}\\
		&\le\limsup_{x\rightarrow0+0}\frac{A(-x,-x)-A(0,-x)}{-x}+\limsup_{x\rightarrow0+0}\frac{A(0,-x)-A(0,-x)}{-x}\\
		&\le\lim_{y\rightarrow0+0}\frac{A(-y,0)-A(0,0)}{-y}+\Tn'(0;f_0,f_2,\varphi)\\
		&=a_1+a_2.
	\end{split}
\end{equation*}
Thus, Lemma \ref{l:2} holds.
\end{proof}

We give a maximal property for $\max\{a_1\log|z_1|,\ldots,a_r\log|z_r|\}$.

\begin{Lemma}\label{l:3}
	Let $\psi_1\coloneqq \max\{a_1\log|z_1|,\ldots,a_r\log|z_r|\}$, and $\psi_2$ be a plurisubharmonic function near $o$ such that $\psi_2\ge\psi_1+O(1)$ near the origin $o$, where $a_j>0$ for $j=1,\ldots,r$. Let $(f,o)\in\mathcal{O}_o^*(\coloneqq \mathcal{O}_o\setminus\{0\})$. If $c_o^f(\psi_1)=c_o^f(\psi_2)$, then we have $\psi_1=\psi_2+O(1)$ near $o$.
\end{Lemma}
\begin{proof}
Without loss of generality, assume that $\psi_2$ is a negative plurisubharmonic function on $\Delta_{\rho}^n$ $(\rho>0)$ and $c_o(\psi_1)=1$. Then
$$\sum_{j=1}^{r}\frac{2}{a_j}=1.$$
Example \ref{example} shows that $\Tn(t;1,f,\psi_1)=at+b$ for any $t\ge0$, where $a$ and $b$ are constant. Since $\psi_2\ge\psi_1+O(1)$, the Tian functions $\Tn(t;1,f,\psi_1)\le\Tn(t;1,f,\psi_2)$ for any $t\ge 0$. Note that $\Tn(1;1,f,\psi_1)=\Tn(1;1,f,\psi_2)$ and $\Tn(t;1,f,\psi_2)$ is concave in $t$. It follows that $\Tn(0;1,f,\psi_1)=\Tn(0;1,f,\psi_2)$, i.e.
$$c_o(\psi_1)=c_o(\psi_2)=1.$$

Set
$$G_1(t)\coloneqq \inf\left\{\int_{\{\tilde\psi_1<-t\}}|F|^2\colon F(o)=1\right\}$$
and
$$G_2(t)\coloneqq \inf\left\{\int_{\{\tilde\psi_2<-t\}}|F|^2\colon F(o)=1\right\},$$
where 
\[\tilde\psi_1=\max\big\{a_1\log|z_1|-a_1\log \rho,\ldots,a_r\log|z_r|-a_r\log \rho\big\}\]
and $\tilde\psi_2=\max\{\tilde\psi_1,\psi_2\}$. We have
$c_o(\tilde\psi_1)=c_o(\tilde\psi_2)=1$ and $\tilde\psi_1\le\tilde\psi_2$.
The strong openness property (\cite{GZopen-c}) shows that $e^{-\tilde\psi_1}$ and $e^{-\tilde\psi_2}$ are not integrable near $o$.  Thus, $G_i(-\log x)$ is concave with respect to $x$ for $i=1,2$ (cf. \cite{guan-effect}). Note that
$$G_1(-\log x)=\pi^n \rho^{2n}x^{\sum_{j=1}^{r}\frac{2}{a_j}}=\pi^n \rho^{2n}x$$ is linear, $G_1(0)=G_2(0)$, and $G_1(t)\ge G_2(t)$ for any $t\ge0$. Then we have $G_1(t)=G_2(t)$ for any $t\ge 0$. We obtain $$\lambda\big(\{\tilde\psi_1<-t\}\big)=G_1(t)=G_2(t)\le\lambda\big(\{\tilde\psi_2<-t\}\big)\le\lambda\big(\{\tilde\psi_1<-t\}\big),$$
for any $t\ge0$. Thus, we get $\tilde\psi_2\equiv\tilde\psi_1$ on $\Delta_{\rho}^n$, which implies $\psi_2=\psi_1+O(1)$ near $o$.
\end{proof}

Assume $b\coloneqq \Tn'_-(0;f_0,f,\varphi)>0$. We can construct tuples that satisfy the condition (T) by taking ``$\max$". 

\begin{Lemma}\label{l:4}
	$(1)$ $(f,\max\{\frac{\log|f|}{b},\varphi\};|f_0|^2)$ satisfies condition (T) and $c_o^{f_0}(\varphi)=c_o^{f_0}(\max\{\frac{\log|f|}{b},\varphi\})$. 
	
	$(2)$ If $(g,\varphi;|f_0|^2)$ satisfies condition (T) for another holomorphic function $g$, then $(g,\max\{\frac{\log|f|}{b},\varphi\};|f_0|^2)$ also satisfies condition (T).
\end{Lemma}
\begin{proof}
	Using Lemma \ref{l:1}, we have \[\Tn(t;f_0,f,\varphi)=\Tn\Big(t;f_0,f,\max\big\{\frac{\log|f|}{b},\varphi\big\}\Big)\]
	for any $t\le 0$, which implies
	$\Tn'_-(0;f_0,f,\max\{\frac{\log|f|}{b},\varphi\})=b$ and $$c_o^{f_0}(\varphi)=c_o^{f_0}\Big(\max\big\{\frac{\log|f|}{b},\varphi\big\}\Big).$$
	As $|f|e^{-b\max\{\frac{\log|f|}{b},\varphi\}}\le 1$ , by definition we have $\Tn'_-(t;f_0,f,\max\{\frac{\log|f|}{b},\varphi\})\ge b$  and $\Tn'_+(t;f_0,f,\max\{\frac{\log|f|}{b},\varphi\})\ge b$ for any $t\ge0$. Note that the function $\Tn(t;f_0,f,\max\{\frac{\log|f|}{b},\varphi\})$ is concave in $t$. We obtain that $(f,\max\{\frac{\log|f|}{b},\varphi\};|f_0|^2)$ satisfies condition (T).
	
	Note that 
	\[\Tn(0;f_0,f,\varphi)=\Tn\Big(0;f_0,f,\max\big\{\frac{\log|f|}{b},\varphi\big\}\Big)=\Tn\Big(0;f_0,g,\max\big\{\frac{\log|f|}{b},\varphi\big\}\Big).\]
	If $(g,\varphi;|f_0|^2)$ satisfies condition (T), since $\Tn(t;f_0,g,\max\{\frac{\log|f|}{b},\varphi\})$ is concave in $t$, then we have that $(g,\max\{\frac{\log|f|}{b},\varphi\};|f_0|^2)$ also satisfies condition (T).
\end{proof}

The following lemma will be used in the proof of Proposition \ref{p:quasi}.

\begin{Lemma}\label{l:5}
	Let $\varphi=c\log|z_1^{m_1}\ldots z_n^{m_n}|$ for some $c>0$, and $\psi$ be a plurisubharmonic function near $o\in\mathbb{C}^n$ s.t. $\psi\ge\varphi+O(1)$. Let $f$ be a holomorphic function near $o$. Assume that $c_o^{f}(\varphi)=c_o^{f}(\psi)=1$. Then there exists a plurisubharmonic function $\psi_1=\max\{a_1\log|z_{k_1}|,\ldots,a_r\log|z_{k_r}|\}$ ($1\le k_1<\ldots<k_r\le n$ and $a_j>0$) near $o$ s.t. $\psi_1\ge\psi+O(1)$ and $c_o^{f}(\psi_1)=1$.
\end{Lemma}
\begin{proof}
Firstly, we prove the case that $f=z^{\alpha}$ is monomial, where $\alpha=(\alpha_1,\ldots,\alpha_n)\in\mathbb{Z}_{\ge0}^n$.

Using Lemma \ref{l:4}, by finite steps, we obtain a plurisubharmonic function $\tilde\psi\coloneqq \max\big\{\psi,\frac{\log|z_{k_1}|}{b_1},\ldots,\frac{\log|z_{k_r}|}{b_r}\big\}$ ($1\le k_1<\ldots<k_r\le n$ and $b_j>0$) satisfying the following statements:
\begin{enumerate}
    \item $(z_{k_j},\tilde\psi;|f|^2)$ satisfies condition (T) for $j=1,\ldots,r$;
    \item $\Tn'_-(0;f,z_l,\tilde\psi)=0$ for any $l\not\in S$, which implies that $(z_l,\tilde\psi;|f|^2)$ satisfies condition (T), where $S\coloneqq \{k_1,\ldots,k_r\}$;
    \item $c_o^{f}(\tilde\psi)=c_o^{f}(\psi)=1$.
\end{enumerate}
Thus, we have $\sigma(z_{k_j},\tilde\psi)=b_{j}$ for $j=1,\ldots,r$ and $\sigma(z_l,\tilde\psi)=0$ for $l\not\in S.$ It follows from Lemma \ref{l:2} that
\begin{equation}
	\nonumber
	\begin{split}
	\sigma(\varphi,\tilde\psi)&=c\sigma(\log|z_1^{m_1}\ldots z_n^{m_n}|,\tilde\psi)\\
	&=c\sum_{l=1}^n m_l\sigma(z_l,\tilde\psi)=c\sum_{j=1}^rm_{k_r}b_j.
	\end{split}
\end{equation}
Since $\tilde\psi\ge\varphi+O(1)$ and $c_o^{f}(\tilde\psi)=c_o^{f}(\varphi)=1$, we have $\sigma(\varphi,\tilde\psi)=1$, which implies
\begin{equation}
	\label{eq:240909a}
	c\sum_{j=1}^rm_{k_j}b_j=1.
\end{equation}
By a direct calculation, $c_o^{f}(\varphi)=1$ implies that
\begin{equation}
	\label{eq:240909b}
	\alpha_l+1\ge cm_l, \quad l=1,\ldots,n.
\end{equation}
Let $\psi_1=\max\big\{\frac{\log|z_{k_1}|}{b_1},\ldots,\frac{\log|z_{k_r}|}{b_r}\big\}$. It follows from \eqref{eq:240909a} and \eqref{eq:240909b} that
$$1=c_o^{f}(\tilde\psi)\ge c_o^{f}(\psi_1)=\sum_{j=1}^r b_j(\alpha_{k_j}+1)\ge 	c\sum_{j=1}^rm_{k_j}b_j=1,$$
which shows that $c_o^{f}(\tilde\psi)= c_o^{f}(\psi_1)=1$. Following from $\tilde\psi\ge\psi_1+O(1)$ and Lemma \ref{l:3}, we have $\psi_1=\tilde\psi+O(1)\ge \psi+O(1)$ near $o$.
Thus, taking $a_j=\frac{1}{b_j}$ for $j=1,\ldots,r$, Lemma \ref{l:5} holds for this case.

Now, we prove the general case. Write $f=\sum_{\alpha}c_{\alpha}z^{\alpha}$ near $o$. Denote $W\coloneqq \{\alpha\colon c_{\alpha}\not=0\}$.  It is clear that $c_o^{z^{\alpha}}(\varphi)\ge 1$ for any $\alpha\in W$, which implies $c_o^{z^{\alpha}}(\psi)\ge 1$ for any $\alpha\in W$ by $\psi\ge\varphi+O(1)$ near $o$.
As $c_o^{f}(\psi)=1$, there exists $\alpha_0\in W$ such that $c_o^{z^{\alpha_0}}(\psi)=1.$
 According to the previous arguments, we obtain a plurisubharmonic function $\psi_1=\max\{a_1\log|z_{k_1}|,\ldots,a_r\log|z_{k_r}|\}$  s.t. $\psi_1\ge\psi+O(1)$ and $c_o^{z^{\alpha_0}}(\psi_1)=1$. Note that 
 \[\int_{B}|f|^2e^{-2c\psi_1}\ge\int_B|c_{\alpha_0}z^{\alpha_0}|^2e^{-2c\psi_1}\]
 for any small ball $B$ centered at $o$ and $c>0$. Then we have
$c_o^{f}(\psi_1)\le 1$. As $\psi_1\ge\psi+O(1)$ and $c_o^f(\psi)=1$, we have $c_o^{f}(\psi_1)=1$. Thus, Lemma \ref{l:5} holds.
\end{proof}

\subsection{Zhou weight}

In this section, we recall some results on local Zhou weights in \cite{BGMY-valuation}, whose definition can be seen in Definition \ref{def:max_relat}.

\begin{Theorem}[{see \cite[Theorem 1.11]{BGMY-valuation}}]
	\label{thm:multi-valua}
	Let $u,$ $v$ be two plurisubharmonic functions near $o$. Then the following  statements are equivalent:
	
	$(1)$ there exists a plurisubharmonic function $\varphi_{0}$ near $o$ ($\not\equiv-\infty$) and two sequences of numbers  $\{t_{i,j}\}_{j\in\mathbb{Z}_{\ge0}}$ ($t_{i,j}\rightarrow+\infty$ when $j\rightarrow+\infty$, $i=1,2$) such that $\lim_{j\rightarrow+\infty}\frac{t_{1,j}}{t_{2,j}}=1$ and $$\mathcal{I}(\varphi_{0}+t_{1,j}v)_o\subset\mathcal{I}(\varphi_{0}+t_{2,j}u)_o$$ for any $j$;
	
	$(2)$ for any plurisubharmonic function $\varphi_0$ near $o$ and any $t>0$, we have $$\mathcal{I}(\varphi_0+tv)_o\subset\mathcal{I}(\varphi_0+tu)_o;$$

	$(3)$ for any local Zhou weight $\Phi_{o,\max}$ near $o$, we have $$\sigma(u,\Phi_{o,\max})\le\sigma(v,\Phi_{o,\max}).$$
\end{Theorem}

Let $f_{0}=(f_{0,1},\cdots,f_{0,m})$ be a vector of holomorphic functions near $o$. Let $\varphi_{0}$ be a plurisubharmonic function near $o$,
such that $|f_{0}|^{2}e^{-2\varphi_{0}}$ is integrable near $o$. 
Recall the existence of local Zhou weight: 

\begin{Remark}[{see \cite[Remark 1.3]{BGMY-valuation}}]
	\label{rem:max_existence}Let $\varphi$ be a plurisubharmonic function near $o$. 
	Assume that $|f_{0}|^{2}e^{-2\varphi_{0}}|z|^{2N_{0}}e^{-2\varphi}$ is integrable near $o$
	for large enough $N_{0}\gg0$,
	and $|f_0|^2e^{-2\varphi-2\varphi_0}$ is not integrable near $o$.
	
	Then there exists a local Zhou weight $\Phi_{o,\max}$ related to $|f_{0}|^{2}e^{-2\varphi_{0}}$ near $o$
	such that $\Phi_{o,\max}\geq\varphi$.
\end{Remark}

Let $\Phi_{o,\max}$ be a local Zhou weight  related to $|f_0|^2e^{-2\varphi_{0}}$ near $o$. 
\begin{Lemma}[{see \cite[Remark 1.4]{BGMY-valuation}}]\label{l:zhounumber-semiconti}
    If plurisubharmonic functions $\psi_{j}\to\psi$ in $L^{1}_{\mathrm{loc}}$ when $j\rightarrow+\infty$,
	then
	\begin{equation}
\nonumber		\limsup_{j\rightarrow+\infty} \sigma(\psi_{j},\Phi_{o,\max})\leq \sigma(\psi,\Phi_{o,\max}).\end{equation}
\end{Lemma}

For any plurisubharmonic function $\psi$ near $o$, recall the Tian function
\[\Tn (t)\coloneqq \sup\big\{c\ge 0 \colon |f_0|^2e^{-2\varphi_{0}}e^{2t\psi}e^{-2c\Phi_{o,\max}} \ \text{is integrable near} \ o\big\}.\]
Then we have

\begin{Proposition}[{see \cite[Proposition 3.6]{BGMY-valuation}}]
	\label{prop:concave_B}
	The Tian function $\Tn (t)$ is differentiable at $t=0$, and
	\begin{equation*}
		\Tn (t)=\Tn (0)+\sigma(\psi,\Phi_{o,\max})t
	\end{equation*}
	holds for any $t\ge0$. 
\end{Proposition}


 We recall the following Theorem on the relation between the jumping number $c_o^{G}(\Phi_{o,\max})$ and the Zhou valuation $\nu(G,\Phi_{o,\max})$ for any $(G,o)\in\mathcal{O}_o$.

\begin{Theorem}[{see \cite[Theorem 1.10]{BGMY-valuation}}]
	\label{thm:valu-jump}
	For any holomorphic function $G$ near $o$,	
	we have the following relation between the jumping number $c^G_o(\Phi_{o,\max})$ and the Zhou valuation $\nu(\cdot,\Phi_{o,\max})$,
	\begin{equation*}
		\begin{split}
			\nu(G,\Phi_{o,\max})+c_o(\Phi_{o,\max})
			&\le c^G_o(\Phi_{o,\max})\\
			&\le \nu(G,\Phi_{o,\max})-\sigma(\log|f_0|,\Phi_{o,\max})+1+\sigma(\varphi_0,\Phi_{o,\max}).
		\end{split}
	\end{equation*}
\end{Theorem}

\subsection{Tame maximal weight}

In this section, we discuss the tame maximal weights near $o$. The definitions of ``tame weight" and ``maximal weight" can be seen in Introduction.

The following lemma is used to determine whether a plurisubharmonic weight is maximal.

\begin{Lemma}[see \cite{Blo-note}; see also \cite{BT76,Blo93}]
	\label{l:max1}Let $\varphi\in \mathrm{PSH}(\Omega)\cap L^{\infty}_{\mathrm{loc}}(\Omega)$ on an open subset $\Omega$ of $\mathbb{C}^n$. If for any $u\in \mathrm{PSH}(\Omega)$ such that $\varphi\ge u$ outside a compact subset of $\Omega$ we have $\varphi\ge u$ on $\Omega$, then $(dd^c\varphi)^n=0$ on $\Omega$.
\end{Lemma}

The ``tame" property will imply an equality about the Tian  function.

\begin{Remark}\label{r:tame}Let $\varphi$ be a plurisubharmonic function near $o$ satisfying that the Lelong number $\nu_o(\varphi)>0$ (which implies $c_o^{f}(\varphi)<+\infty$ for any $(f,o)\in\mathcal{O}_{o}^*$) and
\[\sup_{(f,o)\in\mathcal{O}_{o}^*}\big(c_o^{f}(\varphi)-\sigma(\log|f|,\varphi)\big)<+\infty.\]
Then for any $(f,o)\in\mathcal{O}_{o}^*$, considering the Tian function 
\[\Tn(t)\coloneqq \sup\big\{c\in\mathbb{R}\colon |f|^{2t}e^{-2c\varphi} \text{ is integrable near } o\big\},\]
we have 
$$\lim_{t\rightarrow+\infty}\frac{\Tn(t)}{t}=\Tn'_-(t)=\Tn'_+(t)=\sigma(\log|f|,\varphi),$$
where $\Tn'_-(t)$ and $\Tn'_+(t)$ are the left and right derivatives of function $\Tn(t)$, respectively.

\begin{proof}
By the definitions of $\Tn(t)$ and $\sigma(\log|f|,\varphi)$, we have 
\begin{equation}
    \nonumber\begin{split}
    \Tn(t+\Delta t)&=\sup\{c:|f|^{2t+\Delta t}e^{-2c\varphi} \text{ is integrable near } o\}\\
    &\ge\sup\{c:|f|^{2t}e^{\Delta t(\sigma(\log|f|,\varphi)-\epsilon)\varphi}e^{-2c\varphi} \text{ is integrable near } o\}\\
    &=\Tn(t)+\Delta t(\sigma(\log|f|,\varphi)-\epsilon)
    \end{split}
\end{equation}
for any $\epsilon>0$ and $\Delta t>0$, which shows
$\Tn'_+(t)\ge \sigma(\log|f|,\varphi)$. It follows from the concavity of $\Tn(t)$ that the functions $\frac{\Tn(t)}{t}$, $\Tn'_-(t)$ and $\Tn'_+(t)$ are decreasing, and $$\lim_{t\rightarrow+\infty}\frac{\Tn(t)}{t}=\lim_{t\rightarrow+\infty}\Tn'_-(t)=\lim_{t\rightarrow+\infty}\Tn'_+(t).$$
As $\sup_{(f,o)\in\mathcal{O}_{o}^*}\big(c_o^{f}(\varphi)-\sigma(\log|f|,\varphi)\big)<+\infty$, there exists $C>0$ such that 
\[\Tn(m)=c_o^{f^m}(\varphi)\le m\sigma(\log|f|,\varphi)+C, \quad \forall\,m\in\mathbb{Z}_{\ge 1},\]
which implies 
\[\lim_{m\rightarrow+\infty}\frac{\Tn(m)}{m}\le \sigma(\log|f|,\varphi).\]
As $\Tn(t)/t\ge \sigma(\log|f|,\varphi)$ for any $t\ge0$ and $\lim_{t\rightarrow+\infty}\Tn(t)/t$ exists, we have 
$\lim_{t\rightarrow+\infty}\Tn(t)/t=\sigma(\log|f|,\varphi)$.
\end{proof}
\end{Remark}

Let us recall a useful sufficient condition for a weight being tame.

\begin{Lemma}[see \cite{BFJ08}]
	\label{l:holder}
	Let $\varphi$ be a plurisubharmonic function near $o$ such that $e^{\varphi}$ is $\alpha$-H\"older for some $\alpha>0$. Then for any $(f,o)\in\mathcal{O}_o^*$, we have 
	\[c_o^{f}(\varphi)\le\sigma(\log|f|,\varphi)+\frac{n}{\alpha}.\]
\end{Lemma}
\begin{proof}
	For the convenience of readers, we recall the proof.
	If $c_o^{f}(\varphi)>a>0$, then $(f,o)\in\mathcal{I}(a\varphi)_o$.
	Assume that $\varphi$ is defined on a ball $B$ and $\mathcal{I}(a\varphi)_o$ is generated by $\{f_1,\ldots,f_l\}$ such that $\int_B\sum_{1\le j\le l}|f_j|^2e^{-2a\varphi}<+\infty$.
	
	By Lemma \ref{l:appro-Berg}, we have $$\varphi_a(p)\le\sup_{|z-p|<r}\varphi(z)+\frac{1}{a}\log\frac{c_2}{r^n},$$ where $\varphi_a$ is defined as in Lemma \ref{l:appro-Berg}. As $e^{\varphi}$ is $\alpha$-H\"older, we have
	$$\sup_{|z-p|<r}e^{\varphi( z)}\le e^{\varphi(p)}+Cr^{\alpha}.$$
	Choosing $r=e^{\frac{\varphi(p)}{\alpha}}$, we get $\sup_{|z-p|<r}\varphi(z)\le\varphi(p)+O(1)$, which implies
	$$\varphi_a\le \Big(1-\frac{n}{\alpha a}\Big)\varphi+O(1)$$
	near $o$. As $\int_B\sum_{1\le j\le l}|f_j|^2e^{-2a\varphi}<+\infty$ and $(f,o)\in\mathcal{I}(a\varphi)_o$, we have 
	\[\log|f|\le \log\Big(\sum_{1\le j\le l}|f_j|\Big)+O(1) \le a\varphi_a+O(1)\]
	near $o$, which shows
	$$\log|f|\le \Big(a-\frac{n}{\alpha}\Big)\varphi+O(1).$$
	Thus, Lemma \ref{l:holder} holds.
\end{proof}

Let $\varphi_1$ be a plurisubharmonic function  on a bounded domain $D\subset\mathbb{C}^n$ satisfying that:
\begin{enumerate}
    \item $\sup_{(f,o)}\big(c_o^{f}(\varphi_1)-\sigma(\log|f|,\varphi_1)\big)<+\infty$; and
    \item $C_1\log|z|+O(1)\le\varphi_1\le C_2\log|z|+O(1)$ near $o$ for some positive constants $C_1,C_2$.
\end{enumerate}
Denote
\begin{equation*}
        \varphi_2\coloneqq \Big(\sup\big\{\varphi\in \mathrm{PSH}^-(D)\colon\sigma(\varphi,\Phi)\ge\sigma(\varphi_1,\Phi),  \ \forall\text{ local Zhou weight $\Phi$ near $o$}\big\}\Big)^*,
\end{equation*}
where $\mathrm{PSH}^-(D)$ denotes all negative plurisubharmonic functions on $D$. Then $\varphi_2\in \mathrm{PSH}^-(D)$ (see Proposition \ref{pro:Demailly}). In the following, we will prove that $\varphi_2$ is a tame maximal weight.

\begin{Lemma}\label{l:relative type}
	$\varphi_2$ is a maximal  weight, and for any negative plurisubharmonic function $\psi$ on $D$, we have $\psi\le\sigma(\psi,\varphi_2)\varphi_2$
	on $D$.
\end{Lemma}
\begin{proof}
	As $\psi\le p\varphi_2+O(1)$ near $o$, for any $p\in \big(0,\sigma(\psi,\varphi_2)\big)$ and  any local Zhou weight $\Phi$ near $o$,
	$$\sigma\Big(\max\{\varphi_2,\frac{\psi}{p}\},\Phi\Big)=\min\Big\{\sigma(\varphi_2,\Phi),\sigma(\frac{\psi}{p},\Phi)\Big\}=\sigma(\varphi_2,\Phi).$$
	By definition of $\varphi_2,$ we have $\max\{\varphi_2,\frac{\psi}{p}\}\le\varphi_2$ on $D$, which shows $\psi\le p\varphi_2$
	on $D$. Taking $p\rightarrow \sigma(\psi,\varphi_2)$, $\psi\le\sigma(\psi,\varphi_2)\varphi_2$ holds.
	
	It is clear that  $\sigma(\log|z|,\varphi_2)\ge\sigma(\log|z|,\varphi_1)>0$. It follows from   $\log|z|+C\le\sigma(\log|z|,\varphi_2)\varphi_2$ on $D$ for some $C$ that $\varphi_2$ is locally bounded on  $D\setminus\{o\}$. Hence,  $\varphi_2$ is a maximal weight by Lemma \ref{l:max1} and the definition of $\varphi_2$.
\end{proof}

The jumping numbers (and relative types) for $\varphi_1$ and $\varphi_2$ coincide:
\begin{Lemma}
	\label{l:tame-maximal}
	$c_o^{f}(\varphi_1)=c_o^{f}(\varphi_2)$ and $\sigma(\log|f|,\varphi_1)=\sigma(\log|f|,\varphi_2)$ for all $(f,o)\in\mathcal{O}_o$.
\end{Lemma}

\begin{proof}Note that $\sigma(\max\{u_1,u_2\},\Phi)=\min\{\sigma(u_1,\Phi),\sigma(u_2,\Phi)\}$.
By the definition of $\varphi_2$, Lemma \ref{lem:Choquet} and Proposition \ref{pro:Demailly}, there exists a sequence of negative plurisubharmonic functions $\{\tilde\varphi_j\}$ on $D$, which satisfies $\sigma(\tilde\varphi_j,\Phi)\ge\sigma(\varphi_1,\Phi)$ for any $j$ and any local Zhou weight $\Phi$ near $o$, and $\{\tilde\varphi_j\}$ increasingly converges to $\varphi_2$ almost everywhere on $D$. Hence, it follows from Lemma \ref{l:zhounumber-semiconti} that $$\sigma(\varphi_2,\Phi)\ge\sigma(\varphi_1,\Phi).$$
Note that $\varphi_2\ge \varphi_1$. Then $\sigma(\varphi_2,\Phi)=\sigma(\varphi_1,\Phi)$ for any local Zhou weight $\Phi$. Combining with Theorem \ref{thm:multi-valua}, we get
\[c_o^{f}(\varphi_1)=c_o^f(\varphi_2), \quad \forall\,(f,o)\in\mathcal{O}_o^*.\]
	
For every nonzero holomorphic germ $(f,o)$ and every $k\in\mathbb{Z}_{\ge 1}$, as the weight $\varphi_1$ satisfies $\sup_{(f,o)}(c_o^{f}(\varphi_1)-\sigma(\log|f|,\varphi_1))<+\infty$, there exists $M>0$ (independent of $f$ and $k$) such that
$$c_o^{f^k}(\varphi_2)=c_o^{f^k}(\varphi_1)\le k\sigma(\log|f|,\varphi_1)+M.$$
It is clear that $c_o^{f^k}(\varphi_2)\ge k\sigma(\log|f|,\varphi_2)$,
which implies
$$\sigma(\log|f|,\varphi_2)\le\frac{c_o^{f^k}(\varphi_2)}{k}\le\frac{ k\sigma(\log|f|,\varphi_1)+M}{k}.$$
Letting $k\rightarrow+\infty$, we have $\sigma(\log|f|,\varphi_2)\le \sigma(\log|f|,\varphi_1)$. Note that $\varphi_2\ge \varphi_1+O(1)$ near $o$. Then we have $\sigma(\log|f|,\varphi_2)= \sigma(\log|f|,\varphi_1)$ for any $(f,o)\in\mathcal{O}_o$.
\end{proof}

For any $t>0$, denote
\begin{equation*}
        \varphi_{2,t}\coloneqq \Big(\sup\big\{\varphi\in \mathrm{PSH}^-(D_t)\colon\sigma(\varphi,\Phi)\ge\sigma(\varphi_1,\Phi),  \ \forall\text{ local Zhou weight $\Phi$ near $o$}\big\}\Big)^*,
\end{equation*}
where $D_t\coloneqq \{z\in D\colon\varphi_2<-t\}$. Similarly, we know that $\varphi_{2,t}$ is a maximal  weight such that $c_o^{f}(\varphi_1)=c_o^{f}(\varphi_{2,t})$ and $\sigma(\log|f|,\varphi_1)=\sigma(\log|f|,\varphi_{2,t})$ for any $(f,o)\in\mathcal{O}_o$.
\begin{Lemma}\label{l:varphi_2+t}
	For any $t>0$,  $\varphi_{2,t}=\varphi_2+t$ on $D_t$.
\end{Lemma}

\begin{proof}
	By the definition of $\varphi_{2,t}$, we have
	$$\varphi_{2,t}\ge\varphi_2+t$$
	on $D_t$.
	Let
	\begin{equation*}
		\phi\coloneqq \left\{
		\begin{array}{ll}
			\varphi_{2,t}-t & \ \text{on} \ D_t, \\
			\varphi_2 & \ \text{on} \ D\setminus D_t.
		\end{array}
		\right.
	\end{equation*}
	Then we have $\phi\ge\varphi_2$ on $D$. As $\varphi_{2,t}$ is plurisubharmonic on $D_t$ and $\varphi_2$ is plurisubharmonic on $D$, we have
	\[\limsup_{\tilde z\rightarrow z} \phi(\tilde{z})\le\phi(z), \quad  \forall\, z\in D\setminus\partial D_t.\]
	For any $z\in \partial D_t\cap D$, we have $\phi(z)=\varphi_2(z)\ge-t$, which shows
	\begin{flalign*}
		\begin{split}
			\limsup_{\tilde z\rightarrow z}\phi(\tilde{z})&=\max\left\{\limsup_{D_t\ni\tilde z\rightarrow z}\phi(\tilde{z}),\limsup_{(D\setminus D_t)\ni\tilde z\rightarrow z}\phi(\tilde{z})\right\}\\
			&\le\max\{-t, \varphi_2(z)\}\le\phi(z).
		\end{split}
	\end{flalign*}
	Thus, $\phi$ is an upper semicontinuous function on $D$. For any $z\in \partial D_t\cap D$, as $\phi\ge\varphi_2$ on $D$, we have
	$$\phi(z)=\varphi_2(z)\le\frac{1}{\lambda(\mathbb{B}(z,r))}\int_{B(z;r)}\varphi_2 \ d\lambda\le\frac{1}{\lambda(\mathbb{B}(z,r))}\int_{B(z;r)}\phi \ d\lambda,$$
	where $\mathbb{B}(z,r)\subset D$ is a disk of dimension $1$ with center $z$ and a radius of $r$ and $\lambda$ is the Lebesgue measure on $\mathbb{C}$. Hence, $\phi$ is a negative plurisubharmonic function on $D$.
	By definition of $\varphi_2$, we have $\varphi_2\ge \phi$ on $D$, which implies
	$$\varphi_{2,t}=\varphi_2+t$$
	on $D_t$.
\end{proof}

Using Lemma \ref{l:varphi_2+t} and Demailly’s approximation theorem (Lemma \ref{l:appro-Berg}), we prove the continuity of $e^{\varphi_2}.$ Thus, $\varphi_2$ is a tame maximal weight.

\begin{Lemma}\label{l:tame-maximal2}
	Assume that $D$ is a hyperconvex domain. Then $e^{\varphi_2}$ is continuous on $D$ such that $e^{\varphi_2}\rightarrow1$ when $z\rightarrow\partial D$. Especially, $\varphi_2$ is a tame weight.
\end{Lemma}
\begin{proof}
	Note that there exists $N_0\gg0$ such that $\varphi_2\ge N_0\log|z|+O(1)$ near $o$. As $D$ is a hyperconvex domain, there exists a continuous exhaustion plurisubharmonic function $\rho<0$ on $D$. Then there exist  $r_1>r_2>0$ and $N_1>0$ such that $\{r_2<|z|<r_1\}\Subset D$ and $$\inf_{\{r_2<|z|<r_1\}}\varphi_2(z)\ge  N_1\sup_{\{r_2<|z|<r_1\}}\rho(z).$$
	Let
	\begin{equation*}
		\tilde\varphi_2\coloneqq \left\{
		\begin{array}{ll}
			\varphi_2 & \ \text{on} \ \{|z|<r_1\}, \\
			\max\{\varphi_2,N_1\rho\} & \ \text{on} \  D\setminus\{|z|<r_1\}.
		\end{array}
		\right.
	\end{equation*}
	Then $\tilde\varphi_2$ is a negative plurisubharmonic function on $D$ and $\tilde\varphi_2\ge \varphi_2$. By the definition of $\varphi_2$, we have $\tilde\varphi_2= \varphi_2$, i.e. $\varphi_2\ge N_1\rho$ on $D\setminus\{|z|<r_1\}$, which implies that $\varphi_2\rightarrow0$ when $z\rightarrow\partial D$.
	
	For any positive integer $m$, let $\{\sigma_{m,k}\}_{k=1}^{\infty}$ be an orthonormal basis of
	\[A^2\big(D,2m\varphi_2\big)\coloneqq \left\{f\in\mathcal{O}(D):\int_{D}|f|^2e^{-2m\varphi_2}d\lambda<+\infty\right\}.\]
	Denote by
	$$\psi_m\coloneqq \frac{1}{2m}\log\sum_{k=1}^{\infty}|\sigma_{m,k}|^2$$
	on $D$.
	
	Fixed any $z\in D$, there exists a holomorphic function $f_{z,m}$ on $D$ such that
	\[\int_{D}|f_{z,m}|^2e^{-2m\varphi_2}d\lambda=1 \ \&  \ \frac{1}{2m}\log|f_{z,m}(z)|^2=\psi_m(z).\]
	For any $t>0$ satisfying $z\in\{\varphi_2<-t\}$, we have
	\begin{equation*}
		\begin{split}
			1=\int_{D}|f_{z,m}|^2e^{-2m\varphi_2}d\lambda&
			\ge \int_{\{\varphi_2<-t\}}|f_{z,m}|^2e^{-2m\varphi_2}d\lambda\\
			&\ge e^{2mt}\int_{\{\varphi_2<-t\}}|f_{z,m}|^2d\lambda.
		\end{split}
	\end{equation*}
	It follows from $\varphi_2(z)\rightarrow0$ when $z\rightarrow\partial D$ and $\int_{D}|f_{z,m}|^2e^{-2m\varphi_2}d\lambda=1$ that, for $t>0$ small enough, $\int_{\{\varphi_2<-t\}}|f_{z,m}|^2d\lambda\le e^{-2mt}$ holds. By the construction of $\varphi_2$, we know that there exist positive numbers $b\ge a$ such that $a\log|z|\le \varphi_2\le b\log|z|$ on $D$. Let $0<t_1<t_2$ be small enough. Note that $\{\varphi_2<-t_2\}\subset \{a\log|z|<-t_2\}\subset \{b\log|z|<-t_1\}\subset \{\varphi<-t_1\}$. One has, for any point $w\in\{\varphi_2<-t_2\}$, there exists an positive number $r>0$ (independent of $w$), such that $B(w,r)\subset \{\varphi<-t_1\}$. It follows from the subharmonicity of $|f_{z,m}|^2$ that, for any $w\in\{\varphi_2<-t_2\}$,
	\begin{equation}
		\label{eq:estimate for module}
		\begin{split}
			|f_{z,m}(w)|^2 &\le \frac{1}{\sigma_n r^{2n}}\int_{B(w,r)}|f_{z,m}(w)|^2d\lambda\\
			&\le
			\frac{1}{\sigma_n r^{2n}}\int_{\{\varphi_2<-t_1\}}|f_{z,m}|^2d\lambda\\
			&\le \frac{e^{-2mt_1}}{\sigma_n r^{2n}},
		\end{split}
	\end{equation}
	where $\sigma_n=\frac{\pi^n}{n!}$.	
	By \eqref{eq:estimate for module}, we can find  $M\gg0$ such that
	$|f_{z,m}|<1$ on $\{\varphi_2<-t\}$ for any $m>M$. It follows from Lemma \ref{l:relative type} and Lemma \ref{l:varphi_2+t}  that
	$$\psi_m(z)=\frac{1}{m}\log|f_{z,m}(z)|\le\frac{\sigma\big(\log|f_{z,m}|,\varphi_2+t\big)}{m}(\varphi_2+t)(z)$$
	for $z\in\{\varphi_2<-t\}$. Note that $c_o^{f_{z,m}}(\varphi_2)\ge m$, then Lemma \ref{l:tame-maximal} shows that $\sigma\big(\log|f_{z,m}|,\varphi_2\big)\ge m+M$, where $M$ is a constant independent of $m$. Hence, we have
	\begin{equation}
		\label{eq:0911a}
		\frac{\psi_m(z)}{1+\frac{M}{m}}-t\le \varphi_2(z).
	\end{equation}
	Using inequality \eqref{eq:0911b} in Lemma \ref{l:appro-Berg}, we get
	\begin{equation}
		\label{eq:0911c}
		\varphi_2\le \psi_m+\frac{c_1}{m}
	\end{equation}
	on $D$. Since $e^{\psi_m}$ is smooth on $D$, letting  $m\to+\infty$ firstly and then $t\to0$ in the inequalities \eqref{eq:0911a} and \eqref{eq:0911c}, we get that $e^{\varphi_2}$ is continuous on $D$. Combining with Lemma  \ref{l:tame-maximal}, $\varphi_2$ is a tame weight.
\end{proof}

\subsection{Convergence results for valuations and for relative types}

\

Let $\{\nu_{j}\}_{j\in\mathbb{Z}_{>0}}$ be a sequence of valuations on $\mathcal{O}_o$ satisfying the following conditions:

\begin{enumerate}
    \item There exists a sequence of tame maximal weights $\{\Phi_j\}_{j\in\mathbb{Z}_{>0}}$ near $o$ such that $\nu_j(f)=\sigma(\log|f|,\Phi_j)$ for any $(f,o)\in\mathcal{O}_o$;
    \item There exists $A_1>0$ such that $A_1\log|z|\le\Phi_{j}+O(1)$ near $o$ holds for all $j$;
    \item $\sup_{j}\sup_{(f,o)\in\mathcal{O}_o\setminus\{0\}}\big(c_o^f(\Phi_{j})-\nu_j(f)\big)<+\infty$.
\end{enumerate}

We provide a necessary and sufficient condition for statement (1) in Theorem \ref{thm:valuationexistsweight}.

\begin{Remark}
  Condition (3) implies that $\sup_j\nu_j(f)<+\infty$ for any $(f,o)\in \mathcal{O}_o^*$. In fact, by definition, we have $c_o(\Phi_{j})\le c_o^f(\Phi_{j})-\sigma(\log|f|,\Phi_j)\le M$ for every $j$, which implies the Lelong number $\nu_o(\Phi_{j})\ge \frac{1}{M}$ by Lemma \ref{l:Lelong}, i.e., $\Phi_j\le \frac{1}{M}\log|z|+O(1)$ near $o$. Then $\nu_j(f)=\sigma(\log|f|,\Phi_{j})\le M\nu_o(\log|f|)$ for any $(f,o)\in \mathcal{O}_o^*$ and every $j$.
\end{Remark}

We give a convergence result for valuations, which will be used to solve the valuative interpolation problem.

\begin{Proposition}\label{thm:converge}
	There exist a subsequence of $\{\nu_j\}$ denoted by $\{\nu_{j_l}\}$ and a tame maximal weight $\varphi_{0}$ near $o$, which satisfy that $\{\nu_{j_l}\}$ converges to a valuation $\nu_0$ on $\mathcal{O}_o$ and $\nu_0(f)=\sigma(\log|f|,\varphi_{0})$ for any $(f,o)\in\mathcal{O}_o$.
\end{Proposition}
\begin{proof}
	Denote $I_k^j\coloneqq \{(f,o)\in\mathcal{O}_o\colon \nu_j(f)\ge k\}$ for every $j,k\in\mathbb{Z}_{>0}$. We prove Proposition \ref{thm:converge} in three steps.
	
	\
	
	\emph{Step 1.} Find a subsequence $\{I_k^{j_l}\}$ of $I_k^j$ such that $\cup_{l_1\ge1}\big(\cap_{l\ge l_1}I_k^{j_l}\big)$ is maximal for any $k$. It means that for any subsequence $\{I_k^{j'_l}\}$ of $\{I_k^{j_l}\}$, $\cup_{l_1\ge1}\big(\cap_{l\ge l_1}I_k^{j'_l}\big)=\cup_{l_1\ge1}\big(\cap_{l\ge l_1}I_k^{j_l}\big)$.
	
	Note that
	$$\cup_{j_1\ge1}\big(\cap_{j\ge j_1}I_k^{j}\big)=\{(f,o)\in\mathcal{O}_o\colon\exists\, j_1\text{ s.t. $\nu_j(f)\ge k$ holds for any $j\ge j_1$}\}.$$
	
	As $\mathcal{O}_o$ is  a Noetherian ring, there exists a subsequence $\{I_1^{j^{(1)}(l)}\}$ of $I_1^j$ such that $\cup_{l_1\ge1}\big(\cap_{l\ge l_1}I_1^{j^{(1)}(l)}\big)$ is maximal. Similarly, there exists a subsequence $\{I_1^{j^{(2)}(l)}\}$ of $I_1^{j^{(1)}(l)}$ such that $\cup_{l_1\ge1}\big(\cap_{l\ge l_1}I_2^{j^{(2)}(l)}\big)$ is maximal. By induction, for any $r\in\mathbb{Z}_{\ge1}$, there exists a sequence $\{I_1^{j^{(r)}(l)}\}$  such that $\cup_{l_1\ge1}\big(\cap_{l\ge l_1}I_r^{j^{(r)}(l)}\big)$ is maximal and $\{I_1^{j^{(r+1)}(l)}\}$ is a subsequence of $\{I_1^{j^{(r)}(l)}\}$. Using the diagonal method, we obtain a subsequence $\{I_k^{j_l}\}$ of $I_k^j$ such that $\cup_{l_1\ge1}\big(\cap_{l\ge l_1}I_k^{j_l}\big)$ is maximal for any $k$.

	\

	\emph{Step 2.}
	We prove 
	\[\lim_{l\rightarrow+\infty}\nu_{j_l}(f)=\nu_0(f)=\sigma(\log|f|,\varphi_{0}), \quad \forall\, (f,o)\in\mathcal{O}_o.\]
	
	Firstly, we prove that $\lim_{l\rightarrow+\infty}\nu_{j_l}(f)$ exists for any $(f,o)\in\mathcal{O}_o$.
	
	We prove it by contradiction: if not, there exist two subsequences of valuations $\{\nu_{j'_{l}}\}$ and $\{\nu_{j''_l}\}$ of $\{\nu_{j_l}\}$ such that $$\lim_{l\rightarrow+\infty}\nu_{j'_l}(f)=a_1<a_2=\lim_{l\rightarrow+\infty}\nu_{j''_l}(f).$$
	Note that $\nu_j(f^m)=m\nu_j(f)$. Then, without loss of generality, assume that there exists a positive integer $k_0\in (a_1,a_2)$. Then we have
	$$(f,o)\in \cup_{l_1\ge1}\big(\cap_{l\ge l_1}I_{k_0}^{j'_l}\big)$$
	and $$(f,o)\not\in \cup_{l_1\ge1}\big(\cap_{l\ge l_1}I_{k_0}^{j_l}\big),$$
	which contradicts to the maximality of $\cup_{l_1\ge1}\big(\cap_{l\ge l_1}I_{k_0}^{j_l}\big)$.
	Thus, $\lim_{l\rightarrow+\infty}\nu_{j_l}(f)$ exists for any $(f,o)\in\mathcal{O}_o$.

	Secondly, we construct the plurisubharmonic function $\varphi_{0}$.
	
	As $A_1\log|z|\le\Phi_{j}+O(1)$ near $o$ for any $j$, we have $\nu_j(f)\ge\frac{N}{A_1}$ for any $(f,o)\in J_N$, where $J_N$ is an ideal of $\mathcal{O}_o$ is generated by $\{z^{\alpha}\colon|\alpha|=N\}$. Then $$J_{\lceil A_1k\rceil}\subset \cap_j I_k^j,$$ which implies that $\mathcal{O}_o/\big(\cup_{j_1\ge1}\big(\cap_{l\ge j_1}I_k^{j}\big)\big)$ is a finite dimensional $\mathbb{C}$-linear space for any $k.$
	
	Denote $I_k\coloneqq \cup_{l_1\ge1}\big(\cap_{l\ge l_1}I_k^{j_l}\big)$. It is clear that $I_{k_1}\cdot I_{k_2}\subset I_{k_1k_2}$ for any $k_1,k_2\in\mathbb{Z}_{>0}$. As $J_{\lceil A_1k\rceil}\subset I_k,$ there exist finitely many polynomials $f_{k,1},\ldots,f_{k,r_k}$ which generate $I_k$. Denote $|I_k|^2\coloneqq \sum_{1\le r\le r_k}|f_{k,r}|^2$ (without loss of generality, we assume $|I_k|<1$ in the polydisc $\Delta^n$).
	
	Let $\varphi_0\coloneqq (\max_{k}\frac{1}{k}\log|I_k|)^*$ on $\Delta^n$, where $(h)^*$ denotes the upper semicontinuous regularization of $h$. Then $\varphi_0$ is a plurisubharmonic function on $\Delta^n$ and equals almost everywhere to $\max_{k}\frac{1}{k}\log|I_k|$.
	
	Thirdly, we prove \[\sigma(\log|f|,\varphi_0)=\lim_{l\rightarrow+\infty}\nu_{j_l}(f), \quad \forall\,(f,o)\in\mathcal{O}_o.\]
	
	If $\lim_{l\rightarrow+\infty}\nu_{j_l}(f)>a$ for some $a>0$, then there exist $s,k\in\mathbb{Z}_{>0}$ such that $\lim_{l\rightarrow+\infty}\nu_{j_l}(f^s)>k>sa$,
	which implies $(f^s,o)\in I_k$. By the definition of $\varphi_{0}$, we have
	$$\frac{s}{k}\log|f|\le \frac{1}{k}\log|I_k|+O(1)\le\varphi_{0}+O(1)$$
	near $o$. Then we obtain the $\sigma(\log|f|,\varphi_0)\ge\frac{k}{s}>a$, which shows $\sigma(\log|f|,\varphi_0)\ge\lim_{l\rightarrow+\infty}\nu_{j_l}(f)$.
	
	On the other hand, if $\sigma(\log|f|,\varphi_{0})>a$ for some $a>0$, then $\log|f|\le a\varphi_{0}+O(1)$ near $o$, which shows that $|f|^{2s}e^{-2as\varphi_{0}}$ is integrable near $o$ for any integer $s>0$. By the strong openness property of multiplier ideal sheaves, there exists $k>0$ such that $|f|^{2s}e^{-\frac{2as}{k}\log|I_k|}$ is integrable near $o$. There exists $l_1\ge1$ such that $I_k=\cap_{l\ge l_1} I_{k}^{j_{l_1}}\subset  I_{k}^{j_{l}}$ for any $l\ge l_1$. Note that 
	\[\frac{1}{k}\log|I_{k}|\le\frac{1}{k}\log|I_{k}^{j_{l}}|+O(1)\le \Phi_{j_l}+O(1)\]
	near $o$. Then we obtain that for any $s$, there exists $l_s$ such that for any $l\ge l_s$,
	$|f|^{2s}e^{-2as\Phi_{j_l}}$
	is integrable near $o$. Set 
	\[M\coloneqq \sup_{j}\sup_{(h,o)\in\mathcal{O}_o\setminus\{0\}}\big(c_o^h(\Phi_{j})-\nu_j(h)\big)<+\infty.\]
	Now, we have
	$$\nu_{j_l}(f)=\frac{\nu_{j_l}(f^s)}{s}\ge\frac{c_o^{f^s}(\Phi_{j_l})-M}{s}\ge\frac{as-M}{s}$$
	for any $s$ and $l\ge l_s$. Letting $l\rightarrow+\infty$ and $s\rightarrow+\infty$, we get $ \lim_{l\rightarrow+\infty}\nu_{j_l}(f)\ge a$, which induces $ \lim_{l\rightarrow+\infty}\nu_{j_l}(f)\ge \sigma(\log|f|,\varphi_{0})$.
	
	Thus, we have $\lim_{l\rightarrow+\infty}\nu_{j_l}(f)=\nu_0(f)=\sigma(\log|f|,\varphi_{0})$ for any $(f,o)\in\mathcal{O}_o$.
	
	\
	
	\emph{Step 3.} Find a tame maximal weight near $o$.

	If $|f|^2e^{-2a\varphi_{0}}$ is integrable near $o$, by the strong openness property of multiplier ideal sheaves, there exists $k>0$ such that $|f|^{2}e^{-\frac{2a}{k}\log|I_k|}$ is integrable near $o$. There exists $l_1\ge1$ such that $I_k=\cap_{l\ge l_1} I_{k}^{j_{l_1}}\subset  I_{k}^{j_{l}}$ for any $l\ge l_1$. Note that 
	\[\frac{1}{k}\log|I_{k}|\le\frac{1}{k}\log|I_{k}^{j_{l}}|+O(1)\le \Phi_{j_l}+O(1)\]
	near $o$. Then for any $l\ge l_1$, $|f|^{2}e^{-2a\Phi_{j_l}}$ is integrable near $o$.
	By the definition of $M$, we have
	$$\nu_{j_l}(f)\ge c_o^{f}(\Phi_{j_l})-M\ge a-M$$
	for any $l\ge l_1$. Letting $l\rightarrow+\infty$, we have 
	\[\sigma(\log|f|,\varphi_{0})= \lim_{l\rightarrow+\infty}\nu_{j_l}(f)\ge a-M.\]
	Thus, $c_o^{f}(\varphi_0)-\sigma(\log|f|,\varphi_0)\le M$ when $\sigma(\log|f|,\varphi_0)<+\infty$.
	As there exists $(f_0,o)$ such that $\sup_j\nu_j(f_0)<+\infty$ and $\lim_{l\rightarrow+\infty}\nu_{j_l}(f_0)=\nu_0(f)=\sigma(\log|f_0|,\varphi_{0})$, we have $c_o^{f_0}(\varphi_0)<+\infty$, which implies that there exists $C_2>0$ such that $\varphi_{0}\le C_2\log|z|+O(1)$ near $o$. Then we have
	$$\sup_{(f,o)\in\mathcal{O}_o\setminus\{0\}}\big(c_o^{f}(\varphi_0)-\sigma(\log|f|,\varphi_0)\big)\le M.$$
	It is clear that there exists $C_1>0$ such that $\varphi_{0}\ge C_1\log|z|+O(1)$ near $o$.
	Using Lemma \ref{l:relative type}, Lemma \ref{l:tame-maximal} and Lemma  \ref{l:tame-maximal2}, we can find a tame maximal weight $\tilde\varphi_{0}$ satisfying  $\nu_0(f)=\sigma(\log|f|,\tilde\varphi_{0})$ for any $(f,o)\in\mathcal{O}_o$.
\end{proof}

\begin{Remark}
\label{r:proof of 1.2}Following the above proof, we know that if $\{\nu_j\}$ does not satisfy conditions $(1)-(3)$, and only satisfies  $\sup_j\nu_j(g)<+\infty$ for any holomorphic function $g$ near $o$, we can also obtain that: {There exists a subsequence of $\{\nu_j\}$ denoted by $\{\nu_{j_l}\}$ such that $\{\nu_{j_l}\}$ converges to a valuation $\nu_0$ on $\mathcal{O}_o$.}
\end{Remark}

Similarly, for relative types, we can also obtain the following convergence result.

Let $\{\Phi_{j}\}_{j\in\mathbb{Z}_{>0}}$ be a sequence of tame maximal weights near $o$ satisfying the following statements:
\begin{enumerate}
    \item There exists $A_1>0$ such that $A_1\log|z|\le\Phi_{j}+O(1)$ near $o$ holds for all $j$;
    \item $M\coloneqq\sup_{j}\sup_{(f,o)\in\mathcal{O}_o\setminus\{0\}}\big(c_o^f(\Phi_{j})-\sigma(\log|f|,\Phi_{j})\big)<+\infty$.
\end{enumerate}

\begin{Proposition}
\label{thm:converge2}
There exists a subsequence of $\{\Phi_{j}\}$ denoted by $\{\Phi_{j_l}\}$ and a tame maximal weight $\varphi_{0}$ near $o$ such that 
\[\lim_{l\rightarrow+\infty}\sigma(\psi,\Phi_{j_l})=\sigma(\psi,\varphi_{0})\]
for any plurisubharmonic function $\psi$ near $o$, and
\[\sup_{(f,o)\in\mathcal{O}_o\setminus\{0\}}\big(c_o^f(\varphi_{0})-\sigma(\log|f|,\varphi_{0})\big)\le M.\]
\end{Proposition}
\begin{proof}
	Denote $I_k^j\coloneqq \{(f,o)\in\mathcal{O}_o\colon\sigma(\log|f|,\Phi_{j})\ge k\}$ for any $j,k\in\mathbb{Z}_{>0}$. Note that
	\begin{equation*}
	    \begin{split}
	        \sigma(\log|f_1f_2|,\Phi_{j})&\ge\sigma(\log|f_1|,\Phi_{j})+\sigma(\log|f_2|,\Phi_{j})\\ &\ge\max\{\sigma(\log|f_1|,\Phi_{j}),\sigma(\log|f_2|,\Phi_{j})\}
	    \end{split}
	\end{equation*}
	and 
	\[\sigma(\log|f_1+f_2|,\Phi_{j})\ge\sigma(\log(|f_1|+|f_2|),\Phi_{j})=\min\{\sigma(\log|f_1|,\Phi_{j}),\sigma(\log|f_2|,\Phi_{j})\}.\]
	Thus $I_k^j$ is an ideal of $\mathcal{O}_o$.
	
	Following the proof of Proposition \ref{thm:converge}, we obtain that:
	\begin{enumerate}
	    \item There exists a subsequence $\{I_k^{j_l}\}$ of $I_k^j$ such that $\cup_{l_1\ge1}\big(\cap_{l\ge l_1}I_k^{j_l}\big)$ is maximal for any $k$. It means that for any subsequence $\{I_k^{j'_l}\}$ of $\{I_k^{j_l}\}$, $\cup_{l_1\ge1}\big(\cap_{l\ge l_1}I_k^{j'_l}\big)=\cup_{l_1\ge1}\big(\cap_{l\ge l_1}I_k^{j_l}\big)$.
	    \item Denote $I_k\coloneqq \cup_{l_1\ge1}\big(\cap_{l\ge l_1}I_k^{j_l}\big)$. It is clear that $I_{k_1}\cdot I_{k_2}\subset I_{k_1k_2}$ for any $k_1,k_2\in\mathbb{Z}_{>0}$. As $A_1\log|z|\le \Phi_j$ for any $j$, there exist finitely many polynomials $f_{k,1},\ldots,f_{k,r_k}$ which generate $I_k$. Denote $|I_k|^2\coloneqq \sum_{1\le r\le r_k}|f_{k,r}|^2$ (without loss of generality, we assume $|I_k|<1$ on the polydisc $\Delta^n$).
	    Let $\varphi_0\coloneqq (\max_{k}\frac{1}{k}\log|I_k|)^*$ on $\Delta^n$, which is a plurisubharmonic function on $\Delta^n$ and equal almost everywhere to $\max_{k}\frac{1}{k}\log|I_k|$.
	    \item For any $(f,o)$, 
	    \[\lim_{l\rightarrow+\infty}\sigma(\log|f|,\Phi_{j_l})=\sigma(\log|f|,\varphi_0)\]
	    and 
	    \[\sup_{(f,o)\in\mathcal{O}_o\setminus\{0\}}\big(c_o^f(\varphi_{0})-\sigma(\log|f|,\varphi_{0})\big)\le M.\]
	    \item Using Lemma \ref{l:relative type}, Lemma \ref{l:tame-maximal} and Lemma  \ref{l:tame-maximal2}, we can find a tame maximal weight $\tilde\varphi_{0}$ satisfying  $\sigma(\log|f|,\tilde\varphi_{0})=\sigma(\log|f|,\varphi_{0})$  and $c_o^{f}(\tilde\varphi_0)=c_0^{f}(\varphi_{0})$ for any $(f,o)\in\mathcal{O}_o$.
	\end{enumerate}
	
	In the following, we prove $\lim_{l\rightarrow+\infty}\sigma(\psi,\Phi_{j_l})=\sigma(\psi,\tilde\varphi_{0})$ for any plurisubharmonic function $\psi$ near $o$.
	
	Fixed a plurisubharmonic function $\psi$ near $o$, for any $m\ge0$, denote $\psi_m\coloneqq \frac{1}{m}\log|\mathcal{I}(m\psi)_o|$. As $\psi\le\psi_m+O(1)$ near $o$, we have $\sigma(\psi_m,\Phi_{j_l})\le\sigma(\psi,\Phi_{j_l})$	and  $\sigma(\psi_m,\tilde\varphi_{0})\le\sigma(\psi,\tilde\varphi_{0})$.
	
	If $\sigma(\psi,\tilde\varphi_{0})>a>0$, then $|\mathcal{I}(m\psi)_o|^2e^{-2ma\tilde\varphi_{0}}$ is integrable near $o$ for any $m$. By $\sup_{(f,o)\in\mathcal{O}_o\setminus\{0\}}\big(c_o^f(\tilde\varphi_{0})-\sigma(\log|f|,\tilde\varphi_{0})\big)\le M$, we have
	$$\sigma(\psi_m,\tilde\varphi_{0})=\frac{1}{m}\sigma(\log|\mathcal{I}(m\psi)_o|,\tilde\varphi_{0})\ge\frac{ma-M}{m},$$
	which shows
	\begin{equation}
		\label{eq:0722a}\sigma(\psi_m,\tilde\varphi_{0})\le\sigma(\psi,\tilde\varphi_{0})\le\sigma(\psi_m,\tilde\varphi_{0})+\frac{M}{m}
	\end{equation}
	for any $m$.
	Similarly, we have
	\begin{equation}
		\label{eq:0722b}\sigma(\psi_m,\Phi_{j_l})\le\sigma(\psi,\Phi_{j_l})\le\sigma(\psi_m,\Phi_{j_l})+\frac{M}{m}
	\end{equation}
	for any $m$ and $l$.
	
	Since  $\lim_{l\rightarrow+\infty}\sigma(\log|f|,\Phi_{j_l})=\sigma(\log|f|,\tilde\varphi_0)$ holds for any $(f,o)$, we have
	\begin{equation}
		\nonumber\begin{split}
			\lim_{l\rightarrow+\infty}\sigma(\psi_m,\Phi_{j_l})&=\frac{1}{m}\lim_{l\rightarrow+\infty}\sigma(\log|\mathcal{I}(m\psi)_o|,\Phi_{j_l})\\
			&=\frac{1}{m}\lim_{l\rightarrow+\infty}\min_{1\le k\le p_k}\{\sigma(\log|f_{m,k}|,\Phi_{j_l})\}\\
			&=\frac{1}{m}\min_{1\le k\le p_k}\{\sigma(\log|f_{m,k}|,\tilde\varphi_{0})\}\\
			&=\sigma(\psi_m,\tilde\varphi_{0}),	
		\end{split}
	\end{equation}
	where $\mathcal{I}(m\psi)_o$ is generated by $\{(f_{m,k},o)\}_{1\le k\le p_k}$. Combining with \eqref{eq:0722a} and \eqref{eq:0722b}, we obtain $\lim_{l\rightarrow+\infty}\sigma(\psi,\Phi_{j_l})=\sigma(\psi,\tilde\varphi_{0})$.
	
	Thus, Proposition \ref{thm:converge2} holds.
\end{proof}

 For an analogous compactness result in algebraic settings, see \cite[Proposition 5.9]{JON-Mus2012}.

\section{Proofs of Theorem \ref{thm:interpolation} and its corollaries}

In this section, we prove Theorem \ref{thm:interpolation} by using Theorem \ref{thm:converge} and some results on Zhou valuations and Tian functions, and then we prove Corollary \ref{c:interpolation-n functions},  \ref{c:valuation-tame valuation},  \ref{C:interpolation-comp-poly}, \ref{C:interpolation-real}, \ref{C:interpolation-real-poly}.

\begin{proof}[\textbf{Proof of Theorem \ref{thm:interpolation}}]
	
	If there exists a valuation $\nu$  such that $\nu(f_j)=a_j$ for any $j$, then by Lemma \ref{l:relativetype}, 
	\[\sigma(\log|F|,\varphi)\le \nu(F)=\sum_{0\le j\le m}\nu(f_j)=\sum_{1\le j\le m}a_j.\]
	On the other hand, by the definition of relative type, clearly we also have
	\[\sigma(\log|F|,\varphi)\ge\sum_{0\le j\le m}\sigma(\log|f_j|,\varphi)\ge \sum_{1\le j\le m}a_j.\]
	Thus, $\sigma(\log|F|,\varphi)=\sum_{1\le j\le m}a_j$ holds.
	
	We prove this theorem in two steps.
	
	\
	
	\emph{Step 1. $(2)\Rightarrow(1)$.}
	
	As $\sigma(\log|F|,\varphi)<+\infty,$ we have $f_j(o)=0$ for any $1\le j\le m.$
	Denote $\varphi_N\coloneqq \max\{\varphi,N\log|z|\}$ for any $N>0$. 
	Then we have $\sigma(\log|f_j|,\varphi_N)\ge\sigma(\log|f_j|,\varphi)\ge a_j$ for all $0\le j\le m$. As $e^{\varphi}$ is H\"older continuous, there exists $C>0$ such that $c_o^{f}(\varphi)\le \sigma(\log|f|,\varphi)+C$ for any $(f,o)\in\mathcal{O}_o$. Thus,  for any $k>0$, there exists $N_k$ such that
	\begin{equation}
		\nonumber
		\begin{split}
			\sigma(\log|F|,\varphi)&\ge\frac{c_o^{F^k}(\varphi)-C}{k}\\
			&\ge\frac{c_o^{F^k}(\varphi_{N_k})-1-C}{k}\\
			&\ge\frac{k\sigma(\log|F|,\varphi_{N_k})-1-C}{k}\\
			&\ge \lim_{N\rightarrow+\infty}\sigma(\log|F|,\varphi_{N})-\frac{1+C}{k},
		\end{split}
	\end{equation}
	where the second ``$\ge$" follows from Lemma \ref{l:lct-N}. Then we have 
	\[ \lim_{N\rightarrow+\infty}\sigma(\log|F|,\varphi_{N})=\sigma(\log|F|,\varphi).\]
	
	Consider the Tian function
	$$\Tn (t)\coloneqq \sup\big\{c\ge 0 \colon |F|^{2t}e^{-2c\varphi_N}\text{ is integrable near }o\big\}.$$
	As $\varphi_N$ is a tame weight, it follows from Remark \ref{r:tame} that 
	\begin{equation}
	    \label{eq:0918a}
	    \lim_{t\rightarrow+\infty}\Tn '_-(t)=\lim_{t\rightarrow+\infty}\Tn '_+(t)=\sigma(\log|F|,\varphi_N).
	\end{equation}
	Then for every positive integer $s$, $|F|^{2s}e^{-2\Tn (s)\varphi_N}$ is not integrable near $o$ (because of the strong openness property of multiplier ideal sheaves, see \cite{GZopen-c}), then by \Cref{rem:max_existence}, there exists a local Zhou weight $\Phi_{N,s}$ near $o$ related to $|F|^{2s}$ such that
	$$\Phi_{N,s}\ge \Tn (s)\varphi_N.$$
	Consider the Tian function
	$$\widetilde{\Tn }(t)\coloneqq \sup\Big\{c\ge0\colon |F|^{2t}e^{-2c\frac{\Phi_{N,s}}{\Tn (s)}}\text{ is integrable near }o\Big\}.$$
	Thus, we have $\widetilde{\Tn }(s)=\Tn (s)$ and $\widetilde{\Tn }(t)\ge \Tn (t)$ for any $t\ge0$. It follows from \Cref{prop:concave_B} that $\widetilde{\Tn }'(s)=\sigma\Big(\log|F|,\frac{\Phi_{N,s}}{\Tn (s)}\Big)$. Then we have
	\begin{equation}
	    \label{eq:0918b}
	    \Tn '_+(s)\le\sigma\Big(\log|F|,\frac{\Phi_{N,s}}{\Tn (s)}\Big)\le\Tn '_-(s).
	\end{equation}

	Denote the corresponding Zhou valuation of $\Phi_{N,s}$ by $\nu_{N,s}$. Then we have
	$$\Tn (s)\nu_{N,s}(f_j)=\Tn (s)\sigma(\log|f_j|,\Phi_{N,s})\ge\sigma(\log|f_j|,\varphi_N)\ge a_j$$
	for all $0\le j\le m$ and
	$$\Tn (s)\nu_{N,s}(F)=\Tn (s)\sigma(\log|F|,\Phi_{N,s})=\sigma\Big(\log|F|,\frac{\Phi_{N,s}}{\Tn (s)}\Big).$$
	Note that 
	\[\varphi_{N}=\log\Big(\sum_{1\le j\le m}|f_j|^{\frac{1}{a_j}}+\sum_{1\le k\le n}|z_k|^N\Big)+O(1)\]
	near $o$. By Lemma \ref{l:holder}, there exists a constant $C$ (independent of $N$) such that $\sup_{(f,o)\in\mathcal{O}_o^*}(c_o^{f}(\varphi_{N})-\sigma(\log|f|,\varphi_{N}))< C$. 
	As $\widetilde{\Tn }(t)$ is a concave function, we have
	\begin{equation}
		\label{eq:0712a}
		c_o\Big(\frac{\Phi_{N,s}}{\Tn (s)}\Big)=\widetilde{\Tn }(0)\le \Tn (s)-s\Tn '_+(s).
	\end{equation}
	It follows from  the concavity of $\Tn (t)$ and equality \eqref{eq:0918a} that $\Tn (s)-s\Tn '_+(s)\le c_o^{F^s}(\varphi_N)-\sigma(\log|F^s|,\varphi_N)< C$. Then the inequality \eqref{eq:0712a} implies $c_o\big(\frac{\Phi_{N,s}}{\Tn (s)}\big)< C$ for any $N$ and $s$, which shows that the Lelong number $\nu_o\big(\frac{\Phi_{N,s}}{\Tn (s)}\big)\ge\frac{1}{C}$ by Lemma \ref{l:Lelong}. Thus,
	\begin{equation}
		\label{eq:0712b}
		\Tn (s)\nu_{N,s}(f)\le C\mathrm{ord}_o(f)
	\end{equation}
	for any $N$, $s$ and $(f,o)$.
	
	By equality \eqref{eq:0918a}, inequalities \eqref{eq:0918b}, \eqref{eq:0712b}, and the proof of Theorem \ref{thm:converge} (see \cref{r:proof of 1.2}), there exists a subsequence of valuations of $\{\Tn (s)\nu_{N,s}\}_{s}$, which converges to a valuation $\nu_N$. Then we have
	$$\nu_N(f_j)\ge a_j, \quad j=0,1,\ldots,m,$$
	and
	$$\nu_N(F)=\sigma(\log|F|,\varphi_N).$$
	By the same reason, there exists a subsequence of valuations of $\nu_N$, which converges to a valuation $\nu$. Then we have
	$\nu(f_j)\ge a_j$ for any $0\le j\le m$ and
	\[\sum_{0\le j\le m}\nu(f_j)=\nu(F)=\lim_{N\rightarrow+\infty}\sigma(\log|F|,\varphi_N)=\sigma(\log|F|,\varphi)=\sum_{1\le j\le m}a_j,\]
	which implies $\nu(f_j)= a_j$ for $j=0,1,\ldots,m$ since $a_0=0$.

	\
	
	\emph{Step 2. If $\sigma(\log|F|,\varphi)=\sum_{1\le j\le m}a_j$ and $o$ is an isolated point in $\cap_{1\le j\le m}\{f_j=0\}$, there exists a valuation $\nu$ and a tame maximal weight $\varphi_{\nu}$ such that $\nu(f_j)=a_j$ for any $j$ and $\nu(f)=\sigma(\log|f|,\varphi_{\nu})$ for any $(f,o)\in\mathcal{O}_o$.}

	Consider the Tian function
	$$\Tn (t)\coloneqq \sup\big\{c\ge0:|F|^{2t}e^{-2c\varphi}\text{ is integrable near }o\big\}.$$
	As $\varphi$ is a tame weight, we have
	\begin{equation}
		\label{eq:0704b}\lim_{t\rightarrow+\infty}\Tn '_-(t)=\lim_{t\rightarrow+\infty}\Tn '_+(t)=\sigma(\log|F|,\varphi).
	\end{equation}
		As $o$ is an isolated point in $\cap_{1\le j\le m}\{f_j=0\}$, there exists $A_1>0$ such that $\varphi\ge A_1\log|z|+O(1)$ near $o$.
Then for any positive integer $s$, by Remark \ref{rem:max_existence}, there exists a local Zhou weight $\Phi_{s}$ near $o$ with respect to $|F|^{2s}$ such that
	$$\Phi_{s}\ge \Tn (s)\varphi$$
	and
	$$\Tn '_+(s)\le\sigma\Big(\log|F|,\frac{\Phi_{s}}{\Tn (s)}\Big)\le\Tn '_-(s).$$
	Denote the corresponding Zhou valuation of $\Phi_{s}$ by $\nu_{s}$. Then we have
	\begin{equation}
	    \label{eq:0918c}
	    \Tn (s)\nu_{s}(f_j)=\Tn (s)\sigma(\log|f_j|,\Phi_{s})\ge\sigma(\log|f_j|,\varphi)\ge a_j
	\end{equation}
	for any $0\le j\le m$ and
	\begin{equation}
		\label{eq:0704a}
		\begin{aligned}
		\Tn (s)\nu_{s}(F)&=\Tn (s)\sigma(\log|F|,\Phi_{s})\\
		&=\sigma\Big(\log|F|,\frac{\Phi_{s}}{\Tn (s)}\Big)\in[\Tn '_+(s),\Tn '_-(s)],   
		\end{aligned}
	\end{equation}
	which shows $\sup_s\Tn (s)\nu_{s}(F)<+\infty$.

	As $\Phi_s$ is a local Zhou weight with respect to $|F|^{2s}$, by Theorem \ref{thm:valu-jump}, we have
	$$c_o^{f}(\Phi_s)-\nu_s(f)\le 1-\nu_s(F)s$$
	for any $(f,o).$
	It follows from the concavity of $\Tn (s)$, \eqref{eq:0704b} and \eqref{eq:0704a} that
	$$\Tn (s)-s\Tn (s)\nu_s(F)\le \Tn (s)-s\Tn '_+(s)\le \Tn (s)-s\sigma(\log|F|,\varphi).$$
	Note that $\sup_{(f,o)\in\mathcal{O}_o\setminus\{0\}}(c_o^f(\varphi)-\sigma(\log|f|,\varphi))<\infty$. Then we have
	\begin{equation}\nonumber
		\begin{split}
			\sup_s\sup_{(f,o)\in\mathcal{O}_o\setminus\{0\}}\left(c_o^{f}(\frac{\Phi_s}{\Tn (s)})-\Tn (s)\nu_s(f)\right)&\le\sup_s (\Tn (s)-\Tn (s)s\nu_s(F))\\
			&\le\sup_s(\Tn (s)-s\sigma(\log|F|,\varphi))\\
			&<+\infty.
		\end{split}
	\end{equation}
	Note that $\frac{\Phi_s}{\Tn(s)}\ge\varphi\ge A_1\log|z|+O(1)$ near $o$ for any $s$, where $\frac{\Phi_s}{\Tn(s)}$ is the corresponding weight of valuation $\Tn (s)\nu_{s}$.
	By Theorem \ref{thm:converge}, there exist a subsequence of $\{\Phi_s\}$ denoted by $\{\Phi_{s_l}\}$ ($\lim_{l\rightarrow+\infty}s_l=+\infty$), a valuation $\nu$ and a tame maximal weight $\varphi_{\nu}$ such that
	$$\lim_{l\rightarrow+\infty}\Tn (s_l)\nu_{s_l}(f)=\nu(f)=\sigma(\log|f|,\varphi_{\nu})$$
	for any $(f,o)$. Combining equalities \eqref{eq:0704b}, \eqref{eq:0918c} and \eqref{eq:0704a}, 
	 we have 
	\[\nu(F)=\lim_{l\rightarrow+\infty}\Tn '_+(s_l)=\sigma(\log|F|,\varphi)=\sum_{1\le j\le m}a_j\]
	and $\nu(f_j)\ge a_j$ for all $j$, which implies $\nu(f_j)=a_j$ for all $j=0,\dots,m$.
	
	Thus, Theorem \ref{thm:interpolation} holds.
\end{proof}

\begin{proof}[\textbf{Proof of Corollary \ref{c:interpolation-n functions}}]
	By Theorem \ref{thm:interpolation}, it suffices to prove that $\sigma(\log|F|,\varphi)=\sum_{1\le j\le n}a_j$, where $F\coloneqq \prod_{1\le j\le n}f_j$ and $\varphi\coloneqq \log(\sum_{1\le j\le n}|f_j|^{\frac{1}{a_j}})$. By definition, we have $\sigma(\log|F|,\varphi)\ge\sum_{1\le j\le n}a_j$. In the following, we prove $\sigma(\log|F|,\varphi)\le\sum_{1\le j\le n}a_j$.
	
	Let $G\colon U\rightarrow\mathbb{C}^n$ be a holomorphic map defined by $G\coloneqq (f_1,\ldots,f_n)$, where $U$ is an open subset of $\mathbb{C}^n$.
	Let $\tilde F\coloneqq \prod_{1\le j\le n}w_j$ and $\tilde\varphi\coloneqq \log(\sum_{1\le j\le n}|w_j|^{\frac{1}{a_j}})$ on $\mathbb{C}^n$, where $(w_1,\ldots,w_n)$ is the coordinate system on $\mathbb{C}^n$. By a direct calculation, we have
	\begin{equation}
		\label{eq:0705a}
		\lim_{l\rightarrow+\infty}\frac{c_o^{\tilde{F}^l}(\tilde\varphi)}{l}=\sum_{1\le j\le n}a_j.
	\end{equation}
	As $\cap_{1\le j\le n}\{f_j=0\}=\{o\}$, we have 
	\[c_o^{\tilde{F}^l}(\tilde\varphi)=\sup\{c\ge 0\colon |F|^{2l}|J_{G}|^2e^{-2c\varphi} \ \text{is integrable near} \ o\} \eqqcolon c_l,\]
	where $J_G$ is the Jacobian determinant of the map $G$. By definition, we know
	$$\lim_{l\rightarrow+\infty}\frac{c_l}{l}\ge \sigma(\log|F|,\varphi).$$
	Thus, it follows from \eqref{eq:0705a} that $\sigma(\log|F|,\varphi)\le\sum_{1\le j\le n}a_j$. Corollary \ref{c:interpolation-n functions} holds.
\end{proof}

\begin{proof}[\textbf{Proof of Corollary \ref{c:valuation-tame valuation}}]
	We can find holomorphic functions $\{g_k\}_{1\le k\le r}$ near $o$ such that 
	\[\big(\cap_{1\le j\le m}\{f_j=0\}\big)\cap \big(\cap_{1\le k\le r}\{g_k=0\}\big)=\{o\}.\]
	Denote 
	\[\varphi\coloneqq \log\Big(\sum_{1\le j\le m}|f_j|^{\frac{1}{a_j}}+\sum_{1\le k\le r}|g_k|^{\frac{1}{b_k}}\Big)\]
	and 
	\[F\coloneqq \prod_{1\le j\le m}f_j\times\prod_{1\le k\le r}g_k,\]
	where $b_k\coloneqq \nu(g_k)>0$.
	
	It follows from Lemma \ref{l:relativetype} and the definition of relative types that 
	\[\sigma(\log|F|,\varphi)=\sum_{1\le j\le m}a_j+\sum_{1\le k\le r}b_k.\]
	Using Theorem \ref{thm:interpolation}, we obtain that Corollary \ref{c:valuation-tame valuation} holds.
\end{proof}

\begin{proof}[\textbf{Proof of Corollary \ref{C:interpolation-comp-poly}}]
	If $\sigma(\log|F|,\varphi)=\sum_{1\le j\le m}a_j$, Theorem \ref{thm:interpolation} shows that  there exists a valuation $\nu$ on $\mathcal{O}_o$ such that $\nu(f_j)=a_j$ for any $0\le j\le m$. Note that $\mathbb{C}[z_1,\ldots,z_n]$ is a subring of $\mathcal{O}_o$. Thus $\nu$ is also a valuation on $\mathbb{C}[z_1,\ldots,z_n]$.
	
	Conversely, if there exists a valuation $\nu$ on $\mathbb{C}[z_1,\ldots,z_n]$ satisfying that $\nu(f_j)=a_j$ for any $0\le j\le m$ and $\nu(z_l)>0$ for any $1\le l\le n$, by Lemma \ref{l:relativetype2} and the definition of $\sigma(\log|F|,\varphi)$, we have  $\sigma(\log|F|,\varphi)=\sum_{1\le j\le m}a_j$.
\end{proof}

\begin{proof}[\textbf{Proof of Corollary \ref{C:interpolation-comp-poly2}}]
	The sufficient part follows from Corollary \ref{C:interpolation-comp-poly}, then it suffices to prove the necessary part. If there exists a valuation $\nu$ on $\mathbb{C}[z_1,\ldots,z_n]$ satisfying that $\nu(f_0)=0$ and $\nu(f_j)=a_j>0$ for any $1\le j\le m$. Let $I$ be the ideal in $\mathbb{C}[z_1,\ldots,z_n]$ generated by $\{f_j\}_{1\le j\le m}$. 
	As $\cap_{1\le j \le m}\{f_j=0\}=\{o\}$, it follows from Hilbert's Nullstellensatz (see \cite{atiyah}) that  $z_k^N\in I$ for large $N$ and any $1\le k\le n$, which implies $\nu(z_k)>0$ for any $k$. Thus, we have $\sigma(\log|F|,\varphi)=\sum_{1\le j\le m}a_j$ by Corollary \ref{C:interpolation-comp-poly}.
\end{proof}

\begin{proof}[\textbf{Proof of Corollary \ref{C:interpolation-real}}]
	If $\sigma(\log|F|,\varphi)=\sum_{1\le j\le m}a_j$, Theorem \ref{thm:interpolation} shows that  there exists a valuation $\tilde\nu$ on $\mathcal{O}_o$ such that $\tilde\nu(P(f_j))=a_j$ for any $0\le j\le m$. Note that $P:C^{\omega}_{o'}\rightarrow\mathcal{O}_o$ is an injective ring homomorphism. Then there exists a valuation $\nu$ on $C^{\omega}_{o'}$ such that $\nu(f_j)=a_j$ for any $j$.
	
	Conversely, if there exists a valuation $\nu$ on $C^{\omega}_{o'}$ such that $\nu(f_j)=a_j$ for any $j$, by Lemma \ref{c:relativetype-real} and the definition of $\sigma(\log|F|,\varphi)$, we have  $\sigma(\log|F|,\varphi)=\sum_{1\le j\le m}a_j$.
\end{proof}

\begin{proof}[\textbf{Proof of Corollary \ref{C:interpolation-real-poly}}]
	If $\sigma(\log|F|,\varphi)=\sum_{1\le j\le m}a_j$, Corollary \ref{C:interpolation-real} shows that  there exists a valuation $\nu$ on $C^{\omega}_{o'}$ such that $\nu(f_j)=a_j$ for any $0\le j\le m$. Note that $\mathbb{R}[x_1,\ldots,x_n]$ is a subring of $C^{\omega}_{o'}$. Thus $\nu$ is also a valuation on  $\mathbb{R}[x_1,\ldots,x_n]$.
	
	Conversely, if there exists a valuation $\nu$ on $\mathbb{R}[x_1,\ldots,x_n]$ satisfying that $\nu(f_j)=a_j$ for any $1\le j\le m$ and $\nu(x_l)>0$ for any $1\le l\le n$,, by Lemma \ref{c:relativetype2-real} and the definition of $\sigma(\log|F|,\varphi)$, we have  $\sigma(\log|F|,\varphi)=\sum_{1\le j\le m}a_j$.
\end{proof}

\section{Proof of Theorem \ref{thm:valuationexistsweight}}

In this section, we prove Theorem \ref{thm:valuationexistsweight}.

	For any valuation $\nu$ on $\mathcal{O}_o$ centered at $o$ and any $m\in\mathbb{Z}_{>1}$, there exist finitely many polynomials $f^{\nu}_{m,1},\ldots,f^{\nu}_{m,l_m}$ (where we may assume $|f^{\nu}_{m,l}|<1$ on $\Delta^n$) which generate the ideal $\mathfrak{a}^{\nu}_{m}\coloneqq \{(f,o)\colon\nu(f)\ge m\}$. Denote by 
	\[\varphi_0\coloneqq \left(\sup_m\frac{1}{m}\log(\max_{l}|f^{\nu}_{m,l}|)\right)^*\]
	on $\Delta^n$.
	
	\
	
	\emph{Step 1. We prove $\lct^{(f)}(\mathfrak{a}^{\nu}_{\bullet})=c_o^f(\varphi_0)$.}
	
	Since 
	\[\lct^{(f)}(\mathfrak{a}^{\nu}_{\bullet})=\lim_{m\rightarrow+\infty}c_o^f\left(\frac{1}{m}\log|\mathfrak{a}^{\nu}_{m}|\right)=\sup_m c_o^f\left(\frac{1}{m}\log|\mathfrak{a}^{\nu}_{m}|\right)\]
	and $\frac{1}{m}\log|\mathfrak{a}^{\nu}_{m}|\le\varphi_0+O(1)$ near $o$, we have $\lct^{(f)}(\mathfrak{a}^{\nu}_{\bullet})\le c_o^f(\varphi_0)$. For any $(f,o)$ and constant $c>\lct^{(f)}(\mathfrak{a}^{\nu}_{\bullet})$, we have that $|f|^2e^{-2c\frac{1}{m}\log|\mathfrak{a}^{\nu}_{m}|}$ is not integrable near $o$ for any $m$. Note that for any $k>1$,
	$$\max_{1\le m\le k}\frac{1}{m}\log\big(\max_{l}|f^{\nu}_{m,l}|\big)\le \frac{1}{k!}\log|\mathfrak{a}^{\nu}_{k!}|+O(1)$$
	near $o$. Then we have that 
	\[|f|^2\exp\left(-2c\max_{1\le m\le k}\frac{1}{m}\log(\max_{l}|f^{\nu}_{m,l}|)\right)\]
	is not integrable near $o$. By Proposition \ref{p:effect_GZ}, Remark \ref{rem:effect_GZ} and Proposition \ref{pro:Demailly}, $|f|^2e^{-2c\varphi_0}$ is not integrable near $o$, which shows $c_o^{f}(\varphi_0)\le c$. Then we get $\lct^{(f)}(\mathfrak{a}^{\nu}_{\bullet})=c_o^f(\varphi_0)$ for any $(f,o)\in\mathcal{O}_o\setminus\{0\}$.
	
	\
	
	\emph{Step 2. $\sup_{(f,o)}\big(\lct^{(f)}(\mathfrak{a}^{\nu}_{\bullet})-\nu(f)\big)<+\infty$ $\Rightarrow$ there exists a tame maximal weight $\varphi_{\nu}$ s.t.  $\nu(f)=\sigma(\log|f|,\varphi_{\nu})$ and $c_o^{f}(\varphi_{\nu})=\lct^{(f)}(\mathfrak{a}^{\nu}_{\bullet})$.}
	
	For any positive rational number $\frac{q}{p}\le\nu(f)$, it is clear that $\nu(f^{p})\ge q$. Hence, 
	\[\frac{1}{q}\log|f^p|\le\frac{1}{q}\log|\mathfrak{a}_q^{\nu}|+O(1)\le\varphi_0+O(1)\]
	near $o$. Then we have
	$\sigma(\log|f|,\varphi_0)\ge\frac{q}{p}$, which implies $\sigma(\log|f|,\varphi_0)\ge \nu(f)$.
	
	Denote by 
	$$M\coloneqq \sup_{(f,o)}(\lct^{(f)}\big(\mathfrak{a}^{\nu}_{\bullet})-\nu(f)\big)<+\infty.$$
	Fixed $c>\nu(f)$,
	as $\lct^{(g)}(\mathfrak{a}^{\nu}_{\bullet})=c_o^g(\varphi_0)$ for any $(g,o)$, we have
	$$c_o^{f^m}(\varphi_0)=\lct^{(f^m)}(\mathfrak{a}^{\nu}_{\bullet})\le M+\nu(f^m)<M+cm$$
	for every positive integer $m$. Letting $m\rightarrow+\infty$, we get
	\[c=\lim_{m\rightarrow+\infty}\frac{M+cm}{m}\ge\lim_{m\rightarrow+\infty}\frac{c_o^{f^m}(\varphi_0)}{m}.\]
	Note that $c_o^{f^{m+1}}(\varphi_0)\ge c_o^{f^m}(\varphi_0)+\sigma(\log|f|,\varphi_0)$ for any $m$. Then we get
	$$c\ge\lim_{m\rightarrow+\infty}\frac{c_o^{f^m}(\varphi_0)}{m}\ge \sigma(\log|f|,\varphi_0).$$
	Since $c>\nu(f)$ is arbitrarily chosen, $\sigma(\log|f|,\varphi_0)\le \nu(f)$.
	
	Thus, we get $\sigma(\log|f|,\varphi_0) = \nu(f)$. As a consequence,  $$\sup_{(f,o)}\big(c_o^{f}(\varphi_0)-\sigma(\log|f|,\varphi_0)\big)=\sup_{(f,o)}\big(\lct^{(f)}(\mathfrak{a}^{\nu}_{\bullet})-\nu(f)\big)<+\infty.$$
	Since $\{(f,o)\colon \nu(f)>0\}$ is the maximal ideal of $\mathcal{O}_o$, we have $\varphi_0\ge N_0\log|z|+O(1)$ near $o$ for some $N_0>0$. By Lemma \ref{l:relative type}, Lemma \ref{l:tame-maximal} and Lemma  \ref{l:tame-maximal2}, there exists a tame maximal weight $\varphi_{\nu}$ such that $\nu(f)=\sigma(\log|f|,\varphi_{\nu})$  and $c_o^{f}(\varphi_{\nu})=\lct^{(f)}(\mathfrak{a}^{\nu}_{\bullet})$ for any holomorphic germ $(f,o).$
	
	\
	
	\emph{Step 3. There exists a tame maximal weight $\varphi_{\nu}$ s.t. $\nu(f)=\sigma(\log|f|,\varphi_{\nu})$ $\Rightarrow$ $\sup_{(f,o)}\big(\lct^{(f)}(\mathfrak{a}^{\nu}_{\bullet})-\nu(f)\big)<+\infty$.}
	
	For any $m>1$, 
	\[\frac{1}{m}\log\big(\max_{l}|f^{\nu}_{m,l}|\big)=\frac{1}{m}\log|\mathfrak{a}^{\nu}_{m}|+O(1)\le\varphi_{\nu}+O(1)\]
	near $o$. Then it follows from Proposition \ref{p:effect_GZ}, Remark \ref{rem:effect_GZ} and Proposition \ref{pro:Demailly} that
	$$\lct^{(f)}(\mathfrak{a}^{\nu}_{\bullet})=c_o^{f}(\varphi_{0})\le c_o^{f}(\varphi_{\nu})$$
	for any $(f,o)$. Hence, $$\sup_{(f,o)}\big(\lct^{(f)}(\mathfrak{a}^{\nu}_{\bullet})-\nu(f)\big)=\sup_{(f,o)}\big(c_o^{f}(\varphi_{0})-\nu(f)\big)\le\sup_{(f,o)}\big(c_o^{f}(\varphi_{\nu})-\nu(f)\big)<+\infty,$$
	where the last ``$<$" holds since $\varphi_{\nu}$ is a tame weight.
	
	Thus, Theorem \ref{thm:valuationexistsweight} holds.

\section{Proofs of Proposition \ref{p:quasi} and Theorem \ref{thm:multi-valua,quasi}}

In this section, we prove Proposition \ref{p:quasi} and Theorem \ref{thm:multi-valua,quasi}.

Firstly, we recall the following desingularization theorem due to Hironaka, which will be used in the proof of Proposition \ref{p:quasi}.
\begin{Theorem}[\cite{Hironaka}, see also \cite{BM-1991}]\label{thm:desing}
	Let $X$ be a complex manifold, and $M$ be an analytic sub-variety in $X$.  Then there is a local finite sequence of blow-ups $\mu_j\colon X_{j+1}\rightarrow X_j$ $(X_1\coloneqq X,j=1,2,\ldots)$ with smooth centers $S_j$ such that:
	\begin{enumerate}
	    \item Each component of $S_j$ lies either in $(M_j)_{\rm{sing}}$ or in $M_j\cap E_j$, where $M_1\coloneqq M$, $M_{j+1}$ denotes the strict transform of $M_j$ by $\mu_j$, $(M_j)_{\mathrm{sing}}$ denotes the singular set of $M_j$, and $E_{j+1}$ denotes the exceptional divisor $\mu^{-1}_j(S_j\cup E_j)$;
	    \item Let $M'$ and $E'$ denote the final strict transform of $M$ and the exceptional divisor respectively. Then:
	\begin{itemize}
	    \item[(a)] The underlying point-set $|M'|$ is smooth;
	    \item[(b)] $|M'|$and $E'$ simultaneously have only normal crossings.
	\end{itemize}
	\end{enumerate}
\end{Theorem}

The statement (b) in the above theorem means that, locally, there is a coordinate system in which $E'$ is a union of coordinate hyperplanes and $|M'|$ is a coordinate subspace.

\begin{proof}[\textbf{Proof of Proposition \ref{p:quasi}}]
As $|f_0|^2e^{-2\Phi_{o,\max}}$ is not integrable near $o$, there exists some $f_{0,s}$ ($1\le s\le m$) such that $|f_{0,s}|^2e^{-2\Phi_{o,\max}}$ is not integrable near $o$. As $c_o^{f_{0}}(\varphi)=1$ and $\varphi\le\Phi_{o,\max}+O(1)$ near $o$, we have $c_o^{f_{0,s}}(\varphi)=1$. By the definition of local Zhou weights, $\Phi_{o,\max}$ is also a local Zhou weight related to $|f_{0,s}|^2$ near $o$. Thus, to prove the proposition, it suffices to consider the case that $f_0$ is a holomorphic function (not a vector). 

Let $U$ be a small open neighborhood of $o$ in $\mathbb{C}^n$. By Theorem \ref{thm:desing}, let $\mu\colon X\to U$ be a proper modification which gives the resolution of the singularity set $\{\varphi=-\infty\}\subset U$ of $\varphi$, which is an analytic subset of $U$ since $\varphi$ has analytic singularities. As $\Phi_{o,\max}$ is a local Zhou weight related to $|f_0|^2$ near $o$, we have $c_o^{f_0}(\Phi_{o,\max})=1$, which implies that there exists some $z_0\in\mu^{-1}(o)$ such that $c_{z_0}^{w^{b}f_0\circ\mu}(\Phi_{o,\max}\circ\mu)=1$. Here $(W ; w_1, \ldots, w_n)$ is a coordinate ball centered at $z_0$ in $X$ such that $w^{b}\coloneqq \prod_{j=1}^n w_j^{b_j}$ is the zero divisor of the Jacobian $J_{\mu}$ (of $\mu$) and
	$$\varphi\circ\mu(w)=c\log|w^{a}|^2+u(w)$$
	on $W$, where $u\in C^{\infty}(\overline{W})$, $w^{a}=\prod_{j=1}^n w_j^{a_j}$.
	Due to $c_o^{f_0}(\varphi)=1$, $c_{z_0}^{w^{b}f_0\circ\mu}(\Phi_{o,\max}\circ\mu)=1$ and $\Phi_{o,\max}\ge\varphi+O(1)$ near $o$, we have $$c_{z_0}^{w^{b}f_0\circ\mu}(\varphi\circ\mu)=1.$$
	
	By Lemma \ref{l:5}, there exists a plurisubharmonic function 
	\[\psi_1=\max_{1\le j\le r}\big\{c_j\log|w_{k_j}|\big\} \quad (1\le k_1<\ldots<k_r\le n, \ c_j>0)\]
	near $z_0$ s.t. $\psi_1\ge\Phi_{o,\max}\circ\mu+O(1)$ and $c_{z_0}^{w^{b}f_0\circ\mu}(\psi_1)=1$. 
	Choosing arbitrary nonzero holomorphic germ $(g,o)\in \mathcal{O}_{U,o}$, consider the following three functions:
	$$\Tn_1(t)\coloneqq \sup \big\{c\in\mathbb{R}\colon |f_0|^2|g|^{2t}e^{-2c\Phi_{o,\max}}\text{ is integrable near }o\big\},$$
	$$\Tn_2(t)\coloneqq \sup \big\{c\in\mathbb{R}\colon |w^{b}f_0\circ\mu|^2|g\circ\mu|^{2t}e^{-2c\Phi_{o,\max}\circ\mu}\text{ is integrable near }z_0\big\},$$
	$$\Tn_3(t)\coloneqq \sup \big\{c\in\mathbb{R}\colon |w^{b}f_0\circ\mu|^2|g\circ\mu|^{2t}e^{-2c\psi_1}\text{ is integrable near }z_0\big\}.$$
	Then we have $\Tn_i(0)=1$ and $\Tn_i(t)$ is concave with respect to $t$ for $i=1,2,3$. Moreover,
	$$\Tn_3(t)\ge \Tn_2(t)\ge \Tn_1(t), \quad \forall\, t.$$
	\Cref{prop:concave_B} shows that $\Tn_1 (t)$ is differentiable at $t=0$, and $\Tn_1 (t)=\Tn_1 (0)+\sigma(\log|g|,\Phi_{o,\max})t$ holds for all $t\ge0$. Then by the basic properties of concave functions, we get
	$$\Tn_3(t)=\Tn_2(t)=1+\sigma(\log|g|,\Phi_{o,\max})t$$
	holds for all $t\ge0$. Combining with \Cref{example}, we have
	$$\sigma(\log|g\circ\mu|,\psi_1; z_0)=\sigma(\log|g|,\Phi_{o,\max})\eqqcolon\nu(g), \quad \forall\, (g,o)\in \mathcal{O}_{U,o}.$$
	
	Now, according to \cite[Lemma 3.1]{JON-Mus2014} (and using the notations there), we can extract a Kiselman number $\tau_{Z,D,\alpha}$ on $X$ with $Z=\bigcap_{1\le j\le r}(w_{k_j}=0)$, $D=\sum_{1\le j\le r}(w_{k_j}=0)$ and $\alpha=(1/c_1,\ldots,1/c_r)$ such that $\tau_{Z,D,\alpha}(\phi)=\sigma(\phi,\psi_1; z_0)$ for every plurisubharmonic function $\phi$ defined in an open neighborhood of $z_0$. Moreover, by \cite[Lemma 3.8]{JON-Mus2014}, there exists a quasimonomial valuation $\widetilde{\nu}$ on $\mathcal{O}_{U,o}$ which is adapted to $(X,D)$ such that 
	\[\widetilde{\nu}(g)=\tau_{Z,D,\alpha}\big((\log |g|)\circ \mu\big)=\sigma(\log|g\circ\mu|,\psi_1; z_0)=v(g), \quad \forall (g,o)\in\mathcal{O}_{U,o}.\]
	We conclude that the Zhou valuation $v$ which is associated with the Zhou weight $\Phi_{o,\max}$ is quasimonomial.
\end{proof}

\begin{proof}[\textbf{Proof of Theorem \ref{thm:multi-valua,quasi}}]
	By Theorem \ref{thm:multi-valua}, it is enough to prove $(3)\Rightarrow(2)$. In the following, we
	assume that the statement $(3)$ holds.
	
	Firstly, we prove that  
	\[\mathcal{I}\big(\max\{v,N\log|z|\}\big)_o\subset\mathcal{I}\big(\max\{u,N\log|z|\}\big)_o\]
	holds for any fixed $N>0$.
	We prove this by contradiction: if not, there is a holomorphic function $f$ near $o$ such that $(f,o)\in\mathcal{I}(\max\{v,N\log|z|\})_o$ and $(f,o)\not\in\mathcal{I}(\max\{u,N\log|z|\})_o$. By the strong openness property (see \cite{GZopen-c}), there exists $\epsilon_0>0$ such that  
	\[(f,o)\in\mathcal{I}\big((1+\epsilon_0)\max\{v,N\log|z|\}\big)_o.\]
	By Theorem \ref{thm:analytic weight}, there exists $N_1>0$ such that $\frac{1}{N_1}<\epsilon_0$ and 
	$$\mathcal{I}(\max\{u,N\log|z|\})_o=\mathcal{I}\left(\frac{1}{N_1}\log|J_{N_1+1}|\right)_o,$$
	where 
	\[J_{N_1+1}\coloneqq \mathcal{I}\big((N_1+1)\max\{u,N\log|z|\}\big)_o.\]
	By Remark \ref{rem:max_existence}, there exists a local Zhou weight $\Phi_{o,\max}$ related to $|f|^2$ near $o$ satisfying that
	$$\Phi_{o,\max}\ge c_o^{f}(\log|J_{N_1+1}|)\log|J_{N_1+1}|+O(1)$$
	near $o$, and $|f|^2e^{-2\Phi_{o,\max}}$ is not integrable near $o$. Since $|f|^2e^{-2\frac{1}{N_1}\log|J_{N_1+1}|}$ is not integrable near $o$, we have $c_o^{f}(\log|J_{N_1+1}|)\le \frac{1}{N_1}$. Note that $\max\{u,N\log|z|\}\le \frac{1}{N_1+1}\log|J_{N_1+1}|+O(1)$ near $o$ (by Lemma \ref{l:appro-Berg}). Then
	$$\sigma\big(\max\{u,N\log|z|\},\Phi_{o,\max}\big)\ge \sigma\left(\frac{1}{N_1+1}\log|J_{N_1+1}|,\Phi_{o,\max}\right)\ge\frac{N_1}{N_1+1}.$$
	By Proposition \ref{p:quasi}, the Zhou valuation associated with $\Phi_{o,\max}$ is quasimonomial. Hence, it follows from the statement $(3)$ and the above inequality that 
	\[\sigma(\max\{v,N\log|z|\},\Phi_{o,\max})\ge \frac{N_1}{N_1+1}>\frac{1}{1+\epsilon_0},\]
	i.e.,
	$$(1+\epsilon_0)\max\{v,N\log|z|\}\le \Phi_{o,\max}+O(1)$$
	near $o$, which contradicts to $(f,o)\in\mathcal{I}((1+\epsilon_0)\max\{v,N\log|z|\})_o$ and $(f,o)\not\in\mathcal{I}(\Phi_{o,\max})_o$. Then we have
	\begin{equation}
		\nonumber
		\mathcal{I}\big(\max\{v,N\log|z|\}\big)_o\subset\mathcal{I}\big(\max\{u,N\log|z|\}\big)_o
	\end{equation}
	holds for any fixed $N>0$. Similarly, we have
	\begin{equation}
		\label{eq:250413a}\mathcal{I}\big(t\max\{v,N\log|z|\}\big)_o\subset\mathcal{I}\big(t\max\{u,N\log|z|\}\big)_o
	\end{equation}
	holds for any fixed $N>0$ and any $t>0$.
	
	Now, we prove the statement $(2)$. Since  $\sigma(t\varphi_1+\varphi_2,\Phi_{o,\max})=t\sigma(\varphi_1,\Phi_{o,\max})+\sigma(\varphi_2,\Phi_{o,\max})$ for any plurisubharmonic functions $\varphi_1$, $\varphi_2$, and Zhou weight $\Phi_{o,\max}$ near $o$, it is enough to prove that $\mathcal{I}(v)_o\subset\mathcal{I}(u)_o.$ The inequality \eqref{eq:250413a} implies that
	$$c_o^h\big(\max\{v,N\log|z|\}\big)\le c_o^h\big(\max\{u,N\log|z|\}\big)$$
	for any $(h,o)\in\mathcal{O}_o$ and any $N>0$. By Lemma \ref{l:lct-N}, we have
	$$c_o^h(v)\le c_o^h(u), \quad \forall\,(h,o)\in\mathcal{O}_o,$$
	which implies $\mathcal{I}(v)_o\subset\mathcal{I}(u)_o$.
	
	Thus, Theorem \ref{thm:multi-valua,quasi} holds.
\end{proof}

	\section{Proofs of  Theorem \ref{thm:derivative} and Proposition \ref{p:differentiable-zhou valuation}}
	
	In this section, we prove Theorem \ref{thm:derivative} and Proposition \ref{p:differentiable-zhou valuation}.

\begin{proof}[\textbf{Proof of Theorem \ref{thm:derivative}}]
Let $\psi_{1}$ and $\psi_{2}$ be plurisubharmonic functions near $o$.	

Let $a>0$, by definition, we know that 
\begin{equation}
    \label{eq:multiplicative}
    \Tn '(0;a\psi_1)=a\Tn '(0;\psi_1).
\end{equation}

Denote 
$$G(t,s)\coloneqq \sup\big\{c\colon |f_{0}|^{2}e^{-2\varphi_{0}}e^{2t\psi_1}e^{2s\psi_2}e^{-2c\varphi} \ \text{is integrable near} \ o\big\}.$$
By H\"{o}lder's inequality, $G(t,s)$ is a concave function on $(-\epsilon,+\infty)\times (-\epsilon,+\infty)$ for some $\epsilon>0$.
 Note that $G(t,0)=\Tn (t;\psi_1)$, $G(0,s)=\Tn (s;\psi_2)$ and 
$$G(t,t)=\Tn (t;\psi_1+\psi_2).$$
It follows from $\Tn (t)$ is differentiable at  $t= 0$ for any plurisubharmonic function that $G_t'(0,0)$, $G_s'(0,0)$  and 
$G_{1,1}'(0,0)$
exist, where $G_{1,1}'(0,0)$ is the derivative of $G(t,t)$ at $t=0$. 
As $G(t,s)$ is concave, we have that
\begin{equation}
\label{eq:concave}
    \begin{split}
      G\left(\frac{1}{2}t,\frac{1}{2}t\right) \ge \frac{1}{2}G(t,0)+\frac{1}{2}G(0,t).
    \end{split}
\end{equation}

When $t>0$, the inequality \eqref{eq:concave} tells
\begin{equation*}
    \begin{split}
     \frac{ G(\frac{1}{2}t,\frac{1}{2}t) -G(0,0)}{\frac{1}{2}t}\ge \frac{G(t,0)-G(0,0)}{t}+\frac{G(0,t)-G(0,0)}{t}.
    \end{split}
\end{equation*}
Since $G_t'(0,0)$, $G_s'(0,0)$  and 
$G_{1,1}'(0,0)$
exist, letting $t\to 0$, we have
\begin{equation*}
    G_{1,1}'(0,0)\ge G_t'(0,0)+G_s'(0,0).
\end{equation*}
When  $t<0$, the inequality \eqref{eq:concave} implies
\begin{equation*}
    \begin{split}
     \frac{ G(\frac{1}{2}t,\frac{1}{2}t) -G(0,0)}{\frac{1}{2}t}\le \frac{G(t,0)-G(0,0)}{t}+\frac{G(0,t)-G(0,0)}{t}.
    \end{split}
\end{equation*}
Letting $t\to0$, as  $G_t'(0,0)$, $G_s'(0,0)$  and 
$G_{1,1}'(0,0)$
exist, we have
\begin{equation*}
    G_{1,1}'(0,0)\le G_t'(0,0)+G_s'(0,0).
\end{equation*}
 Thus, we get 
 \begin{equation}
 \label{eq:additive}
     G_{1,1}'(0,0)= G_t'(0,0)+G_s'(0,0).
 \end{equation}

Combining \eqref{eq:multiplicative} with \eqref{eq:additive},
the statement (1) of Theorem \ref{thm:derivative} holds.

Now we prove the statement (2) of Theorem \ref{thm:derivative}. For $i=1,2$, denote
$$\Tn _i(t)\coloneqq \sup\big\{c\colon |f_{0}|^{2}e^{-2\varphi_{0}}e^{2t\psi_i}e^{-2c\varphi} \ \text{is integrable near} \ o\big\}$$
and
$$\Tn _{\max}(t)\coloneqq \sup\big\{c\colon |f_{0}|^{2}e^{-2\varphi_{0}}e^{2t\max\{\psi_1,\psi_2\}}e^{-2c\varphi} \ \text{is integrable near} \ o\big\}.$$

Let $t>0$. Note that, 
$$e^{2t\max\{\psi_1,\psi_2\}}\le e^{2t\psi_1}+e^{2t\psi_2} \le 2e^{2t\max\{\psi_1,\psi_2\}}.$$
Then, for any $c<  \min\{\Tn _1(t),\Tn _2(t)\}$, one has
\begin{equation*}
    \begin{aligned}
        \int_U |f_{0}|^{2}e^{-2\varphi_{0}}e^{2t\max\{\psi_1,\psi_2\}}e^{-2c\varphi}\le &\ 
\int_U |f_{0}|^{2}e^{-2\varphi_{0}}e^{2t\psi_1}e^{-2c\varphi}\\
&+\int_U |f_{0}|^{2}e^{-2\varphi_{0}}e^{2t\psi_2}e^{-2c\varphi},
    \end{aligned}
\end{equation*}
where $U$ is a small open neighborhood of $o$, and hence
\[\Tn _{\max}(t)\ge \min\{\Tn _{1}(t),\Tn _{2}(t)\}.\] 
 
 For any $c\ge \min\{\Tn _1(t),\Tn _2(t)\}$, by the strong openness property, there exists at least one of $(\psi_i)_{i=1,2}$ (without loss of generality, may assume it is $\psi_1$) such that 
$ \int_V|f_{0}|^{2}e^{-2\varphi_{0}}e^{2t\psi_1}e^{-2c\varphi}=+\infty$, where $V$ is any open neighborhood of $o$.
Then we have
\begin{equation*}
\begin{split}
&2\int_V |f_{0}|^{2}e^{-2\varphi_{0}}e^{2t\max\{\psi_1,\psi_2\}}e^{-2c\varphi}\\
\ge& \ 
\int_V|f_{0}|^{2}e^{-2\varphi_{0}}e^{2t\psi_1}e^{-2c\varphi}+\int_V|f_{0}|^{2}e^{-2\varphi_{0}}e^{2t\psi_2}e^{-2c\varphi}\\
=& \ +\infty.
\end{split}
\end{equation*}
Hence $\Tn _{\max}(t)\le \min\{\Tn _1(t),\Tn _2(t)\}$. Thus, when $t>0$, we have $\Tn _{\max}(t)= \min\{\Tn _1(t),\Tn _2(t)\}$.

Since $\Tn (t)$ is differentiable at  $t= 0$ for any plurisubharmonic function, $\Tn _{\max}(0)=\Tn _1(0)=\Tn _{2}(0)$ and $\Tn _{\max}(t)= \min\{\Tn _1(t),\Tn _2(t)\}$, we have
$\Tn '_{\max}(0)= \min\{\Tn '_1(0),\Tn '_2(0)\}$. This proves the statement (2) of Theorem \ref{thm:derivative}.

Now we prove the statement (3) of Theorem \ref{thm:derivative}.

Let $\psi_i=\log|f_i|$ for $i=1,2$.
Denote $$\Tn _0(t)\coloneqq \sup\big\{c\colon |f_{0}|^{2}e^{-2\varphi_{0}}|f_1+f_2|^{2t}e^{-2c\varphi} \ \text{is integrable near} \ o\big\}.$$
Let $t\ge 0$. Note that
\begin{equation}\nonumber
  |f_1+f_2|^{2t} \le  (|f_1|+|f_2|)^{2t}\le \big(2\max\{|f_1|,|f_2|\}\big)^{2t}.
\end{equation}
Hence $\Tn _0(t)\ge \Tn _{\max}(t)$ for any $t\ge 0$. Since $\Tn (t)$ is differentiable at  $t= 0$ for any plurisubharmonic function, $\Tn _{\max}(0)=\Tn _0(0)$ and $\Tn _0(t)\ge \Tn _{\max}(t)$, we have $\Tn '_0(0)\ge\Tn '_{\max}(0)$. By the statement (2) of Theorem \ref{thm:derivative}, one has $\Tn '_{\max}(0)= \min\{\Tn '_1(0),\Tn '_2(0)\}$. Thus, we get
$$\Tn '_0(0)\ge\min\{\Tn '_1(0),\Tn '_2(0)\}.$$
Statement (3) of Theorem \ref{thm:derivative} holds.
\end{proof}

\begin{proof}[\textbf{Proof of Proposition \ref{p:differentiable-zhou valuation}}]
	It follows from Lemma \ref{rem:max_existence} that there exists a local Zhou weight $\Phi_{o,\max}$ related to $|f_0|^2e^{-2\varphi_0}$ near $o$ such that $\Phi_{o,\max}\ge\varphi$ near $o$. Denote
	\[A(g,t)\coloneqq \sup\big\{c\ge 0\colon |f_0|^2|g|^{2t}e^{-2\varphi_0-2c\Phi_{o,\max}} \ \text{is integrable near} \ o\big\}\]
	for any holomorphic function $g$ near $o$.
	Then
	\begin{equation}
		\label{eq:0224a}\Tn (t;f_0,\varphi_0,\varphi,\log|g|)\le A(g,t)
	\end{equation}
	for any $g$ and any $t$.
	Proposition \ref{prop:concave_B} shows that $A(g,t)$ is differentiable at $t=0$ and $A'(g,0)=\sigma(\log|g|,\Phi_{o,\max})$ for any $g$. As $\Tn (t;f_0,\varphi_0,\varphi,\log|g|)$ is differentiable at $t=0$ for any $g$, it follows from \eqref{eq:0224a} that
	$$\nu(g)\coloneqq \Tn '(0;f_0,\varphi_0,\varphi,\log|g|)= A'(g,0)=\sigma(\log|g|,\Phi_{o,\max})$$
	for any $g$. Thus $\nu$ is a Zhou valuation on $\mathcal{O}_o$.
\end{proof}

\vspace{.1in} {\em Acknowledgements}.
The second author was supported by National Key R\&D Program of China 2021YFA1003100 and NSFC-12425101. The third author was supported by NSFC-12401099 and the Talent Fund of Beijing Jiaotong University 2024-004. The fourth author was supported by NSFC-12501106.


\appendix \section{Algebraic viewpoints of valuative interpolation}

In this appendix, we present some algebraic viewpoints of the valuative interpolation problem.

\subsection{Equivalence of interpolations in different valuation spaces}
Let $X$ be a separated, regular, connected, excellent scheme over $\mathbb{Q}$. Denote by $\Val_{X}$ the space of valuations on $X$ with a center in $X$, and $\QM_X$ the space of quasimonomial valuations on $X$ (cf. \cite{JON-Mus2012}). For every point $\xi\in X$, we also denote by $\Val_{X,\xi}$ the valuations on $X$ centered at $\xi$. Let $\mathcal{I}$ be the set of nonzero ideals on $X$ (``ideal on $X$" for ``coherent ideal sheaf on $X$"). We consider the valuative interpolation problem in $\mathcal{I}$.

\begin{Lemma}\label{lem-qm.enough}
    Let $m\in\mathbb{Z}_{>0}$. Given nonzero ideals $\mathfrak{a}_1,\ldots,\mathfrak{a}_m\in \mathcal{I}$ on $X$, if there exists $v\in \Val_X$ with $v(\mathfrak{a}_i)>0$ for each $i$, then there exists a quasimonomial valuation $\widetilde{v}\in\Val_X$ such that $\widetilde{v}\le v$, $c_X(\widetilde{v})=c_X(v)$ and $\widetilde{v}(\mathfrak{a}_i)=v(\mathfrak{a}_i)$ for $i=1,\ldots,m$.
\end{Lemma}

\begin{proof}
    Set $\mathfrak{a}=\mathfrak{a}_1\cdot\ldots\cdot \mathfrak{a}_m$, which is an ideal on $X$. Then we can choose a log smooth pair $(Y,D)$ over $X$ which gives a log resolution of $\mathfrak{a}$. Consider the retraction map
    \[r_{Y,D}\colon \Val_X \rightarrow \QM(Y,D),\]
    and let $\widetilde{v}\coloneqq r_{Y,D}(v)$. According to \cite[Corollary 4.8]{JON-Mus2012}, we have $\widetilde{v}\le v$, and since $(Y,D)$ gives a log resolution of $\mathfrak{a}$, we also have $\widetilde{v}(\mathfrak{a})=v(\mathfrak{a})$. It follows that $\widetilde{v}(\mathfrak{a}_i)=v(\mathfrak{a}_i)$ for every $i=1,\ldots,m$ since each $v(\mathfrak{a}_i)>0$. Thus, the quasimonomial valuation $\widetilde{v}\in\QM(Y,D)$ satisfies all we need except that $c_X(\widetilde{v})=c_X(v)$. However, note that $\widetilde{v}\le v$ gives $\overline{c_X(\widetilde{v})}\ni c_X(v)$, so we will have $c_X(\widetilde{v})=c_X(v)$ once $\widetilde{v}(\mathfrak{m}_{X,\xi})>0$, where $\mathfrak{m}_{X,\xi}$ is the ideal defining $\bar{\xi}$. Thus, to obtain a valuation satisfying the desired conditions, we only need to add $\mathfrak{a}_{m+1}\coloneqq \mathfrak{m}_{X,\xi}$ in the sequence, and repeat the previous process for the sequence of ideals $\mathfrak{a}_1,\ldots,\mathfrak{a}_m,\mathfrak{a}_{m+1}\in \mathcal{I}$.
\end{proof}

Notably, Corollary \ref{c:interpolation-n functions} can be derived by combining Lemma \ref{prop-equiv.qm.Zhou.center} with \cite[Theorem 5.13]{BFJ08}.

Given $(\mathfrak{a}_1,\ldots,\mathfrak{a}_m)\in\mathcal{I}^{\oplus m}$ and $(b_1,\ldots,b_m)\in\mathbb{R}_{>0}^m$, it also follows from the proof of Lemma \ref{lem-qm.enough} and Zorn's lemma that the set
\[\big\{v\in\Val_X\colon v(\mathfrak{a}_i)=b_i, \ \forall\,i=1,\ldots,m\big\}\]
always has minimal elements once this set is non-empty, and every minimal element must be quasimonomial.

Lemma \ref{lem-qm.enough} shows that the solvability of valuative interpolation problem in $\Val_X$ is equivalent to which in the space of quasimonomial valuations on $X$, and we can also show it is equivalent to be solvability in $\ZVal_X$, the space of Zhou valuations on $X$ (see \cite{BGZ24} for the definition of Zhou valuation on general schemes).

\begin{Proposition}\label{prop-equiv.qm.Zhou.center}
    Let $m\in\mathbb{Z}_{>0}$. Let $(\mathfrak{a}_1, \ldots, \mathfrak{a}_m)\in \mathcal{I}^{\oplus m}$, $(b_1,\ldots, b_m)\in\mathbb{R}_{>0}^{m}$ and $\xi\in X$ a point. The following statements are equivalent.
    \begin{enumerate}
        \item There exists a valuation $v\in\Val_X$ (resp. $v\in \Val_{X,\xi}$) such that $v(\mathfrak{a}_i)=b_i$ for $i=1,\ldots,m$.
        \item There exists a quasimonomial valuation $v\in\QM_X$ (resp. $v\in \QM_{X}\cap\Val_{X,\xi}$) such that $v(\mathfrak{a}_i)=b_i$ for $i=1,\ldots,m$.
        \item There exists a Zhou valuation $v\in\ZVal_X$ (resp. $v\in \ZVal_X\cap\Val_{X,\xi}$ if $\xi$ is a closed point) such that $v(\mathfrak{a}_i)=b_i$ for $i=1,\ldots,m$.
    \end{enumerate}
\end{Proposition}

\begin{proof}
    (3) $\Rightarrow$ (1) is clear, and Lemma \ref{lem-qm.enough} shows (1) $\Rightarrow$ (2). We now prove (2) $\Rightarrow$ (3).
    
    Let $v\in\QM_X$ be a quasimonomial valuation with $v(\mathfrak{a}_i)>0$ for each $i$ (in particular $A(v)<+\infty$). According to \cite[Proposition 2.5]{JON-Mus2014}, we can find a nonzero ideal $\mathfrak{q}$ on $X$ and a graded sequence $\mathfrak{c}_{\bullet}$ on $X$ such that $v$ computes $\lct^{\mathfrak{q}}(\mathfrak{c}_{\bullet})$. Then it follows from \cite[Theorem 7.8]{JON-Mus2012} that $v$ also computes $\lct^{\mathfrak{q}}(\mathfrak{c}_{\bullet}^v)$, where $\mathfrak{c}_{\bullet}^v$ is the graded sequence associated to the valuation $v$ (cf. \cite[Example 2.2]{JON-Mus2012}). Using \cite[Proposition 6.3]{BGZ24}, we have that the Tian function
    \[t \longmapsto \lct(\mathfrak{q}, t\cdot \mathfrak{a}; \mathfrak{c}_{\bullet}^v)\]
    is linear on $[0,+\infty)$ for every nonzero ideal $\mathfrak{a}$ on $X$. Set $\mathfrak{a}=\mathfrak{a}_1\cdot\ldots\cdot \mathfrak{a}_m$. By \cite[Lemma 6.4]{BGZ24}, there exists $w\in\ZVal_X$ such that $w\ge v$ and $w(\mathfrak{a})=v(\mathfrak{a})$. Since every $v(\mathfrak{a}_i)>0$, we get $w(\mathfrak{a}_i)=v(\mathfrak{a}_i)$ for $i=1,\ldots,m$. If $v\in \QM_X\cap \Val_{X,\xi}$ and $\xi$ is a closed point, since $w\ge v$, we also have $w\in\Val_{X,\xi}$.
\end{proof}

\subsection{Two dimensional case via valuative tree}
Now we turn to the $2$ dimensional case, where $R=\mathbb{C} [[x,y]]$. In this subsection, the valuations can take values in $\overline{\mathbb{R}_{\ge 0}}\coloneqq \mathbb{R}_{\ge 0}\cup\{+\infty\}$ on $R^{\times}$. In addition, we only consider the \emph{normalized centered} valuations on $R = \mathbb{C} [[x,y]]$, i.e. the valuations $v$ on $R$ satisfying $v(\mathfrak{m})=1$, where $\mathfrak{m}=(x,y)$ is the maximal ideal of $R$. We recall the valuative tree theory which can be referred to \cite{FJ04, FJ05, FJ05b}.

Denote by $\V$ the space of normalized centered valuations on $R$. Then $\V$ admits a natural \emph{partial ordering} $\le$ given by $v_1\le v_2$ if and only if $v_1(f)\le v_2(f)$ for all $f\in R$; and a natural \emph{(weak) topology}, which is the weakest topology such that the functional $\varphi_f\colon \V \to \overline{\mathbb{R}_{\ge 0}}$ given by $\varphi_f(v)\coloneqq v(f)$ is continuous for all $f\in R$. 

Let $m(f)$ be the \emph{multiplicity} of any $f\in R$:
\[m(f)\coloneqq\max\{k\ge 0\colon f\in\mathfrak{m}^k\}.\]
Then the \emph{multiplicity valuation} $v_{\mathfrak{m}}\in\V$, where $v_{\mathfrak{m}}(f)\coloneqq m(f)$.

Let $C$ be an irreducible (possible formal) curve across the origin $o$ in $\mathbb{C}^2$, and define
\[v_C(f)\coloneqq \frac{C\cdot (f=0)}{m(C)},\]
where ``$\cdot$" denotes the intersection multiplicity and $m$ denotes the multiplicity of the curve. Such a valuation $v_C$ is called a \emph{curve valuation} associated to $C$. If $C$ is defined by $\phi\in \mathfrak{m}$, write $v_C=v_{\phi}$. Note that $v_{\phi}(\phi)=+\infty$.

For every curve valuation $v_C$ associated to an irreducible curve $C$ and every $t\in [1,+\infty)$, define
\[v_{C,t}(f)\coloneqq \min\left\{v_{D}(f)\colon v_{D} \text{ curve valuation with } \frac{m(C)m(D)}{C\cdot D}\le t^{-1}\right\}, \quad \forall\, f\in\mathfrak{m}.\]
The set of \emph{quasimonomial valuations} on $R$, denoted by $\V_{\qm}$, consists of valuations in $\V$ which can be written as some $v_{C,t}$. We have (cf. \cite[Proposition 3.55]{FJ04} and \cite[Section 1.2.3]{FJ05})
\begin{equation}\label{eq-qm.equal}
    v_{C,s}=v_{D,t} \ \Longleftrightarrow \ s=t\le\frac{C\cdot D}{m(C)m(D)}.
\end{equation}
If $C$ is defined by $f$, write $v_{C,t}=v_{f,t}$.

It was shown in \cite{FJ04} that the valuation space $\V$ is a complete tree (called the \emph{valuative tree}; see \Cref{figure-valuative.tree}) rooted at $v_{\mathfrak{m}}$, and parameterized by the \emph{skewness}: For every $v\in \V $, the skewness of $v$ is defined by
\[\alpha(v)\coloneqq\sup\left\{\frac{v(f)}{m(f)}\colon f\in\mathfrak{m}\right\}\in [1,+\infty].\]
In particular, for a quasimonomial valuation $v_{C,t}\in\V_{\qm}$. we have $\alpha(v_{C,t})=t$. The skewness function $\alpha\colon \V \to [1,+\infty]$ is strictly increasing and restricts to a bijection on any segment $[v_{\mathfrak{m}}, v]\coloneqq \{w\in \V\colon v_{\mathfrak{m}}\le w\le v\}$ in $\V$ onto its image $[1,\alpha(v)]$.

Under this parameterization, the valuative tree ends at the elements of $\V\setminus\V_{\qm}$ which are either curve valuations or \emph{infinitely singular valuations} (cf. \cite[Theorem 3.26]{FJ04}).

\begin{figure}[H]
\centering
\begin{tikzpicture}[
    point/.style={draw=none, fill=none},
    linenode/.style={
        draw,          
        line width=0.8pt, 
        inner sep=0pt,  
        minimum height=#1,  
        minimum width=0pt, 
        outer sep=0pt,     
        shape=rectangle,   
        fill=black       
    },
    linenode/.default=6pt
]

\draw[->] (0,0) -- (11,0) node[below] {$\alpha$};
\draw[->] (0,0) -- (0,6) node[above] {$\mathcal{C}$};

\node[linenode, label=below:{$3/2$}] at (2,0) {};
\node[linenode, label=below:{$5/3$}] at (5,0) {};
\node[linenode, label=below:{$2$}] at (6.5,0) {};
\node[linenode, label=below:{$+\infty$}] at (9,0) {};

\node[point, label=left:{$v_{\mathfrak{m}}$}] at (0,3) {};

\node[point, label=right:{$y=0$}] at (9,0.5) {};
\node[point, label=right:{$y=x^2$}] at (9,2) {};
\node[point, label=right:{$(y^2-x^3)^3=x^{10}$}] at (9,3) {};
\node[point, label=right:{$y^2=x^3$}] at (9,5) {};
\node[point, label=right:{$x=0$}] at (9,6) {};

\draw[thick] (0,3) -- (9,0.5);
\draw[thick] (0,3) -- (9,6);
\draw[thick] (2,22/9) -- (9,5);
\draw[thick] (5,223/63) -- (9,3);
\draw[thick] (6.5,43/36) -- (9,2);
\end{tikzpicture}
\caption{The valuative tree; see \cite[p.5]{FJ05}.}
\label{figure-valuative.tree}
\end{figure}

The tree structure on $\V$ induces an \emph{infimum} $\wedge_i v_i\in\V$ for any collection $(v_i)_{i\in I}$ of valuations in $\V$ (see \cite[Corollary 3.15]{FJ04}):
\[\big(\wedge_i v_i\big)(\phi)\coloneqq\inf_{i} v_i(\phi), \quad \forall\, \phi\in\mathfrak{m} \ \text{irreducible},\]
and $v(f g)=v(f)+v(g)$ for all $f,g\in\mathfrak{m}$. We remark that this fails in higher dimensions, i.e. the valuative tree structure does not exist for $\dim\ge 3$; see \cite[p.54]{FJ04}.

\begin{Lemma}[{see \cite[Proposition 3.25 \& Corollary 3.30]{FJ04}}]\label{lem-v.equal.alpha}
    If $v\in\V$ and $\phi\in\mathfrak{m}$ is irreducible, then
    \[v(\phi)=\alpha(v\wedge v_{\phi})m(\phi).\]
    Moreover, if $v(\phi)$ is irrational, then $v=v_{\phi,t}$
    with $t=v(\phi)/m(\phi)$.
\end{Lemma}

For every $v\in\V$, since $v(\mathfrak{m})=1$, we must have 
\[v(f)\ge v\big(\mathfrak{m}^{v_{\mathfrak{m}}(f)}\big)=v_{\mathfrak{m}}(f)=m(f), \quad \forall\, f\in\mathfrak{m}.\]

We consider the following interpolation problem on $R = \mathbb{C} [[x,y]]$.

\begin{Problem}\label{prob-val.inter.2.dim}
Let $k\in\mathbb{Z}_{>0}$. Given $k$ nonzero distinct \underline{irreducible} elements $f_1, \ldots,$ $f_k \in \mathfrak{m}$, and $k$ finite rational numbers $b_1,\ldots, b_k\in\mathbb{Q}_{\ge 1}$ satisfying $b_j\ge m(f_j)$ for $j=1,\ldots, k$, find the necessary and sufficient condition such that there exists a valuation $v\in \V$ which satisfies $v(f_j)=b_j$ for $j=1,\ldots,k$.
\end{Problem}

Lemma \ref{lem-v.equal.alpha} shows the interpolation problem is easier when the interpolated number set includes irrational numbers, hence this case will not be considered in this section.

For $1\le j\le k$, we denote
\[B_j\coloneqq \frac{b_j}{m(f_j)}\ge 1,\]
and
\[\widetilde{B}\coloneqq \max_{1\le j\le k} B_j.\]
Divide $\{1,2,\ldots,k\}$ into two parts:
\[I_1\coloneqq \big\{j\colon B_j=\widetilde{B}\big\}, \quad I_2\coloneqq \big\{j\colon B_j<\widetilde{B}\big\}.\]

\begin{Proposition}\label{prop-val.inter.2.dim}
    There exists a valuation $v\in \V$ such that $v(f_j)=b_j$ for $j=1,\ldots,k$ if and only if the following statements hold:
    \begin{enumerate}
        \item $\alpha\big(\wedge_{i\in I_1} v_{j}\big)\ge \widetilde{B}$; and
        \item $\alpha(v_i \wedge v_j)=B_j$ for $i\in I_1$ and $j\in I_2$.
    \end{enumerate}
    Here $v_j=v_{f_j}$ is the curve valuation associated to $f_j$ for $1\le j\le k$.
\end{Proposition}

\begin{proof}
    Assume that $v\in\V$ satisfies $v(f_j)=b_j$ for $1\le j\le k$. Note that by Lemma \ref{lem-v.equal.alpha}, we have
    \[v(f_j)=\alpha(v\wedge v_j)m(f_j),\]
    which implies that the interpolation problem is equivalent to solving the following equation system for $v$:
    \begin{equation*}
        \alpha(v\wedge v_j)=\frac{v(f_j)}{m(f_j)}=\frac{b_j}{m(f_j)}=B_j, \quad j=1,\ldots,k.
    \end{equation*}
    
    We demonstrate the following simple fact, which will be used repeatedly.
    
    \begin{Lemma}\label{lem-simple.lem}
        Let $w_1,\ldots, w_r\in \V$. If there exists $w\in\V$ satisfying $w_j\le w$ for $1\le j\le r$, then
        \[\alpha\Big(\wedge_{1\le j\le r} w_j\Big)=\min_{1\le j\le r} \alpha(w_j).\]
    \end{Lemma}
    
    \begin{proof}
        Since $w_j\le w$ for $1\le j\le r$, all the valuations $w_j$ lie on the segment $[v_{\mathfrak{m}}, w]$, so there is a full order among these valuations and parameterized by the skewness. Hence the lemma holds.
    \end{proof}
    
    Assume $I_2=\emptyset$ first.
    
    If $|I_1|=1$, i.e. $k=1$, then the ``only if" part is obvious since $\alpha\big(\wedge_{i\in I_1} v_{j}\big)=\alpha(v_1)=+\infty$. For the ``if" part, note that we can take the valuation $v=v_{f_{1},t}\in\V_{\qm}$ with $t=B_1$ which satisfies $\alpha(v)=\alpha(v_{f_1,t})=t=B_1$, and thus $v(f_1)=b_1$.
    
    When $I_2=\emptyset$ and $|I_1|>1$, by Lemma \ref{lem-simple.lem},
    \[\alpha\big(v\wedge (\wedge_{j\in I_1}v_j)\big)=\alpha\big(\wedge_{j\in I_1}(v\wedge v_j)\big)=\min_{j\in I_1}\alpha(v\wedge v_j)=\widetilde{B},\]
    as every $v\wedge v_j\le v$. Then we have $\alpha\big(\wedge_{i\in I_1} v_{j}\big)\ge \widetilde{B}$, which gives the ``only if" part. For the ``if" part, if $\alpha\big(\wedge_{i\in I_1} v_{j}\big)\ge \widetilde{B}$, then by (\ref{eq-qm.equal}) and Lemma \ref{lem-v.equal.alpha},
    \begin{equation}\label{eq-vfj.bj.1}
        v_{f_i, \widetilde{B}} = v_{f_j, \widetilde{B}}\eqqcolon v, \quad \forall\,i,j\in I_1=\{1,\ldots,k\},
    \end{equation}
    and the above valuation $v\in\V_{\qm}$ satisfies
    \begin{equation}\label{eq-vfj.bj.2}
        v(f_j)=m(f_j)\alpha(v\wedge v_j)=m(f_j)\alpha(v)=m(f_j)\widetilde{B}=b_j, \quad j=1,\ldots,k.
    \end{equation}
    Consequently, the proposition holds when $I_2=\emptyset$.
    
    Now we assume $I_2\neq \emptyset$.
    
    First, we prove the ``only if" part. By the $I_2=\emptyset$ case, we only need to prove $\alpha(v_i\wedge v_j)=B_j$ for $i\in I_1$ and $j\in I_2$. Now if $v\in\V$ is the valuation solving the interpolation problem, again by Lemma \ref{lem-simple.lem}, we have
    \begin{equation*}
    \begin{aligned}
        \alpha\left(v\wedge v_i\wedge v_j\right)&=\alpha\Big((v\wedge v_i)\wedge (v \wedge v_j)\Big)\\
        &=\min\big\{\alpha(v\wedge v_i), \alpha(v\wedge v_j)\big\}=\min\big\{\widetilde{B}, B_j\big\}=B_j,
    \end{aligned}
    \end{equation*}
    and
    \begin{equation*}
    \begin{aligned}
        \alpha\left(v\wedge v_i\wedge v_j\right)&=\alpha\Big((v\wedge v_i)\wedge (v_i \wedge v_j)\Big)\\
        &=\min\big\{\alpha(v\wedge v_i), \alpha(v_i\wedge v_j)\big\}=\min\big\{\widetilde{B}, \alpha(v_i\wedge v_j)\big\},
    \end{aligned}
    \end{equation*}
    which implies $\alpha(v_i\wedge v_j)=B_j$, $\forall\,i\in I_1$, $j\in I_2$. Therefore, the ``only if" part is proved.
    
    Next, we prove the ``if" part. Just as before, we take the valuation $v\in\V_{\qm}$ given by
    \[v_{f_i, \widetilde{B}} = v_{f_j, \widetilde{B}}\eqqcolon v, \quad \forall\,i,j\in I_1.\]
    Similarly to (\ref{eq-vfj.bj.1}) and (\ref{eq-vfj.bj.2}), we have $v(f_i)=b_i$ for all $i\in I_1$. In addition, as $\alpha(v_i\wedge v_j)=B_j$ for $i\in I_1$ and $j\in I_2$ by assumption, we also have (by taking any $i\in I_1$)
    \begin{equation*}
        \begin{aligned}
            \alpha(v\wedge v_j)=\alpha(v_{f_i,\widetilde{B}}&\wedge v_j)=\alpha\big(v_{f_i,\widetilde{B}}\wedge (v_i \wedge v_j)\big)\\
            &=\min\{\alpha(v_{f_i,\widetilde{B}}), \alpha(v_i\wedge v_j)\}=\min\{\widetilde{B}, B_j\}=B_j,
        \end{aligned}
    \end{equation*}
    for all $j\in I_2$. Therefore, the valuation $v\in\V$ satisfies $v(f_j)=b_j$ for all $j\in\{1,\ldots,k\}$.
    
    The proof is complete.
\end{proof}

By the proof, we can see that if the valuative interpolation problem \ref{prob-val.inter.2.dim} is solvable, then the solution $v\in\V$ can be selected from $\V_{\qm}$, and in fact the solution can be selected from $\V_{\mathrm{div}}$ (which is the space of \emph{divisorial valuations}, consisting of valuations in $\V_{\qm}$ with rational skewness; see \cite[Theorem 3.26]{FJ05}) for this case that all the interpolated numbers are rational.

Moreover, we can see that if the valuative interpolation problem (\Cref{prob-val.inter.2.dim}) is solvable, then it has a unique minimal solution $v_{\min} \in \mathcal{V}$ (in fact, $v_{\min} = v_{f_i, \widetilde{B}}$ for all $i \in I_1$). Any other solution $v \in \mathcal{V}$ satisfies $v > v_{\min}$, and must represent a tangent vector distinct from those represented by all $v_j$ (for $1 \le j \le k$), i.e., the segment $(v_{\min}, v]$ does not intersect any of the segments $(v_{\min}, v_j]$.

\begin{Remark}
Since
\begin{equation*}
    \begin{aligned}
        \alpha(\wedge_{i} v_i)=\min_{i\neq j}\alpha(v_i\wedge v_j)&=\min_{i\neq j}\frac{(f_i=0)\cdot (f_j=0)}{m(f_i)m(f_j)}\\
        &=\min_{i\neq j}\frac{\dim_{\mathbb{C}}\mathbb{C}[[x,y]]/(f_i,f_j)}{m(f_i)m(f_j)},
    \end{aligned}
\end{equation*}
we can translate Proposition \ref{prop-val.inter.2.dim} into: There exists $v\in \V$ such that $v(f_j)=b_j$ for $j=1,\ldots,k$ if and only if 
\[\min_{\substack{i,j\in I_1 \\ i\neq j}}\frac{\dim_{\mathbb{C}}\mathbb{C}[[x,y]]/(f_i,f_j)}{m(f_i)m(f_j)}\ge \widetilde{B}\]
and
\[\frac{\dim_{\mathbb{C}}\mathbb{C}[[x,y]]/(f_i,f_j)}{m(f_i)m(f_j)}=B_j, \quad \forall\, i\in I_1, \ j\in I_2.\]

In particular, we have the following corollary. 
\end{Remark}

\begin{Corollary}\label{cor-finite.increasing.exist}
Let $k\in\mathbb{Z}_{>1}$. Given $k$ nonzero distinct irreducible elements $f_1, \ldots$, $f_k \in \mathfrak{m}$, and $k$ rational numbers $b_1,\ldots, b_k\in\mathbb{Q}_{\ge 1}$. If \[\frac{b_k}{m(f_k)}>\max_{1\le j\le k-1}\frac{b_j}{m(f_j)},\]
then there exists a valuation $v\in \V$ such that $v(f_i)=b_i$ for $i=1,\ldots,k$ if and only if 
\[\dim_{\mathbb{C}}\mathbb{C}[[x,y]]/(f_k,f_j)=m(f_k)b_j, \quad j=1,\ldots,k-1.\]
\end{Corollary}

\begin{proof}
    For this case, we have $I_1=\{k\}$, and this corollary follows from the above remark immediately.
\end{proof}

We can also consider the infinite valuative interpolation problem in $\mathbb{C}[[x,y]]$. 

\begin{Proposition}\label{prop-exist.2dim.infinite}
    Let $f_j\in\mathfrak{m}$ be a sequence nonzero distinct irreducible polynomials, and $b_j\in\mathbb{Q}_{\ge 1}$ be a sequence of rational numbers, $j\in\mathbb{Z}_+$. Suppose $b_j/m(f_j)$ is strictly increasing in $j$. Then there exists a valuation $v\in\V$ satisfying
    \begin{equation}\label{eq-val.inter.infinite}
        v(f_j)=b_j, \quad \forall\,j\in\mathbb{Z}_+,
    \end{equation}
    if and only if
    \begin{equation}\label{eq-val.inter.infinite.condition}
        \dim_{\mathbb{C}}\mathbb{C}[[x,y]]/(f_i,f_j)=m(f_i)b_j, \quad \forall\,i,j\in\mathbb{Z}_+ \ \text{with} \ i>j.
    \end{equation}
    
    Moreover, the following statements hold:
    \begin{enumerate}
        \item There exists $v\in\V_{\qm}$ satisfying (\ref{eq-val.inter.infinite}) if and only if there exists $v\in\V_{\mathrm{div}}$ satisfying (\ref{eq-val.inter.infinite}), and if such $v\in\V_{\qm}$ exists then the sequence $\big(b_j/m(f_j)\big)_{j\ge 1}$ must be bounded.
        \item If the sequence $(b_j)_{j\ge 1}$ has unbounded denominators, then every $v\in\V$ satisfying (\ref{eq-val.inter.infinite}) must be infinitely singular.
    \end{enumerate}
\end{Proposition}

\begin{proof}
    If there exists $v\in\V$ satisfying (\ref{eq-val.inter.infinite}), then $v(f_j)=b_j$ for $j=1,\ldots,i$, and thus by Corollary \ref{cor-finite.increasing.exist}, (\ref{eq-val.inter.infinite.condition}) must hold since $b_j/m(f_j)$ is strictly increasing. Conversely, if (\ref{eq-val.inter.infinite.condition}) holds, then for every $k\in\mathbb{Z}_+$, there exists some $v_k\in \V$ such that
    \[v(f_j)=b_j, \quad j=1,\ldots,k.\]
    As $\V$ is complete (cf. \cite[Theorem 3.14 and Proposition 7.13]{FJ04}), we can extract a subsequence of $(v_k)_{k\ge 1}$ which is (weakly) convergent to some $v\in \V$. Clearly we have $v(f_j)=b_j$ for all $j\ge 1$.
    
    If there exists $v\in\V_{\qm}$ satisfying (\ref{eq-val.inter.infinite}), we may write $v=v_{C,t}$ for some irreducible curve $C$ and $t\in [1,+\infty)$. Note that $v_{C,t}\in \V_{\mathrm{div}}$ if and only if $t$ is rational. Otherwise, if $v=v_{C,t}$ for some irrational $t$, then for every $j$,
    \begin{equation*}
        \begin{aligned}
            b_j=v_{C,t}(f_j)&=m(f_j)\alpha(v_{C,t}\wedge v_{f_j})\\
            &=m(f_j)\min\{t, \alpha(v_{C} \wedge v_{f_j})\}
            =m(f_j)\alpha(v_{C} \wedge v_{f_j})=v_C(f_j),
        \end{aligned}
    \end{equation*}
    since $b_j$ is rational. It follows that the curve valuation $v_C$ satisfies (\ref{eq-val.inter.infinite}), and $v_{C,s}\in\V_{\mathrm{div}}$ also satisfies (\ref{eq-val.inter.infinite}) for all rational number $s$ with $s>t$.
    
    If $v=v_{C,t}\in\V_{\qm}$ satisfies (\ref{eq-val.inter.infinite}), we have
    \[\frac{b_j}{m(f_j)}=\min\{t, \alpha(v_C\wedge v_{f_j})\}\le t, \quad \forall\, j\ge 1.\]
    Thus, $\big(b_j/m(f_j)\big)_{j\ge 1}$ is bounded.
    
    If the sequence $(b_j)_{j\ge 1}$ has unbounded denominators and $v\in \V$  satisfies (\ref{eq-val.inter.infinite}), then the value semigroup $v(R^*)$ is not finitely generated over $\mathbb{N}$. However, \cite[Proposition 3.52]{FJ04} shows that any valuation in $\V$ with a finite approximating sequence (i.e. of finite multiplicity) must admit a finitely generated value semigroup over $\mathbb{N}$. We conclude that the multiplicity of the valuation $m(v)=+\infty$, so $v$ is infinitely singular (\cite[Proposition 3.37]{FJ04}).
\end{proof}

\begin{Example}
    Let $k\in\mathbb{Z}_+\cup\{+\infty\}$ and
    \[f_j(x,y)=y-x-2!x^2-\cdots-j!x^j\in \mathbb{C}[[x,y]], \quad 1\le j\le k.\]
    Then for any increasing sequence of rational numbers $1\le b_1<b_2<\ldots<b_k$, there exists a valuation $v\in\V$ satisfying $v(f_j)=b_j$ for $1\le j\le k$ if and only if
    \begin{enumerate}
        \item $b_j=j+1$ for $j=1,\ldots,k-1$ and $b_k>k$, when $k<+\infty$;
        \item $b_j=j+1$ for all $j\ge 1$, when $k=+\infty$,
    \end{enumerate}
    due to Corollary \ref{cor-finite.increasing.exist}, Proposition \ref{prop-exist.2dim.infinite}, and
    \[\dim_{\mathbb{C}}\mathbb{C}[[x,y]]/(f_i,f_j)=\min\{i,j\}+1, \quad i\neq j.\]
    In fact, for $k=+\infty$, the valuation $v\in\V$ given by
    \[v\big(g(x,y)\big)=\mathrm{ord}_t\, g\Big(t\,, \,\sum_{m\ge 1} m!t^m\Big), \quad \forall\,g(x,y)\in\mathbb{C}[[x,y]],\]
    satisfies $v(f_j)=j+1$ for all $j\ge 1$. This is a curve valuation associated with a formal curve, which is not a quasimonomial valuation.
\end{Example}


\end{document}